\documentclass[onefignum,onetabnum]{siamart171218}
\pdfoutput=1

\usepackage{amssymb,amsmath,amsfonts,stmaryrd} 
\usepackage{mathtools} 
\usepackage{graphicx}
\usepackage{color}
\usepackage{tabularx}
\usepackage{tikz} 
\usepackage{enumitem}
\usepackage[labelfont=bf]{caption}
\usepackage{changepage}
\usepackage{afterpage} 
\usepackage{gensymb} 
\usepackage{placeins} 
\usepackage{enumitem}
\usepackage{siunitx} 
\usepackage{soul} 
\usepackage[export]{adjustbox} 
\usepackage{tabularx}
\usepackage{wrapfig}
\usepackage[export]{adjustbox}
\usepackage{mwe}
\usepackage{dashrule}
\usepackage{graphbox}
\usepackage{comment}
\usepackage{gensymb}
\usepackage{lscape} 

\usepackage[belowskip=-5pt]{subcaption}

\usepackage{caption}
\usepackage{lipsum}
\usepackage{amsfonts}
\usepackage{graphicx}
\usepackage{epstopdf}
\usepackage{algorithmic}
\ifpdf
  \DeclareGraphicsExtensions{.eps,.pdf,.png,.jpg}
\else
  \DeclareGraphicsExtensions{.eps}
\fi

\newsiamremark{remark}{Remark}
\newsiamremark{hypothesis}{Hypothesis}
\crefname{hypothesis}{Hypothesis}{Hypotheses}
\newsiamthm{claim}{Claim}

\headers{Additive Polynomial Time Integrators, Part I}{T. Buvoli and B. S. Southworth}

\title{Additive Polynomial Time Integrators, Part I: Framework and Fully-implicit-explicit (FIMEX) Collocation Methods \thanks{Submitted to the editors DATE.
\funding{B. S. S. was supported by the Laboratory Directed Research and Development program of Los Alamos National Laboratory as a Nicholas C. Metropolis Fellow and under project number 20220174ER. TB was funded in part by the NSF grant DMS-2012875.}}}

\author{Tommaso Buvoli\thanks{Department of Mathematics, Tulane University, 6823 St. Charles Avenue, New Orleans, Louisiana 70118, USA (\email{tbuvoli@tulane.edu}).}
  \and Ben S. Southworth 
	\thanks{Theoretical Division, Los Alamos National Laboratory, USA (\email{southworth@lanl.gov}).}
}

\usepackage{amsopn}

\ifpdf
\hypersetup{
  pdftitle={Additive Polynomial Time Integrators, Part I: Framework and Fully-Implicit-Explicit (FIMEX) Collocation Methods},
  pdfauthor={T. Buvoli and B. S. Southworth}
}
\fi

\newcommand{\ODEd}{ODE dataset}
\newcommand{\ODED}{ODE Dataset}
\newcommand{\ODEds}{ODE datasets}

\newcommand{\ODEp}{ODE polynomial}
\newcommand{\ODEP}{ODE Polynomial}

\newcommand{\ODEps}{ODE polynomials}

\definecolor{plot_black}{RGB}{10, 10, 10}
\definecolor{plot_red}{RGB}{192, 57, 43}
\definecolor{plot_orange}{RGB}{230, 126, 34}
\definecolor{plot_yellow}{RGB}{241, 196, 15}
\definecolor{plot_green}{RGB}{39, 174, 96}
\definecolor{plot_blue}{RGB}{41, 128, 185}
\definecolor{plot_violet}{RGB}{155, 89, 182}
\definecolor{plot_grey}{RGB}{149, 165, 166}

\begin{document}

\maketitle

\begin{abstract}
In this paper we generalize the polynomial time integration framework to additively partitioned initial value problems. The framework we present is general and enables the construction of many new families of additive integrators with arbitrary order-of-accuracy and varying degree of implicitness. In this first work, we focus on a new class of implicit-explicit polynomial block methods that are based on fully-implicit Runge-Kutta methods with Radau nodes, and possess high stage order. We show that the new fully-implicit-explicit (FIMEX) integrators have improved stability compared to existing IMEX Runge-Kutta methods, while also being more computationally efficient due to recent developments in preconditioning techniques for solving the associated systems of nonlinear equations. For PDEs on periodic domains where the implicit component is trivial to invert, we will show how parallelization of the right-hand-side evaluations can be exploited to obtain significant speedup compared to existing serial IMEX Runge-Kutta methods. For parallel (in space) finite-element discretizations, the new methods can achieve orders of magnitude better accuracy than existing IMEX Runge-Kutta methods, and/or achieve a given accuracy several times times faster in terms of computational runtime.
\end{abstract}

\begin{keywords}
  Additive integrators, Linearly implicit, Implicit-Explicit, Fully-implicit Runge-Kutta, General linear methods
\end{keywords}

\begin{AMS}
	65L04, 65L05, 65L06
\end{AMS}

\section{Introduction}

Many problems in science and engineering can be modeled using high-dimensional systems of ordinary differential equations. These equations typically arise  %
 from the mathematical description of physical phenomena or from the spatial discretization of partial differential equations. %
Solving these systems amounts to integrating an initial value problem (IVP)
\begin{align}
	y'(t) = f(t,y(t)), \quad y(t_0) = y_0.
	\label{eq:model-ode}
\end{align}
In practice, it is common to additively partition the right-hand-side, $f(t,y)$, into $m$ components,
\begin{align}
	f(t,y) = \sum_{k=1}^m f^{\{k\}}(t,y).
	\label{eq:partitioned-rhs}
\end{align}
A simple example is a discretized advection-diffusion-reaction equation where each physical process is represented by a separate term. 

The solution of any additively partitioned system can be numerically approximated using an additive integrator   that treats each component $f^{\{k\}}(t,y)$ differently \cite{cooper1980additive, kennedy2003additive, sandu2015generalized}. This can be particularly efficient for solving multiscale, multiphysics problems where the optimal method differs across the components, or where it is prohibitively expensive to treat the full operator $f(t,y)$ implicitly. A canonical example is a linear advection-diffusion equation, where the diffusion places severe explicit time-step restrictions, but fully implicit solves for advection-diffusion discretizations are significantly more challenging than for pure diffusion. An additive integrator can treat the diffusion implicitly and the advection explicitly, addressing each of these problems. If the right-hand-side consists of two components and an additive integrator treats the first implicitly and the second explicitly, then the integrator is frequently called an implicit-explicit (IMEX) method.

Two closely related classes of integrators are linearly implicit methods and W-methods, that respectively utilize an exact or approximate local Jacobian of $f(t,y)$ at each timestep \cite[IV.7]{wanner1996solving}\cite{calvo2001linearly,akrivis2003linearly,glandon2020linearly}. Given the solution at the $n$th timestep $y_n = y(t_n)$, we can rewrite the system (\ref{eq:model-ode}) as an additively partitioned system with $m=2$,
\begin{align}
	 f^{\{1\}}(t,y) = J_ny, \quad \text{and} \quad f^{\{2\}}(t,y) = f(t,y) - J_{n}y,
\end{align}
where $J_n$ approximates or is equal to the local Jacobian  $\frac{\partial f}{\partial y}(t_n, y_n)$. Any additive integrator that treats $f^{\{1\}}(t,y)$ implicitly and $f^{\{2\}}(t,y)$ explicitly reduces to a linearly implicit method.

In the past three decades, the construction of additive integrators has been an active area of research that has produced a range of methods, including linear multistep methods (LMMs) \cite{ascher1995implicit, dimarco2017implicit, wang2008variable}, Runge-Kutta (RK) methods \cite{ascher1997implicit,kennedy2003additive,Minion2003IMEX,sandu2015generalized,izzo2017highly,kennedy2019higher}, and general linear methods (GLMs) \cite{cardone2015construction,zhang2016high,soleimani2018superconvergent}, including those based on extrapolation \cite{constantinescu2010extrapolated,cardone2014extrapolation}. Each method class has certain benefits and drawbacks. Additive LMMs have a low computational cost per timestep, but high-order methods experience instabilities on equations with limited diffusion. Additive diagonally implicit RK methods possess good stability and allow for simplified adaptive time-stepping, however they are known to suffer from order-reduction on stiff equations \cite{boscarino2007error}. Moreover, high-order RK method derivations that rely on nonlinear order conditions grow increasingly difficult to construct, and alternative approaches must be considered \cite{constantinescu2010extrapolated,Minion2003IMEX}. GLMs with good stability and no order-reduction exist, however to avoid increasing the number of order conditions further, one must typically consider simplified method formulations. 

In this paper we will generalize the recently introduced polynomial time integration framework \cite{buvoli2018polynomial, buvoli2019constructing, buvoli2019esdc, buvoli2021epbm} to include additively partitioned differential equations. We then demonstrate its utility by introducing a new family of Implicit-Explicit (IMEX) integrators where the implicit integrator is a fully-implicit collocation method. Selecting a fully-implicit integrator may seem peculiar since these methods are often considered too slow to be competitive. However, recent developments in block linear and nonlinear solvers make their use in numerical PDEs quite tractable \cite{chen14,pazner17,farrell2020irksome,jiao2020optimal,rana2020new,irk1,irk2}, even outperforming diagonally implicit RK methods in many cases \cite{irk1,irk2}. Despite these developments, it is not possible to derive IMEX RK methods based on the fully implicit RK methods since the explicit stages would be nonlinearly coupled to the fully implicit stages. More generally, constructing high order IMEX RK schemes is often nontrivial, as mentioned in \cite{kennedy2019higher}. In this work, we will show how the additive polynomial framework provides a natural way to develop fully-implicit-explicit (FIMEX) integrators, which are high-order accurate, and allow us to leverage developments in fully implicit solvers in the context of additive integration.

More generally, there are two main advantages of the additive polynomial framework that will be explored in this paper. First, the framework simplifies the construction of high-order additive GLMs that do not suffer from order-reduction on stiff equations. Specifically, the polynomial framework makes extensive use of interpolating polynomials that trivially satisfy nonlinear order conditions and ensure high stage-order. Furthermore, method construction can be done using geometric arguments that are similar to those used to derive spatial finite difference stencils. Second, the framework can be used to derive efficient methods for solving equations that are either naturally split into multiple terms, or where the right-hand-side has been rewritten using an exact or approximate Jacobian. This allows us to simultaneously introduce a range of new high-order additive and semi-implicit integrators.

This paper is organized as follows. In \cref{sec:additive-integrators,sec:polynomial-integrators} we respectively provide short introductions to additive integrators and the polynomial framework. In Section \cref{sec:additive-polynomial-framework} we generalize the polynomial framework for additively partitioned differential equations, and in Section \ref{sec:constructing-imex-pbms} we develop new classes of fully-implicit-explicit methods and study their stability. Lastly, \cref{sec:numerical-experiments} demonstrates the improved accuracy and efficiency of the new integrators for solving partial differential equations.

\section{Additive integrators}
\label{sec:additive-integrators}

Additive integrators are a class of methods for solving the partitioned initial value problem (\ref{eq:model-ode}, \ref{eq:partitioned-rhs}). In this section, we give a short introduction to additive integrators for equations with two partitions ($m=2$), where
	\begin{align}
		y' = f^{\{1\}}(t,y) + f^{\{2\}}(t,y), \quad y(t_0) = y_0.
		\label{eq:model-ode-m2}	
	\end{align}
We assume that the term $f^{\{1\}}(t,y)$ is stiff (i.e. a small stepsize is required for any explicit method when solving $y'=f^{\{1\}}(t,y)$) while the term $f^{\{2\}}(t,y)$ is nonstiff.

 In such a scenario it is desirable to consider integrators that only treat $f^{\{1\}}(t,y)$ implicitly. One approach for deriving additive methods is to integrate (\ref{eq:model-ode-m2}), and then approximate the resulting integrals for each term separately
\begin{align}
	y(t_{n+1}) = y(t_n) + \underbrace{\int_{t_n}^{t_{n+1}} f^{\{1\}}(t,y(t))dt}_{\text{treat implicitly}} + \underbrace{\int_{t_n}^{t_{n+1}} f^{\{2\}}(t,y(t))dt}_{\text{treat explicitly}}.
	\label{eq:imex-euler-derivation-start}
\end{align}
One of the simplest additive integrators can be derived by taking an implicit one-sided approximation for $f^{\{1\}}(t,y)$, and an explicit one-sided approximation for $f^{\{2\}}(t,y)$. This produces the IMEX Euler method 
\begin{align}
	y_{n+1} = y_{n} + h f^{\{1\}}_{n+1} + h f^{\{2\}}_{n},
	\label{eq:imex-euler-generic}	
\end{align}
where the stepsize $h=t_{n+1}-t_n$. Higher-order IMEX-LMM methods \cite{ascher1995implicit, dimarco2017implicit} use higher-order polynomial approximations constructed from previous solution values, while the output of a higher-order IMEX-RK method \cite{ascher1997implicit,kennedy2003additive,sandu2015generalized,kennedy2019higher} is a linear combination of newly computed stage values. IMEX-GLMs \cite{cardone2015construction,zhang2016high,soleimani2018superconvergent} combine both ideas by using previous solution values and new stages.

When using any IMEX method there are a range of choices for $f^{\{1\}}$ and $f^{\{2\}}$ that affect the computational cost and stability of the integrator. For a semi-linear system ${y' = Ly + N(t,y)}$, with a nonstiff nonlinearity, it is natural to let
\begin{align}
	f^{\{1\}}(t,y) &= Ly, ~~ f^{\{2\}}(t,y) = N(t,y).
	\label{eq:partition-semilinear}
\end{align}
For a more general nonlinear system $y'=A(t,y) + B(t,y)$ that naturally splits into two components $A(t,y)$ and $B(t,y)$, several choices for the numerical splitting include:
\begin{enumerate}[leftmargin=*]
	\item Fully implicit in $A(t,y)$:
	\begin{align}
		& f^{\{1\}}(t,y) = A(t,y), \quad
		f^{\{2\}}(t,y) = B(t, y). 
		\label{eq:partition-fully-implicit}
	\end{align}
	\item Linearly implicit in $A(t,y)$:
	\begin{align}
		& f^{\{1\}}(t,y) = \frac{\partial{A}}{\partial{y}}(t_n, y_n)y, \quad
		f^{\{2\}}(t,y) = B(t,y) + A(t,y) - \frac{\partial{A}}{\partial{y}}(t_n,y_n)y. 
		\label{eq:partition-linearly-in-first}
	\end{align}
	\item Linearly implicit in $A(t,y)$ and $B(t,y)$:
	\begin{align}
		\begin{aligned}
		& f^{\{1\}}(t,y) = \left[ \frac{\partial{A}}{\partial{y}}(t_n, y_n) + \frac{\partial{B}}{\partial{y}}(t_n, y_n) \right]y, \\
		& f^{\{2\}}(t,y) = B(t,y) + A(t,y) - \left[\frac{\partial{A}}{\partial{y}}(t_n,y_n) + \frac{\partial{B}}{\partial{y}}(t_n,y_n)\right]y.
		\end{aligned} 
	\end{align}
	\item Linearly implicit in $J_n$: ($J_n$ approximates the full or partial Jacobian at $t=t_n$)
	\begin{align}
		f^{\{1\}}(t,y) &= J_ny, &
		f^{\{2\}}(t,y) &= B(t, y) + A(t, y) - J_ny.
		\label{eq:partition-approximate-jacobian}
	\end{align}
\end{enumerate}
For additional clarity we write the formulas for the IMEX-Euler method \cref{eq:imex-euler-generic} using three of the proposed partitionings:
\begin{align*}
	&\text{partitioning (\ref{eq:partition-semilinear})} && y_{n+1} = (I - h L)^{-1}(y_{n} + hN(t_n, y_n)), \\
	&\text{partitioning (\ref{eq:partition-fully-implicit})} && y_{n+1} = y_{n} + h A(t_{n+1},y_{n+1}) + h B(t_{n},y_{n}), \\
	&\text{partitioning (\ref{eq:partition-approximate-jacobian})} && y_{n+1} = y_{n} + (I - h J_n)^{-1}h(B(t_n, y_n) + A(t_n, y_n)).
\end{align*}
Selecting a fully implicit partitioning generally leads to improved stability but requires a nonlinear solve at each timestep. Conversely, a linearly implicit choice only requires a linear solve at each step, but the method may have inferior stability properties, especially if $J_n$ does not closely approximate the local Jacobian. 

In the sections that follow we will first review the polynomial time integration framework, and then generalize it so that we can construct new high-order additive integrators for (\ref{eq:imex-euler-generic}) using any of the proposed partitionings (\ref{eq:partition-semilinear})-(\ref{eq:partition-approximate-jacobian}).
\section{Polynomial time integrators}
\label{sec:polynomial-integrators}

The polynomial time integration framework \cite{buvoli2018polynomial,buvoli2019constructing} is based on continuous polynomials in time that are constructed by fitting through solution or derivative values. The values may be known (i.e. input values) or unknown (i.e. future stage values or outputs), with the latter leading to implicit equations and ultimately implicit methods. Within the family of classical time integration methods, we can interpret polynomial time integrators  as parametrized general linear methods (GLMs) whose inputs and outputs approximate the solution at a set of scaled nodes $\{z_j\}$. The nodes $z_j$ may be positive, negative, or zero, and the scaling factor is a positive number called the {\em node radius} $r$. Note that the node radius introduces an additional degree of freedom that complements the stepsize $h$. 

\begin{figure}[!ht]
	\centering
	\includegraphics[width=0.5\linewidth]{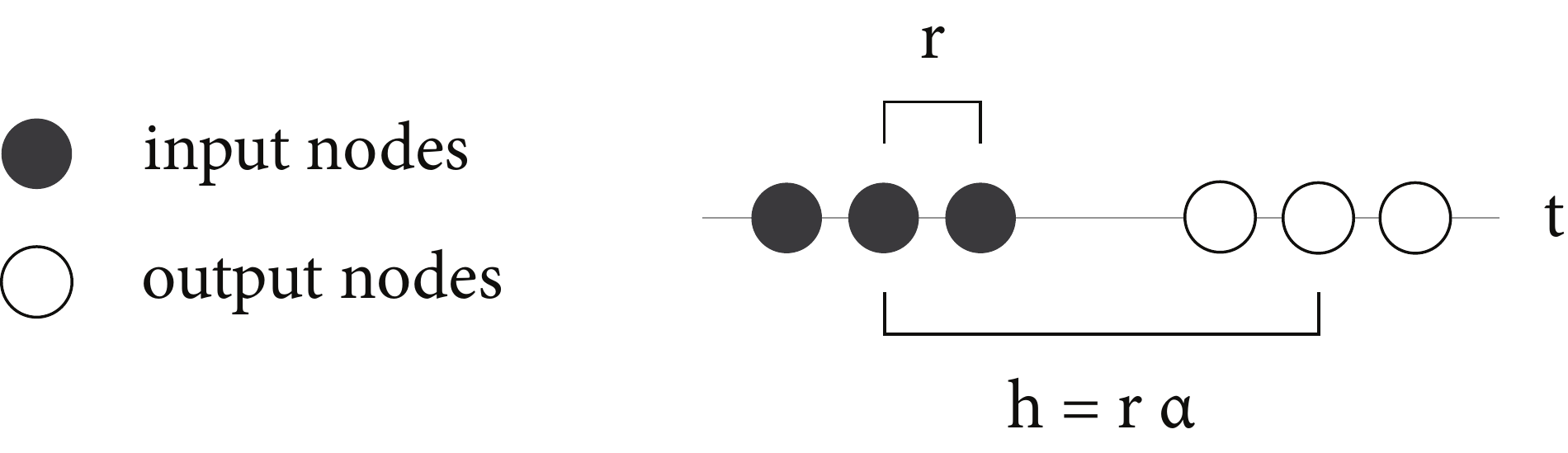}
	\caption{A diagram showing the input and output nodes for a polynomial integrator with three equispaced input and output nodes (e.g. ${z_j} = \{-1,0,1\}$). The node set is scaled by the node radius and the distance between the inputs and outputs is proportional to the extrapolation factor $\alpha$.}
	\label{fig:alpha_visualization}	
\end{figure}

The input and output values of a polynomial integrator, along with the associated derivatives, are represented using the notation
	\begin{center}
		\renewcommand*{\arraystretch}{1.5}
		\begin{tabular}{rl}
			input (solutions):  & \hspace{0.65em} $y_j^{[n]} \approx y\left(t_n + r z_j\right)$  \\
			output (solutions): & $y_j^{[n+1]} \approx y\left(t_n + r z_j + h\right)$ \\
			input derivatives:  & \hspace{0.65em} $f_j^{[n]} \coloneqq f(t_n+rz_j, y_j^{[n]}) \approx y'(t_n+rz_j)$ \\
			output derivatives: & $f_j^{[n+1]} \coloneqq f\left(t_n + r z_j + h, y_j^{[n+1]}\right) \approx y'(t_n+rz_j + h)$. \\[0.5em]
		\end{tabular}			
	\end{center}
Since the polynomials that make up polynomial methods are expressed in local coordinates (\ref{eq:odep_local_coordinates}), it is convenient to also parametrize the stepsize $h$ in terms of the node radius. We therefore let
	\begin{align}
		h = r \alpha	
	\end{align}
where the constant $\alpha$ is called the {\em extrapolation factor}. In Figure \ref{fig:alpha_visualization} we show a visualization of the parameters $r$, $h$, and $\alpha$  for a method with three real-valued, equispaced nodes $\{z_j\}$.

In this work we look at methods with real-valued quadrature nodes $\{z_j\}$ that are scaled relative to the interval $[-1,1]$. This choice is inspired by the fact that nearly all of the theoretical analysis relating to orthogonal polynomials and interpolation is conducted on this interval. However, when implementing a polynomial integrator it is typically simpler to translate the nodes to the interval $[0, 2]$ such that the left-most node at $\tau=0$ is located at $t=t_0$ where the initial condition for \eqref{eq:model-ode} is provided. Translating the nodes only affects temporal locations of the inputs, as determined by \eqref{eq:odep_local_coordinates}, and has no effect on the method coefficients.
	
The name polynomial integrator originates from the fact that all the method coefficients are derived using interpolating polynomials, known as {\em ODE polynomials}. All \ODEps{} approximate the Taylor series of the solution in local coordinates $\tau$, where
\begin{align}
	t(\tau) = r\tau + t_n.	
	\label{eq:odep_local_coordinates}
\end{align}
The general form for an \ODEp{} of degree $g$ with expansion point $b$ is
	\begin{align}
		p(\tau; b) &= \sum_{j=0}^{g}  \frac{a_{j}(b)(\tau - b)^j}{j!}
		\label{eq:solution_ode_polynomial}
	\end{align}
where the constants $\{a_{j}(b)\}$ are called {\em approximate derivatives} since they approximate the derivatives of the solution in local coordinates such that $a_{j}(b) \approx \left. \frac{d^j}{d\tau^j} y(t(\tau)) \right|_{\tau=b} = r^{j}y^{(j)}(t(b))$ (the factor of $r$ originates from the transformation into local coordinates \eqref{eq:odep_local_coordinates}). Each approximate derivative $a_{j}(b)$ is computed by differentiating interpolating polynomials that are constructed using any subset of the method's inputs, outputs, stages, and the corresponding derivatives. The order of the polynomial method is directly related to the degree of its \ODEps{} and the order of the interpolating polynomials that are used to construct the associated approximate derivatives \cite[Sec. 3.6]{buvoli2018polynomial}. As a consequence, the order of accuracy is always bounded below by $\min(g_{min}, \delta - 1)$ where $g_{min}$ is the minimum  degree of the \ODEps{} and $\delta$ is the minimum degree of all the polynomials that determine the approximate derivatives.

A general formulation for the approximate derivatives $a_j(b)$ is described in \cite{buvoli2018polynomial,buvoli2019constructing}; here we will only describe two important sub-families:
	\begin{enumerate}[leftmargin=*]
		\item {\em Adams \ODEps{}} are constructed using two Lagrange interpolating polynomials $L_y(\tau)$ and $L_f(\tau)$ that respectively satisfy solution values or derivative values at the method's input, output, or stage nodes. Specifically, $L_y(\tau) \approx y(t(\tau))$ interpolates at least one solution value, and $L_f(\tau) \approx ry'(t(\tau))$ is a polynomial of degree $g-1$ that interpolates $g$ derivative values. The approximate derivatives are then
			\begin{align}
				a_0(b) = L_y(b), 
				\quad \text{and} \quad
				a_{j}(b) = \frac{d^{j-1} L_f}{d\tau^{j-1}}(b),
				\quad j=1,\ldots,g.
				\label{eq:adams_approximate_derivatives}
			\end{align}					
		By substituting \eqref{eq:adams_approximate_derivatives} into \eqref{eq:solution_ode_polynomial} and noting that  $p'(\tau;b) = L_f(\tau)$, we can express an Adams \ODEp{} in the equivalent integral form
		\begin{align}
			p(\tau;b) = L_y(b) + \int_{b}^\tau L_f(\xi) d\xi.
		\label{eq:adams_poly_integral_form}
		\end{align}
		We now see that an Adams \ODEp{} approximates the integral equation of an initial value problem where the expansion point $b$ is the location of the initial condition in local coordinates.
		\vspace{0.5em}
		\item {\em BDF \ODEps{}} are constructed using a polynomial ${H_y(\tau) \approx y(t(\tau))}$ of degree $g$ that interpolates $g$ solution values and whose derivative $H_y'(\tau)$ interpolates a single derivative value. The approximate derivatives are given by 
			\begin{align}
				a_j(b) = \frac{d^{j} H_y}{d\tau^j}(b) \quad \implies \quad p(\tau;b) = H_y(\tau) ~ \forall b.
				\label{eq:bdf_approximate_derivatives}	
			\end{align}
			Note that BDF polynomials do not depend on the expansion point $b$.			
	\end{enumerate}

When presenting polynomial methods, it is convenient to introduce a set containing all the data values that can be used to construct the interpolating polynomials that determine the approximate derivatives $a_j(b)$ in the \ODEp{} \eqref{eq:solution_ode_polynomial}. This set is called the {\em \ODEd{}}, and in the case of classical polynomial methods, it simply contains the method's inputs, outputs, stages, and their derivatives, along with the corresponding temporal nodes. An \ODEd{} of size $w$ is denoted as
	\begin{align}
		D(r, t_n) = \left\{ \left(\tau_j,~ y_j,~ r f_j \right) \right\}_{j=1}^w
		\quad \text{where} \quad
		y_j \approx y(t(\tau_j)), ~ f_j = f(t(\tau_j), y_j).
		\label{eq:ode_ds}
	\end{align}
	
Using all the previous definitions, the formula for any polynomial method with $s$ stages and $q$ outputs can be written compactly as
\begin{align}
	\begin{aligned}
		Y_i &= p_j(c_j(\alpha), b_j(\alpha)) & j &= 1, \ldots, s,\\
		y^{[n+1]}_j &= p_{j+s}(z_j + \alpha; \hspace{0.125em} b_{j+s}(\alpha)) & j &= 1, \ldots, q,
	\end{aligned}
	\label{eq:polynomial_glm}
\end{align}
where $Y_i$ denote stage values, $c_j(\alpha)$ are stage nodes in local coordinates, and $p_j(\tau; b)$ are \ODEps{} constructed from an \ODEd{} of size $w=2q+s$ that contains the methods inputs, stage values, and outputs. The associated temporal nodes of the method's \ODEds{} \eqref{eq:ode_ds} are respectively
		\begin{align}
			\tau_j = 
			\left\{
			\begin{array}{lll}
 				z_j &  1 \le j \le q, \\
 				c_{j-q}(\alpha) & q+1 \le j \le q + s, \\
 				z_{j-q-s} + \alpha & q + s < j \le 2q + s.
 			\end{array}
 			\right.
		\end{align}
	We remark that our definition of stages differs from the standard convention used for general linear methods and Runge-Kutta methods. Specifically, if one recasts the method \cref{eq:polynomial_glm} as a GLM then the outputs $y^{[n+1]}_j$ also double as additional stage values; in other words, \cref{eq:polynomial_glm} will be a GLM with $s+q$ stages. By avoiding the standard convention we allow for a more compact method definition and avoid method coefficient duplication.

Using the general formulation \eqref{eq:polynomial_glm}, it is possible to derive many different families of polynomial integrators. One example is \emph{polynomial block methods} (PBMs) from \cite{buvoli2019constructing} which are characterized by $s=0$ and will be the primary focus of this paper due to their simpler structure. Additional examples include well-known time integrators like backward difference formulas (BDF), Adams-Moulton methods, and collocation methods, which can all be expressed as PBMs with a fixed $\alpha$.  In the following subsection we show how to write fully implicit collocation methods in the polynomial framework. The resulting formulation will be used again in Section \ref{sec:constructing-imex-pbms} to construct new additive polynomial integrators.

\subsection{Collocation methods and Radau IIA}
\label{sec:polynomial-radau}

Collocation methods \cite[II.7]{hairer1993nonstiff} are time integrators based on polynomial quadrature that can also be expressed as fully-implicit Runge-Kutta methods. Well-known examples include the A-stable Gauss methods \cite[Sec. 342]{butcher2008numerical} and the L-stable Radau IIA methods \cite{hairer1999stiff} that respectively achieve orders of $2\sigma$ and $2\sigma-1$, where $\sigma$ is the number of stages. An additional benefit of these methods is that they both satisfy B-stability \cite{wanner1996solving,wanner1976}.

Suppose that we seek an approximate solution of \eqref{eq:model-ode} at $t=t_n+h$ given an initial condition $y_n = y(t_n)$. To derive a collocation method, we can approximate the solution $y(t)$ using a polynomial $p_y(t)$ that satisfies the initial condition at $t=t_n$ and the differential equation at a set of $m$ collocation points $t_{n,j}$ such that  $p'_y(t_{n,j}) = f(t_{n,j}, p_y(t_{n,j}))$. These constraints lead to the fully-implicit nonlinear system
	\begin{align}
		p_y(t_{n,j}) = y_n + \int_{t_n}^{t_{n,j}} \underbrace{\sum_{j=1}^\sigma \ell_j(t) f(t_{n,j},p_y(t_{n,j}))}_{p_f(t)}dt, \quad j=1\ldots,m
		\label{eq:nonlinear_collocation_system}		
	\end{align}
	where $\ell_j(t) = \prod_{k\ne j} (t-t_{n,k})/(t_{n,j} - t_{n,k})$ is the $j$th Lagrange basis polynomial and $p_f(t)$ is a Lagrange interpolating polynomial for the solution derivative $y'(t)$. If we compare the right-hand-side of \eqref{eq:nonlinear_collocation_system} to the Adams \ODEp{} integral formulation \eqref{eq:adams_poly_integral_form}, we see that $p_y(t(\tau))$ is equivalent to an Adams \ODEp{} $p(\tau; b)$ with $b = \tau(0)=t_n$, $L_y(\tau) = y_n$, and $L_f(\tau) = rp_f(t(\tau))$. 
	Therefore, by appropriately defining an \ODEd{} and node set $\{z_j\}$ we can express any collocation method \eqref{eq:nonlinear_collocation_system} as a one-step polynomial integrator (i.e. \eqref{eq:polynomial_glm} with $q=1$, $s=\sigma$) whose output is computed using a single Adams \ODEp{}.
	
However, for the purposes of this paper we will instead rewrite a collocation method as a multivalued PBM (i.e. \eqref{eq:polynomial_glm} with $q>1$ and $s=0$), that advances the full ``block'' solution at the set of collocation points forward in time by $h$. To simplify the derivation we focus solely on a Radau IIA method with $\sigma-1$ stages, and show how this integrator can be re-expressed as a PBM with fixed $\alpha$ and $q=\sigma$. Though such a formulation may seem unnecessarily verbose, we will use the multivalued inputs (i.e., solutions at collocation points from the previous time step) in Section \ref{sec:constructing-imex-pbms} to construct a new class of additive IMEX integrators.

\begin{figure}[h!]
	
	\centering
	\includegraphics[width=0.8\linewidth]{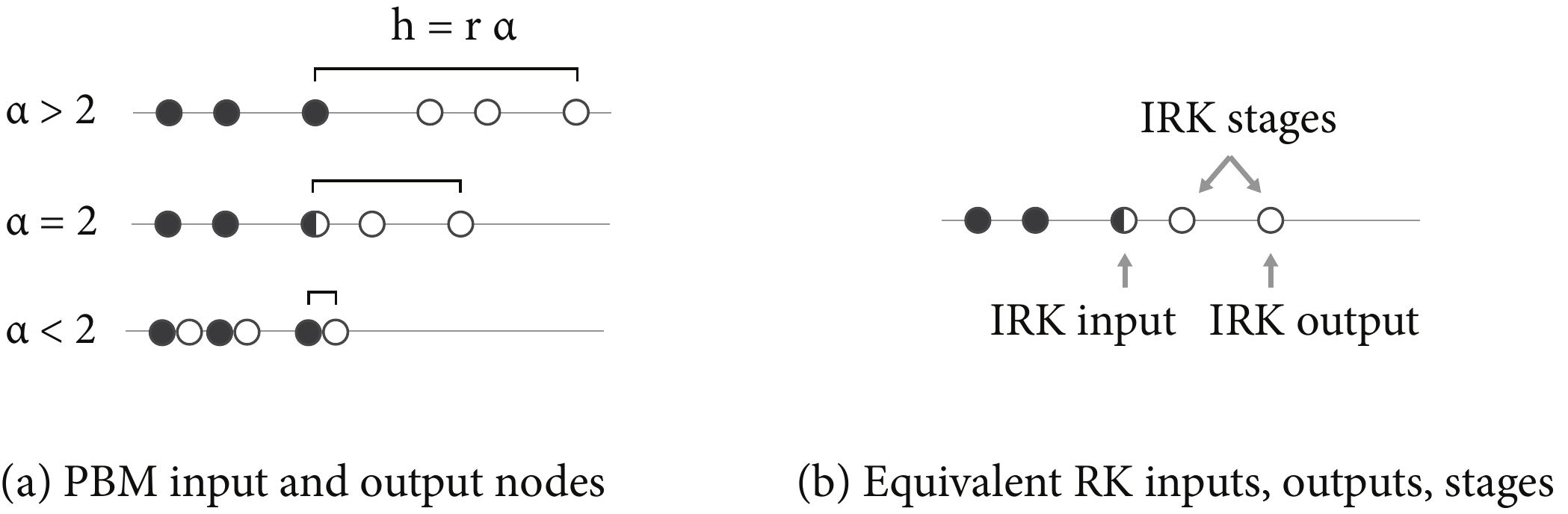}
	\caption{
		Illustration of input nodes 
		\protect\tikz \protect\draw[fill=black,line width=1pt]  circle(.25em);
		and output nodes
		\protect\tikz \protect\draw[fill=white,line width=1pt]  circle(.25em);
		 for the PBM with nodes \eqref{eq:radau-nodes} and $q=3$, such that $\{z_j\} = \{-1, -1/3, 1\}$. The grey lines show the time axis which flows to the right. (a) An illustration of the temporal nodes for the input and outputs for three different values of $\alpha$. When $\alpha=2$ the last input node overlaps with the first output node. (b) A illustration relating the inputs and outputs of the PBM \eqref{eq:fi_pbm} with $b=1$ and $\alpha=2$ to the input, output, and stage values of an IRK method. 
	}
	\label{fig:fi-pbm-nodes}
\end{figure}

We start by selecting the PBM nodes
\begin{align}
	\{z_j\}_{j=1}^q = \{-1, 2x_1 - 1, \ldots, 2x_{q-1} - 1\},
	\label{eq:radau-nodes}
\end{align}
where $x_j$ is the $j$th zero of the polynomial $\frac{d^{q-2}}{dx^{q-2}}(x^{q-2}(x-1)^{q-1})$. Our first quadrature node is $z_1=-1$ and the remaining nodes are the $q-1$ Radau nodes scaled on the interval $x \in [-1,1]$. In Figure \ref{fig:fi-pbm-nodes}(a), we illustrate the input and output nodes for a method with $q=3$.

The outputs of the PBM will all be computed using a single Adams \ODEp{} \eqref{eq:adams_poly_integral_form}. Since we want a method based on implicit Radau quadrature, we select $L_f(\tau)$ to be the interpolating polynomial that matches all the output derivatives at the scaled Radau nodes ($L_f(z_j + \alpha) = f^{[n+1]}_j$ for $j=2,\ldots,q$). To maximize the accuracy of the integration constant $L_y(b)$ we choose  $L_y(\tau)$ to be the interpolating polynomial that matches all the input values ($L_y(z_j) = y^{[n]}_j$ for $j=1,\ldots,q$). If we temporarily leave the expansion point $b$ free, then we obtain the PBM
	\begin{align}
		y_j^{[n+1]} = p(z_j + \alpha; b) = L_y(b) + \int_{b}^{z_j + \alpha} L_f(s) ds, \quad j = 1\ldots,q.
		\label{eq:fi_pbm}
	\end{align}

To obtain a one-step Radau IIA collocation we:
\begin{enumerate}[leftmargin=*]
	\item \ul{Select $\alpha=2$}. The last input node now overlaps with the first output node; see Figure \ref{fig:fi-pbm-nodes}(a). This allows us to use the last input value as an accurate integration constant for Radau quadrature.
	\item \ul{Select $b=1$}. This choice sets the last input $y^{[n]}_q$ as the integration constant since $L_y(1) = y^{[n]}_q$. Furthermore, $y_1^{[n+1]} = p(z_1+2;1) = y_q^{[n]}$, therefore the first output is equal to the last input.
\end{enumerate}
The resulting PBM is equivalent to a $(2q - 3)$th order, one-step Radau IIA collocation method with stepsize $h=2r$. Note that for general $b$ and $\alpha$, the PBM \eqref{eq:fi_pbm} is a $(q-1)$th order, fully-implicit multivalued integrator\footnote{Regardless of the fact that we have Radau nodes, selecting general $b$ will not produce a method with $2\sigma-3$ accuracy. This is due to the fact that our initial condition will only be order $q-1$ accurate, and the lower integration bound is different than the one used for Radau quadrature.}.

The inputs and outputs of the PBM with $\alpha=2$ and $b=1$ are related to the RK inputs and stages according to the following Table which is visualized in Figure \ref{fig:fi-pbm-nodes}(b):
	\begin{center}
	\vspace{2ex}
		\begin{tabular}{r|lll}
			RK step index $\nu$ & RK input $y_\nu$ 		& RK output $y_{\nu+1}$ & RK stages \\ \hline
			$n-1$ 		& $y_{n-1} = y_1^{[n]}$ & $y_{n} = y_q^{[n]}$ & $Y_{j-1} = y_j^{[n]}$, $j=2,\ldots,q$ \\
			$n$	  		& $y_n = y_q^{[n]}$ & $y_{n+1} = y_q^{[n+1]}$ & $Y_{j-1} = y_j^{[n+1]}$, $j=2,\ldots,q$ \\
		\end{tabular}
	\vspace{2ex}
	\end{center}
	
	For clarity we write the PBM for $q=3$. The nodes are $\{z_j\} = \{-1,-1/3,1\}$, and the polynomials are
	\begin{align*}
		L_y(\tau) &= \tfrac{(t-1) (3 t+1)}{4} y_1^{[n]} - \tfrac{9\left(t^2-1\right)}{8} y_2^{[n]} +  \tfrac{(t+1) (3 t+1)}{8} y_3^{[n]}, \\   
		L_f(\tau) &= -\tfrac{3}{4}\left(\tau - \alpha - 1\right)rf_2^{[n+1]} +  \tfrac{3}{4}\left(\tau -\alpha +\tfrac{1}{3}\right) rf_3^{[n+1]}.
	\end{align*}
	Selecting $b=1$ and $\alpha=2$, and substituting $L_y(\tau)$ and $L_f(\tau)$ into (\ref{eq:fi_pbm}) yields
	\begin{align}
		\begin{aligned}
			y^{[n+1]}_1 &= y_3^{[n]} \\
			y^{[n+1]}_2 &= y_3^{[n]} + \tfrac{5}{6} rf_2^{[n+1]} - \tfrac{1}{6} rf_3^{[n+1]}\\
			y^{[n+1]}_3 &= y_3^{[n]} + \tfrac{3}{2} rf_2^{[n+1]} + \tfrac{1}{2} rf_3^{[n+1]}
		\end{aligned}
		\label{eq:radau-IIA-multivalue-form}
	\end{align}
	which reduces to a Radau IIA method with stepsize $h = 2r$ (replacing $r$ with $h/2$ produces the well-known IRK coefficients). Also note that the first two inputs are not used to compute the output, therefore \eqref{eq:radau-IIA-multivalue-form} is equivalent to a one-step method.

\subsection{Iterators}
\label{subsubsec:propagator_iterator}

Before introducing additive integrators, we require one additional concept from the polynomial framework, namely the idea of an iterator. A polynomial integrator only advances the solution if the extrapolation factor $\alpha$ is greater than zero (see Figure \ref{fig:alpha_visualization}). If we let $\alpha = 0$, then we obtain a special method known as an iterator. Iterators recompute the solution at the current timestep and certain method constructions share many similarities with predictor corrector block methods \cite{shampine1969block} and spectral deferred correction iterations \cite{Dutt2000SDC, christlieb2010parallel}. In \cite{buvoli2021epbm} iterators were used in the context of exponential integration to compute initial conditions and create composite methods with improved stability. In this work, we will show how these same ideas can be applied to additive polynomial integrators.

\section{Additive polynomial integrators}
\label{sec:additive-polynomial-framework}
		
		In this section we introduce the additive polynomial time integration framework for solving the partitioned system 
			\begin{align}
				y'(t) = f(t,y(t)) = \sum_{k = 1}^m f^{\{k\}}(t, y(t)).
				\label{eq:model-ode-partitioned}
			\end{align}
			In Subsection \ref{sec:partitioned_odeds_odep} we extend the ODE dataset (\ref{eq:ode_ds}) and the \ODEp{} (\ref{eq:odep_local_coordinates})  from the unpartitioned system \eqref{eq:model-ode} to the partitioned system \eqref{eq:model-ode-partitioned}. Then, in Subsection \ref{sec:additive_pbm} we introduce the class of additive polynomial block methods that will provide a starting point for deriving new integrators in Section \ref{sec:constructing-imex-pbms}.

		\subsection{Partitioned \ODEds{} \& \ODEps{}}
		\label{sec:partitioned_odeds_odep}
			
		A partitioned \ODEd{} contains all the data values that can be used to construct partitioned \ODEps{}. We can trivially extend the \ODEd{} for partitioned equations by replacing the full right-hand-side $f(t,y)$ with all the derivative components $f^{\{j\}}(t,y)$.
			
		\begin{definition}[Partitioned \ODED{}]
			A partitioned \ODEd{} $D(r,s)$ of size $w$ is an ordered set of tuples of the form
			\begin{align}
				D(r,s) = \left\{ \left(\tau_j,~ y_j,~ rf^{\{1\}}_j,~ \ldots,~ rf^{\{m\}}_j \right) \right\}_{j=1}^w	
				\label{eq:local_ode_dataset_partitioned}
			\end{align}
			where $t(\tau) = r\tau + s$, $y_j \approx y(t(\tau_j))$, and $f^{\{k\}}_j = f^{\{k\}}(t(\tau_j), y_j)$.
		\end{definition}

		Before generalizing the \ODEp{} it is convenient to first introduce the notion of a {\em total derivative}.

		\begin{definition}[Total Derivative]
			\label{def:total-derivative}
		 A total derivative $F$ approximates the local ODE solution derivative $ry'(\tau)$ at a point $\tau = \hat{\tau}$, such that $F \approx ry'(t(\hat{\tau}))$, and is the sum of $m$ derivative components, that is
		\begin{align}
			F = rf^{\{1\}}_{j_1} + r f^{\{2\}}_{j_2} + \ldots + rf^{\{m\}}_{j_m}	 \quad \text{where} \quad \hat{\tau} = \tau_{j_1} = \tau_{j_2} = \ldots = \tau_{j_m},
			\label{eq:total-derivative}
		\end{align}
		for indices $j_k \in \{1, \ldots, w\}$ and $k=1,...,m$. 
		\end{definition}
		In words, a total derivative is the sum of $m$ component derivatives from an \ODEd{} where the $k$th component derivative approximates $rf^{\{k\}}(t(\hat\tau),y(t(\hat\tau)))$.

		To generalize the \ODEp{} for partitioned systems, we can reuse (\ref{eq:solution_ode_polynomial}) and simply modify the rules for computing approximate derivatives. 

		\begin{definition}[Partitioned \ODEP{}]
			A partitioned \ODEp{} of degree $g$ with expansion point $b$ can be expressed as	
			\begin{align}
				p(\tau; b) &= \sum_{j=0}^{g}  \frac{a_{j}(b)(\tau - b)^j}{j!},
				\label{eq:local_ode_polynomial_partitioned}
			\end{align}
			where each approximate derivative $a_{j}(b)$ is computed using values from a partitioned \ODEd{} $D(r,s)$ in one of the following ways:

			\begin{enumerate}[leftmargin=*]
				\item {\em  By differentiating a polynomial $h_j(\tau)$ that approximates $y(t(\tau))$}; $h_j(\tau)$ must be a polynomial of least degree that interpolates at least one solution value in $D(r,s)$ and whose derivative $h_j'(\tau)$ interpolates any number of total derivatives from Definition \ref{def:total-derivative}.
		 		The $j$th approximate derivative of the \ODEp{} is then
				\begin{align*}
					a_{j}(b) &= \left. \frac{d^j}{d\tau^j} h_{j}(\tau) \right|_{\tau = b}.
				\end{align*}
				
				\item {\em By differentiating a polynomial $l_j(\tau)$ that approximates $ry'(t(\tau))$}; $l_j(\tau)$ is formed by summing $m$  Lagrange interpolating polynomials so that
					\begin{align*}
						l_j(\tau) = \sum_{k = 1}^{m} l^{\{k\}}_{j}(\tau).
					\end{align*}
					where $l^{\{k\}}_j(\tau) \approx f^{\{k\}}(t(\tau),y(t(\tau))$. Each of the interpolating polynomials $ l^{\{k\}}_{j}$ must satisfy at least one of the derivative component values; specifically, for all $k$ there exists at least one index $\nu \in \{1, \ldots, w\}$, such that $l^{\{k\}}_j(\tau_\nu) = rf^{\{k\}}_\nu$. The $j$th approximate derivative is then
				\begin{align*}
					a_{j}(b) &= \left. \frac{d^{j-1}}{d\tau^{j-1}} l_{j}(\tau) \right|_{\tau = b} &  \text{(only valid for }j \ge 1\text{).}
				\end{align*}	
			\end{enumerate}
		\end{definition}
		
		When the number of components $m=1$, the formulas for the approximate derivatives reduce to those of the classical \ODEps{} described in \cite{buvoli2019constructing}. As is the case for classical \ODEps{}, the family of all additive \ODEps{} is large. Therefore, to simplify method construction, we focus on two special families that depend on significantly fewer free parameters.

	\subsubsection{Special families of partitioned \ODEps{}}
		
		We can generalize the Adams and BDF subfamilies (\ref{eq:adams_approximate_derivatives}) and (\ref{eq:bdf_approximate_derivatives}) for the partitioned \ODEp{} (\ref{eq:local_ode_polynomial_partitioned}). In particular:
		\begin{enumerate}[leftmargin=*]
			\item A {\em partitioned Adams \ODEp{}} has approximate derivatives
			\begin{align}
				a_j(b) = 
				\begin{cases}
					L_y(b) & j = 0 \\
					\left. \frac{d^{j-1}}{d\tau^{j-1}} \sum_{k=1}^m L^{\{k\}}_f(\tau) \right|_{\tau = b} & 1 \le j \le g	
				\end{cases}
				\label{eq:adams_poly_diff_form_partitioned}
			\end{align}
			where:
			\begin{itemize}
				\item $L_y(\tau) \approx y(t(\tau))$ is a Lagrange interpolating polynomial that interpolates only solution data in $D(r,s)$; specifically, 
					for any number of distinct indices ${\nu \in \{1, \ldots, w\}}$, $L_y(\tau_\nu) = y_\nu$.
				\item  Each $L^{\{k\}}_f(\tau) \approx rf^{\{k\}}(t(\tau), y(t(\tau)))$ is a Lagrange interpolating polynomial of degree $g-1$ that satisfies $g$ data values pertaining to the $k$th derivative component; specifically, 
					for $g$ distinct indices $\nu \in \{1, \ldots, w\}$, $L_f^{\{k\}}(\tau_\nu) = rf^{\{k\}}_\nu$.
				
			\end{itemize}
			
			Like the classical Adams \ODEp{} (\ref{eq:adams_poly_integral_form_partitioned}), all partitioned Adams \ODEps{} can be expressed in integral form as
			\begin{align}
			p(\tau;b)	= L_y(b) +  \int_{b}^\tau  \sum_{k=1}^m L^{\{k\}}_f(\xi) d\xi.
			\label{eq:adams_poly_integral_form_partitioned}
			\end{align}

			\item A {\em partitioned BDF \ODEp{}} has approximate derivatives that all satisfy
			\begin{align}
				a_j(b) &= \left. \frac{d^j}{d\tau^j} H_y(\tau) \right|_{\tau = b}	\quad \implies \quad p(\tau; b) = H_y(\tau),
				\label{eq:gbdf_poly_partitioned}
			\end{align}
			where $H_y(\tau) \approx y(t(\tau))$ is an interpolating polynomial of degree $g$ that satisfies $g-1$ solution values $y_{j}$, and whose derivative $H_y'(\tau)$ satisfies one total derivative $F$ from Definition \ref{def:total-derivative}; specifically for $g$ distinct indices $\nu \in \{1, \ldots, w\}$, $H_y(\tau_\nu) = y_\nu$, and for one index $s \in \{1, \ldots, w\}$, $H'_y(\tau_s) = F$.
				
		\end{enumerate}	
		If the number of components $m=1$, then both formulas reduce to the classical Adams \ODEps{} (\ref{eq:adams_approximate_derivatives}) and the classical BDF \ODEps{} (\ref{eq:bdf_approximate_derivatives}).	
		
	\subsection{Additive polynomial block methods}
	\label{sec:additive_pbm}
	
	Block methods \cite{shampine1969block} are multivalued integrators that advance a set (or ``block'') of points at each timestep. If we take $s=0$ in  (\ref{eq:polynomial_glm}), then we obtain a polynomial block method (PBM) \cite{buvoli2019constructing,buvoli2018polynomial,buvoli2020constructing}. PBMs are simpler to derive compared to more general polynomial GLMs since all the \ODEps{} that determine the outputs are  constructed using only input and output data. For this reason, PBMs were used to introduce both classical polynomial methods \cite{buvoli2019constructing} and exponential polynomial methods \cite{buvoli2021epbm}; in this paper we will use PBMs once more to introduce additive polynomial integrators.
	
	An additive PBM depends on the parameters
	\begin{center}
			\renewcommand*{\arraystretch}{1.5}
			\begin{tabular}{llllll}
				$q$ 		& number of inputs/outputs, 	& \hspace{1em} 	& $\left\{ z_j \right\}_{j=1}^q$ & nodes, $z_j \in \mathbb{C}$, $|z_j| \le 1$, \\
				$r$ 		& node radius, $r \ge 0$, 	& 				& $\left\{ b_j \right\}_{j=1}^q$ & expansion points, \\
				$\alpha$ 	& extrapolation factor, 		& 				& \\
			\end{tabular}
	\end{center}
	 and can be written as	
	\begin{align}
		y^{[n+1]}_j &= p_j(z_j + \alpha;~ b_j), & j & =1, \ldots,q,
		\label{eq:additive_polynomial_pbm}
	\end{align}
	where each $p_j(\tau; b)$ is a partitioned \ODEp{} built from the partitioned \ODEd{} 
	\begin{align*}
		D(r,t_n) = 
		\left\{
		\begin{aligned}
			\text{inputs}&: \left\{ \left(z_j,~ y_j^{[n]},~ f_j^{\{1\}[n]}, ~\ldots~,~ f_j^{\{m\}[n]}\right) \right\}_{j=1}^q \\
			\text{outputs}&: \left\{ \left(z_j+\alpha,~ y_j^{[n+1]},~ f_j^{\{1\}[n+1]}, ~\ldots~,~ f_j^{\{m\}[n+1]}\right) \right\}_{j=1}^q. 
		\end{aligned}
		\right. 
	\end{align*}

	Any additive polynomial block method can be written in coefficient form as
\begin{align}
	\mathbf{y}^{[n+1]} = \mathbf{A}(\alpha) \mathbf{y}^{[n]} + r \sum_{k=1}^m \mathbf{B}^{\{k\}}(\alpha) \mathbf{f}^{\{k\}[n]} + r \sum_{k=1}^m  \mathbf{C}^{\{k\}}(\alpha) \mathbf{f}^{\{k\}[n+1]}
	\label{eq:block_coefficient_partitioned}
\end{align}
where the matrices $\mathbf{A}(\alpha)$, $\mathbf{B}^{\{k\}}(\alpha)$, $\mathbf{C}^{\{k\}}(\alpha) \in \mathbb{R}^{q \times q}$, while the solution vectors and $k$th derivative component vectors are  defined using $\mathbf{y}^{[n]} = \left[y_1^{[n]}, \ldots, y_q^{[n]}\right]^T$ and $\mathbf{f}^{\{k\}[n]}_j = \left[f^{\{k\}[n]}_1, \ldots, f^{\{k\}[n]}_q\right]^T$. Any additive polynomial PBM can be recast as an additive general linear method whose parameters are parametrized in terms of the extrapolation parameter $\alpha$. We show this in \cref{sup:additive-polynomial-as-additive-glm} for the doubly partitioned case when $m=2$.

Lastly, additive PBMs are just one example of the more general class of additive polynomial GLMs, described in the following remark.

\begin{remark}[{\bf Additive polynomial GLMs} ]
	The formula for an additive polynomial method with $s$ stages and $q$ outputs is equivalent to \eqref{eq:polynomial_glm} except that each $p_j(\tau; b)$ is now a partitioned \ODEp{} constructed from a partitioned \ODEds{} containing the method's input, output, and stage values. 
The family of additive polynomial GLMs is vast, therefore we cannot provide a detailed exploration of such methods in this work. However, in Section \ref{sec:composite-radau-pbms} we will see one example of an additive polynomial GLM that is constructed by composing multiple PBMs.
\end{remark}

\subsection{Linear stability}
\label{subsec:linear-stability}

Linear stability analysis \cite[IV.2]{hairer1999stiff} characterizes the behavior of a time integrator applied to the Dahlquist test problem $y' = \lambda y$, and is essential for determining which types of problems cause instabilities. To study the linear stability properties of an additive integrator, we use the partitioned Dahlquist equation \cite{ascher1995implicit,sandu2015generalized,izzo2017highly}  
\begin{align}
	y' = \sum_{k=1}^m \lambda_k y, \quad y(0) = y_0,
	\label{eq:partitioned-dahlquist}
\end{align}
where $\lambda_k y$ represents the $k$th component and $\lambda_k \in \mathbb{C}$. The nonlinear problem (\ref{eq:model-ode-partitioned}) reduces to the partitioned Dahlquist equation if all $f^{\{k\}}(t,y)$ are autonomous, diagonalizable linear operators that share the same eigenvectors. We note that this will rarely hold true in practice; nevertheless, this type of analysis is useful for studying and comparing the stability properties of additive methods. 

When applied to  (\ref{eq:partitioned-dahlquist}), the additive polynomial integrator \eqref{eq:block_coefficient_partitioned} reduces to the matrix iteration
	\begin{align}
		\mathbf{y}^{[n+1]} = \mathbf{M}(z_1, z_2, \ldots, z_m, \alpha) \mathbf{y}^{[n]}	
	\end{align}
	where $\mathbf{y}^{[n]}_k = y_k^{[n]}$, $k=1, \ldots, q$, $\mathbf{M}{:}~\mathbb{C}^m\times \mathbb{R}^+ \to \mathbb{C}^{q\times q}$, $z_j = h \lambda_k$, and $h$ is the timestep. The stability region $S$ is the subset of $\mathbb{C}^m \times \mathbb{R}^+$ where $\mathbf{M}(z_1, \ldots, z_m, \alpha)$ is power bounded, so that 
		\begin{align}
			S = \left\{\mathbf{z} \in \mathbb{C}^m, \alpha \in \mathbb{R}^+ ~|~ \sup_{n\in \mathbb{N}} \| \mathbf{M}(\mathbf{z}_1, \ldots, \mathbf{z}_m, \alpha)^n \| < \infty	 \right\}.
		\end{align}
		Therefore, the method is stable if the eigenvalues of $\mathbf{M}(z_1, \ldots, z_m, \alpha)$ lie inside the closed unit disk, and any eigenvalues of magnitude one are non-defective. In Section \ref{sec:constructing-imex-pbms} we will study the stability of additive polynomial integrators with $m=2$, and define several two-dimensional subsets of $S$ that are useful for visualizing the stability region.
				
\section{Constructing FIMEX polynomial integrators based on Radau IIA}
\label{sec:constructing-imex-pbms}

The additive polynomial framework enables the construction of additive PBMs for the general partitioned equation (\ref{eq:model-ode-partitioned}). However, in this introductory work we will focus on constructing Fully-implicit IMEX (FIMEX) methods for the simpler doubly partitioned equation (\ref{eq:model-ode-m2}) where $f^{\{1\}}$ is treated implicitly and $f^{\{2\}}$ is treated explicitly. In particular, we will introduce a family of arbitrary-order FIMEX integrators that are based on the fully-implicit Radau IIA methods described in Section \ref{sec:polynomial-radau}.

To introduce a FIMEX method with an implicit part that is equivalent to Radau IIA, we start by reusing the nodes (\ref{eq:radau-nodes}) where we assume that $q\ge 2$. As in Section \ref{sec:polynomial-radau}, we construct a PBM whose outputs are all computed using a single Adams ODE polynomial $p(\tau; b)$ such that 
\begin{align}
	{y_j^{[n+1]} = p(z_j + \alpha; b_j)}, \quad j = 1, \ldots, q.
	\label{eq:imex-radau-ode-poly-formulation}
\end{align}
Since we are considering an additive integrator for two component systems, the Adams polynomial is now doubly partitioned. To obtain an IMEX integrator, we must treat the first component implicitly and the second component explicitly, that is
\begin{align}
	p(\tau;b)= L_y(b) +  \underbrace{\int_{b}^\tau L^{\{1\}}_f(\xi) d\xi}_{\text{implicit approx.}}  +~ \underbrace{ \int_{b}^\tau L^{\{2\}}_f(\xi) d\xi}_{\text{explicit approx.}}.
			\label{eq:imex_pradau_adams_poly} 
\end{align}

If we want an implicit component that is equivalent to Radau IIA, then we must select $L_y(\tau)$, $L_f^{\{1\}}(\tau)$, $b$ and $\alpha$ identically to what was done in Section \ref{sec:polynomial-radau}, namely: 
	\begin{itemize}
		\item $L_y(\tau)$ is the polynomial of degree $q-1$ that interpolates all $q$ inputs, such that $L_y(z_j) = y_j^{[n]}$, $j = 1, \ldots, q$.
		\item $L^{\{1\}}_f(\tau)$ is the polynomial of degree $q-2$ that interpolates the last $q-1$ output component derivatives $f^{\{1\}[n+1]}_j$ such that $L^{\{1\}}_f(z_j+\alpha) = rf_j^{\{1\}[n+1]}$, ${j = 2, \ldots, q}$.
		\item The extrapolation factor is $\alpha=2$ and the lower integration bound is $b=1$.
	\end{itemize}
To form the explicit approximation, we now make use of the input derivative components $f^{\{2\}[n]}_j$. We propose two different strategies for selecting $L_f^{\{2\}}(\tau)$:
\begin{enumerate}
	\item \ul{Use only the derivative components at the input Radau nodes}; let $L^{\{2\}}_f(\tau)$ be a polynomial of order $q-2$ that satisfies
		\begin{align}
			L^{\{2\}}_f(z_j) &= rf_j^{\{2\}[n]}
			\quad j = 2, \ldots, q.
			\label{eq:imex-radau-L2-poly}
		\end{align}
	
	\item \ul{Use all input derivative components}; let $L^{\{2\}}_f(\tau)$ be a polynomial of order $q-1$ that satisfies
		\begin{align}
			L^{\{2\}}_f(z_j) &= rf_j^{\{2\}[n]}
			\quad j = 1, \ldots, q.
			\label{eq:imex-radau*-L2-poly}
		\end{align}
\end{enumerate}

	Both choices for $L_f^{\{2\}}(\tau)$ lead to the multivalued FIMEX method
	\begin{align}
		y_j^{[n+1]} = y_q^{[n]} + \int_1^{z_j + 2} L_f^{\{1\}}(\xi) + L_f^{\{2\}}(\xi) d\xi, \quad j = 1,\ldots,q,
		\label{eq:poly-imex-radau}
	\end{align}
	where the implicit component is equivalent to Radau IIA; notice that if $f^{\{2\}}(t,y) = 0$, then $L_f^{\{2\}}(\tau) = 0$ and (\ref{eq:poly-imex-radau}) is equivalent to (\ref{eq:fi_pbm}). The order-of-accuracy of \eqref{eq:poly-imex-radau} is equal to the minimum order-of-accuracy of the implicit and explicit component. Because the implicit component is equivalent to Radau IIA with $q-1$ stages, its order-of-accuracy is $2q - 3$. Depending on whether we choose the polynomial \eqref{eq:imex-radau-L2-poly} or \eqref{eq:imex-radau*-L2-poly} the explicit component respectively has an order of $q-1$ or $q$. Therefore, for all $q>2$, the explicit component  determines the overall accuracy and the higher-order polynomial is desirable. However, in Section \ref{sec:radau-linear-stability} we will see that the higher-order polynomial also leads to inferior stability properties.
	
\begin{figure}[ht!]
	
	\centering
	\includegraphics[width=\linewidth]{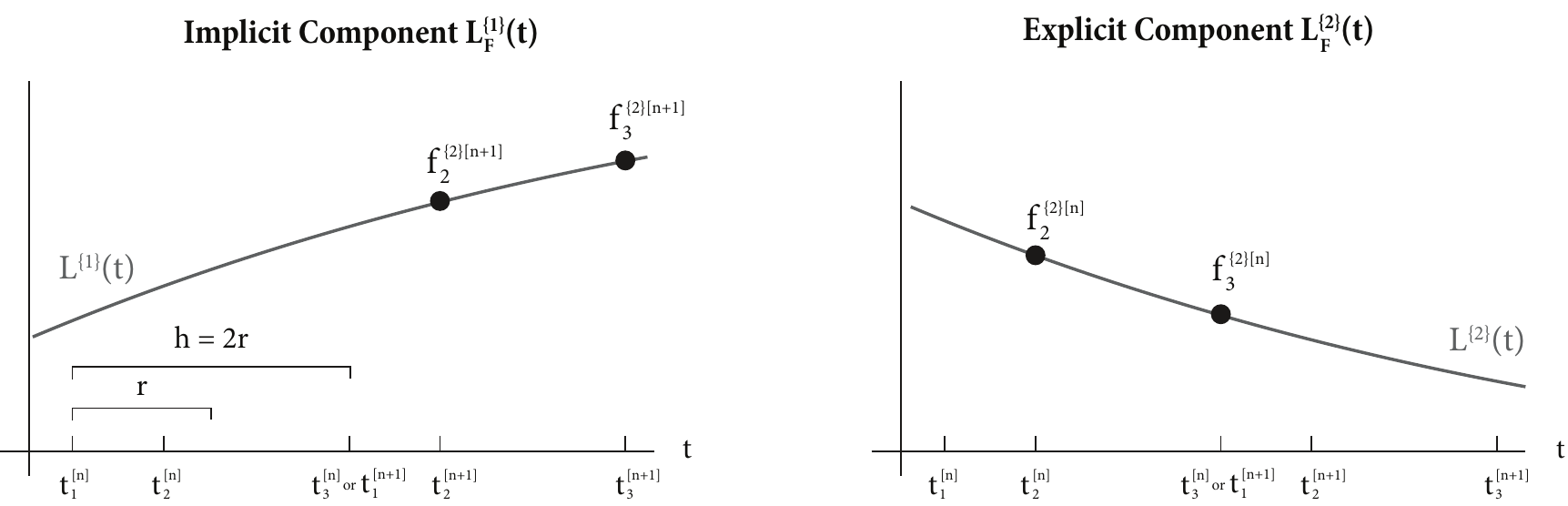}
	\begin{footnotesize}
		\begin{align*}
			\text{\bf Output Formula:} \quad y_j^{[n+1]} = y_{3}^{[n]} + \int_{t^{[n]}_3}^{t^{[n+1]}_j} L^{\{1\}}_f(t) + L^{\{2\}}_f(t) dt
		\end{align*}
	\end{footnotesize}
	\caption{An illustration depicting the polynomials in global time $t$ that determine the outputs of the polynomial FIMEX-Radau method (\ref{eq:poly-imex-radau}) with $q=3$ and nodes $\{z_j\}=\{-1, -1/3, 1\}$. The nodes $t^{[n]}_j$ are the nodes $z_j$ in global coordinates at the $n$th timestep such that $t^{[n]}_j = t_n + r z_j$.
	}
	\label{fig:imex-radau-diagram}
\end{figure}

	From here on, we will refer to the FIMEX method with \eqref{eq:imex-radau-L2-poly} as FIMEX-Radau and the FIMEX method with \eqref{eq:imex-radau*-L2-poly} as FIMEX-Radau*. In Figure \ref{fig:imex-radau-diagram} we present an illustration of the polynomials $L^{\{1\}}_f(\tau)$ and $L^{\{2\}}_f(\tau)$ for a FIMEX-Radau method with $q=3$, and below we also show the coefficients for the method:
	\begin{align}
		\begin{aligned}
			y_1^{[n+1]} &= y_3^{[n]} \\
			y_2^{[n+1]} &= y_3^{[n]} + \tfrac{5}{6} rf_2^{\{1\}[n+1]} -\tfrac{1}{6}  rf_3^{\{1\}[n+1]} - \tfrac{1}{6} rf_2^{\{2\}[n]} + \tfrac{5}{6} rf_3^{\{2\}[n]} \\
			y_3^{[n+1]} &=y_3^{[n]} + \tfrac{3}{2}rf_2^{\{1\}[n+1]}  + \tfrac{1}{2}r f_3^{\{1\}[n+1]} - \tfrac{3}{2}rf_2^{\{2\}[n]} + \tfrac{7}{2}rf_3^{\{2\}[n]}
		\end{aligned}
		\label{eq:radau-eqs}
	\end{align}
	Notice that if $f^{\{2\}}(t,y) = 0$, then (\ref{eq:radau-eqs}) is equivalent to (\ref{eq:radau-IIA-multivalue-form}).
	
	Any FIMEX-Radau or FIMEX-Radau* method with $q$ nodes can be written in the coefficient form
	\begin{align}
			\mathbf{y}^{[n+1]} &= \mathbf{A} \mathbf{y}^{[n]} + r \mathbf{B}^{\{1\}} \mathbf{f}^{\{1\}[n+1]} + r \mathbf{B}^{\{2\}} \mathbf{f}^{\{2\}[n]}
			\label{eq:imex-radau-coefficient-form}
	\end{align}
	where $\mathbf{A}$, $\mathbf{B}^{\{1\}}$, $\mathbf{B}^{\{2\}}$ $\in \mathbb{R}^{q\times q}$. In supplementary materials Section \ref{sup:imex-radau-coefficients} we provide the coefficients for $q=2,3,4$, along with a \textsc{Matlab} script for initializing coefficients for larger $q$. When $q=2$, the FIMEX-Radau method is equivalent to the IMEX-Euler method (\ref{eq:imex-euler-generic}), while FIMEX-Radau* treats the the implicit component with backwards Euler and the explicit component with second-order Adams-Bashforth.  For $q>2$, the nonlinear implicit equations that arise in both FIMEX-Radau methods \eqref{eq:imex-radau-coefficient-form} are analogous to those that arise in standard Radau IIA integration, with a modified right-hand side derived from the explicit component (i.e., linear and nonlinear solvers developed for fully implicit RK \cite{chen14,pazner17,farrell2020irksome,jiao2020optimal,rana2020new,irk1,irk2} naturally apply to FIMEX-Radau).
	
	Note that the method construction described here also applies to different node sets. For example, if we had selected Legendre nodes, then the implicit component would be equivalent to that of the fully implicit Gauss methods. However, Radau nodes have the advantage that they lead to an L-stable implicit component. Lastly, this method construction shares many similarities with the exponential PBM methods based on Legendre nodes from \cite{buvoli2021epbm}. Loosely speaking, we have traded the exponential for a fully implicit approximation and switched from Legendre nodes to Radau nodes. 

\subsection{A FIMEX Radau iterator for obtaining initial conditions}
				
				PBMs with $q>1$, including (\ref{eq:poly-imex-radau}), require multiple inputs at the first time-step. Though a one-step method can be used to compute these solution values, it may not always be possible to match the order of a starting method with that of the PBM. In \cite{buvoli2021epbm} we proposed to use an iterator (a polynomial method with $\alpha=0$ described in Section \ref{subsubsec:propagator_iterator}) to compute the initial conditions of exponential PBMs. The iterator used a discrete exponential Picard iteration to improve the accuracy of an approximate solution. We can reuse the same idea in the context of additive integrators. In fact, it is always possible to construct a polynomial iterator that improves the accuracy of solution values relative to one or more highly accurate inputs  \cite[Ch. 6.2]{buvoli2018polynomial}.

				We can obtain a fully-implicit iterator by  replacing $\alpha=2$ with $\alpha = 0$ in (\ref{eq:poly-imex-radau}) so that the upper integration bound is just $z_j$. However this iterator steps backwards in time, which introduces additional complications. By modifying the method construction so that $b_j = -1$, we can avoid integrating backwards in time. This leads to the iterator
					\begin{align}
						y_j^{[n+1]} = y_1^{[n]} + \int_{-1}^{z_j} L_f^{\{1\}}(s) + L_f^{\{2\}}(s) ds, \quad j = 1\ldots,q,
						\label{eq:poly-imex-radau-iterator}
					\end{align}
				that uses the Radau IIA coefficients for both the implicit and explicit parts. For example, when $q=3$, we obtain the method
				\begin{align}
					\begin{aligned}
						y_1^{[n+1]} &= y_1^{[n]} \\
						y_2^{[n+1]} &= y_1^{[n]} + \tfrac{5}{6} rf_2^{\{1\}[n+1]} - \tfrac{1}{6}  r f_3^{\{1\}[n+1]} + \tfrac{5}{6} rf_2^{\{2\}[n]}  - \tfrac{1}{6} r f_3^{\{2\}[n]} \\
						y_3^{[n+1]} &=y_1^{[n]} + \tfrac{3}{2} rf_2^{\{1\}[n+1]}  +  \tfrac{1}{2} rf_3^{\{1\}[n+1]} + \tfrac{3}{2} rf_2^{\{2\}[n]}  +\tfrac{1}{2}rf_3^{\{2\}[n]}.
					\end{aligned}
					\label{eq:iterator-eqs}
				\end{align}
				Note that both the implicit and explicit coefficients in the iterator \eqref{eq:iterator-eqs} are identical to the implicit coefficients in the propagator \eqref{eq:radau-eqs}.
				If $y_1^{[n]}$ is of order $2q-3$ or higher, then every application of the iterator (\ref{eq:poly-imex-radau-iterator}) will improve the order-of-accuracy of all other solution values by one, up to a maximum order of $2q-3$. Moreover, when a repeatedly applied iteration converges, then the PBM outputs are equivalent to the output and stage values of a fully-implicit Radau IIA  method with $q-1$ stages.
				
				We can therefore use (\ref{eq:poly-imex-radau-iterator}) to compute initial conditions in the following way. We first consider an iterator with quadrature points on the interval $[0, 2]$, instead of $[-1,1]$. The method has identical coefficients, however the first input node at $z_1 = 0$ is now located at $t=t_0$ where the initial condition for \eqref{eq:model-ode} is provided (i.e. $y_1^{[0]} = y(rz_1 + t_0) = y_0$). Next we obtain a zeroth-order estimate for the solution at all the input nodes by temporarily assuming a constant solution such that $y^{[n]}_j = y_0$. Last, we repeatedly apply the iterator to improve the accuracy of the zeroth-order estimate. This procedure provides an accurate initial solution at the times $t=t_0 + rz_j$ for $j=1,\ldots,q$, and can be written abstractly as
				\begin{align}
						\mathbf{y}^{[0]} = M^\kappa (\mathbf{c}), && \mathbf{c}_j = y_0, \quad j = 1, \ldots q,
				\end{align}
				where $\mathbf{c}$ is the initial zeroth-order approximation, and $M^{\kappa}(\mathbf{c})$ is the iterator method applied the $\kappa$ times to initial condition $\mathbf{c}$ such that $M^\kappa(c) = M(M(...(M(c))))$. For a FIMEX-Radau or FIMEX-Radau* method, the iterator should be respectively applied a minimum of $\kappa = q-1$ and $\kappa = q$ times at the first step to match the order of the explicit component.

\subsection{Composite FIMEX Radau methods}
\label{sec:composite-radau-pbms}
				 
				 The method (\ref{eq:poly-imex-radau}) where $L_f^{\{2\}}(\tau)$ is selected using (\ref{eq:imex-radau-L2-poly}) pairs a fully implicit integrator of order $2q - 3$ with an explicit integrator of order $q-1$. Therefore, the combined order is limited to $q-1$. Using (\ref{eq:imex-radau*-L2-poly}) improves the order by one, however the imbalance in accuracy between the implicit and explicit component remains.  One strategy for increasing accuracy further, is to use a composite method that first advances the timestep with the propagator \eqref{eq:radau-eqs} and then corrects the output of the propagator $\kappa$ times using the iterator \eqref{eq:iterator-eqs}. This idea of composite PBMs was introduced in \cite{buvoli2021epbm} for exponential integrators and shares many similarities with spectral deferred correction methods \cite{Dutt2000SDC,Minion2003IMEX,buvoli2019esdc}. In the following sections we show how the same idea leads to composite additive integrators with improved stability and accuracy properties.
				 
				 The composite method can be written abstractly as 
				\begin{align}
					\mathbf{y}^{[n+1]} = M^\kappa ( P  (\mathbf{y}^{[n]}	)),
					\label{eq:composite-pbm}
				\end{align}
				when $P(\cdot)$ denotes the propagator (\ref{eq:poly-imex-radau}) where $L_f^{\{2\}}$ is selected using (\ref{eq:imex-radau-L2-poly}) or (\ref{eq:imex-radau*-L2-poly}), and $M^\kappa(\cdot)$ denotes $\kappa$ applications of the iterator (\ref{eq:poly-imex-radau-iterator}). We will refer to this composite method as FIMEX-Radau($q,\kappa$) if $P$ is (\ref{eq:imex-radau-L2-poly}) and FIMEX-Radau$^*$($q,\kappa$) if $P$ is (\ref{eq:imex-radau*-L2-poly}); $q$ is the number of nodes and $\kappa$ is the number of iterator applications. Since each application of the iterator improves the accuracy order by one, the associated order-of-accuracy for these methods is
				\begin{align}
					\text{FIMEX-Radau}(q,\kappa):& \quad \min(2q-3, q-1+\kappa)	 \label{eq:radau-convergence-order} \\					
					\text{FIMEX-Radau*}(q,\kappa):& \quad \min(2q-3, q+\kappa)	\label{eq:radau*-convergence-order}
				\end{align}
		Last, it is important to note that when using iterative methods to solve the implicit equations (as in numerical PDEs), the implicit solve for the iterator will typically be significantly faster than the propagator, due to a very good initial guess. That is, rather than advance a set of solutions forward in time by $h$ like a propagator, an interator simply increase the accuracy of the \emph{current} solutions by one.

\subsection{Linear stability}
\label{sec:radau-linear-stability}

We now discuss the linear stability properties of FIMEX-Radau$^*$ and FIMEX-Radau. If we select the splitting (\ref{eq:partition-semilinear}), then the stability region of the additive integrator is equivalent to that of the Radau IIA method used for the implicit component. Because all Radau IIA methods are L-stable \cite[IV.5]{hairer1999stiff}, both FIMEX methods will have excellent stability. For the remaining splittings (\ref{eq:partition-semilinear})-(\ref{eq:partition-linearly-in-first}) we must instead consider the five-dimensional stability region  
	\begin{align}
		S = \left\{ z_1 , z_2 \in \mathbb{C} ~|~ \rho \left(\mathbf{M}(z_1,z_2,\alpha=2) \right) \le 1 \right\},
	\end{align}
	where $\rho(\cdot)$ denotes the spectral radius. Due to the high-dimensionality of $S$, it is convenient to consider two-dimensional slices that are formed by fixing $z_1$; 
			\begin{align}
				S(z_1) = \left\{ z_2 \in \mathbb{C} ~|~ \rho \left(\mathbf{M}(z_1,z_2,\alpha=2) \right) \le 1 \right\}.
			\end{align}
			The magnitude and argument of the complex number $z_1$ respectively determine the stiffness of the implicit component, and the degree of diffusion. 
			
			In the PDE setting, $\arg(z_1) = \pi/2$ approximately represents a skew-symmetric advection discretization, while $\arg(z_1)=\pi$ approximately represents a symmetric positive-definite diffusion discretization. An intermediate value, $\arg(z_1) \in (\pi/2, \pi)$, approximately represents a mix of advection and diffusion. However, the stability region $S(z_1)$ will not be symmetric along the real $z_1$ axis (i.e. $S(z_1) \ne S(z_1^*)$). In contrast all real-valued discretizations of PDEs have spectrums that are symmetric about the imaginary axis. We therefore consider the stricter stability region
			\begin{align}
				\hat{S}(z_1) = \left\{ z_2 \in \mathbb{C} ~|~ \max 
				\left[ 
					\rho \left(\mathbf{M}(z_1,z_2,\alpha=2) \right),~
					\rho \left(\mathbf{M}(z_1^*,z_2,\alpha=2) \right)
				 \right] < 1 \right\}.
				\label{eq:stability-region-symmetric}
			\end{align}
			Lastly, it is also useful to consider the region
			\begin{align}
				\tilde{S}(\theta) = \left\{  z_2 \in \mathbb{C}, \gamma \in \mathbb{R}^+ ~|~ \max
				\left[  
					\rho \left(\mathbf{M}(\gamma e^{i\theta},z_2,\alpha=2)\right), ~
					\rho \left(\mathbf{M}(\gamma e^{i\theta},z_2,\alpha=2)\right)
				\right] \right\}	
				\label{eq:stability-region-symmetric-rinf}
			\end{align}
			that contains all the $z_2$ values that ensure stability for any $z_1=\gamma e^{i\omega}$ in the wedge $\gamma \ge0$ and $\theta \le \omega \le 2\pi - \theta$.

\begin{figure}[h] %
	\centering

	\begin{minipage}{0.03\textwidth}
	\end{minipage}
	\begin{minipage}{0.31\textwidth}
		\centering
		\begin{footnotesize}
			\hspace{2em} $\kappa = 0$
		\end{footnotesize}
	\end{minipage}
	\begin{minipage}{0.31\textwidth}
		\centering
		\begin{footnotesize}
			\hspace{2em} $\kappa = 1$
		\end{footnotesize}
	\end{minipage}
	\begin{minipage}{0.31\textwidth}
		\centering
		\begin{footnotesize}
			\hspace{2em} $\kappa = 2$
		\end{footnotesize}
	\end{minipage} \\[1em]

	\begin{minipage}{0.03\textwidth}
		\begin{footnotesize}
			\rotatebox{90}{$\arg(z_1) = 0$}
		\end{footnotesize}
	\end{minipage}
	\begin{minipage}{0.31\textwidth}
		\includegraphics[width=1\textwidth]{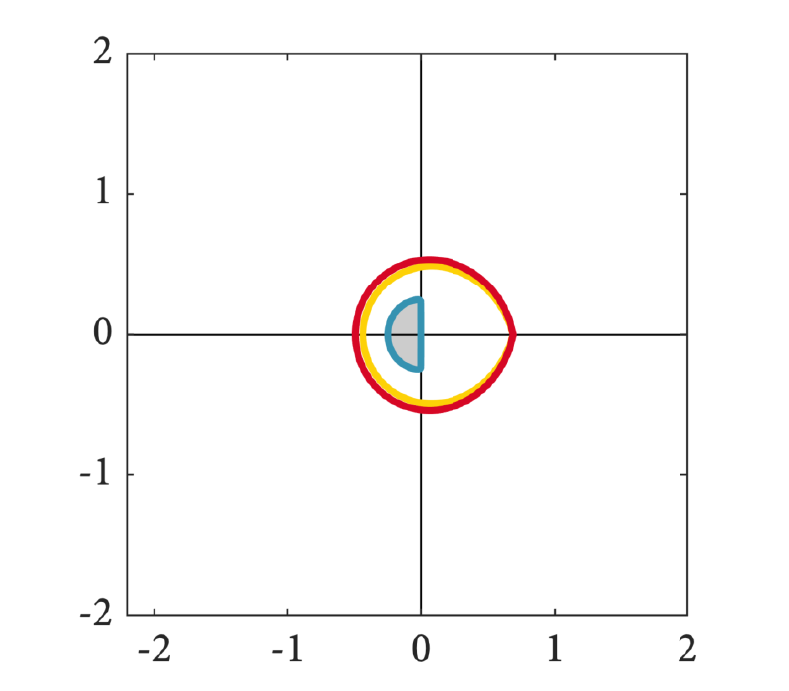}
	\end{minipage}
	\begin{minipage}{0.31\textwidth}
		\includegraphics[width=1\textwidth]{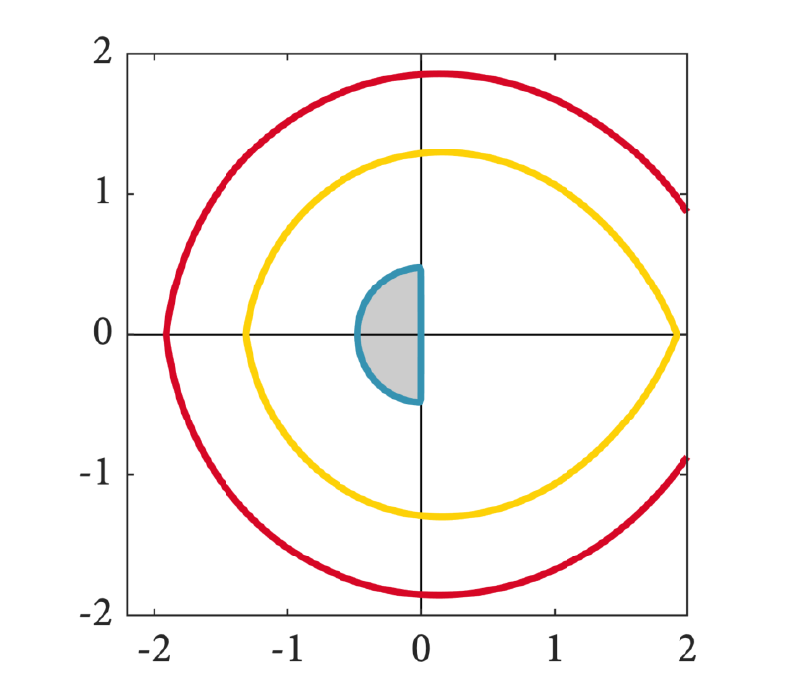}
	\end{minipage}
	\begin{minipage}{0.31\textwidth}
		\includegraphics[width=1\textwidth]{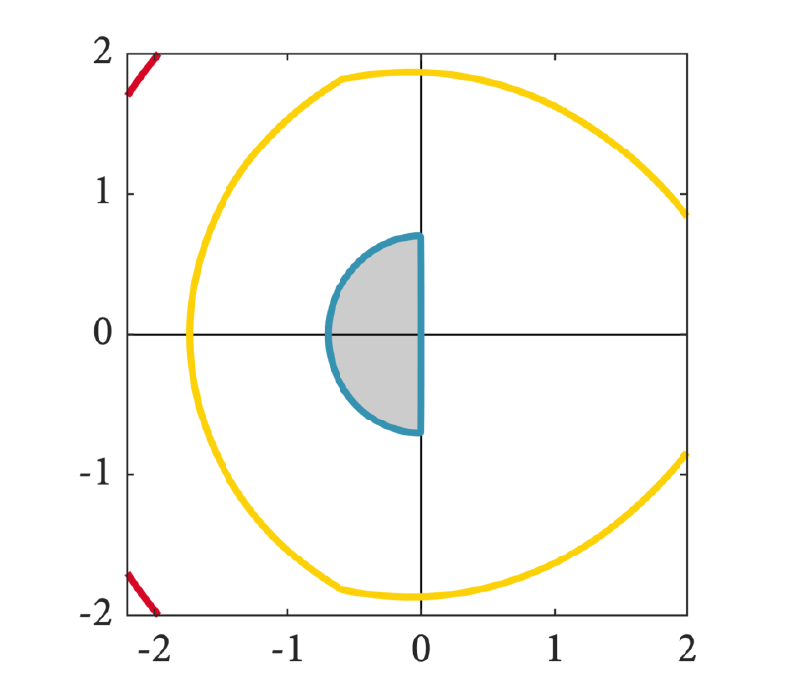}
	\end{minipage}

	\begin{minipage}{0.03\textwidth}
		\begin{footnotesize}
			\rotatebox{90}{$\arg(z_1) = \frac{3\pi}{2}$}
		\end{footnotesize}
	\end{minipage}
	\begin{minipage}{0.31\textwidth}
		\includegraphics[width=1\textwidth]{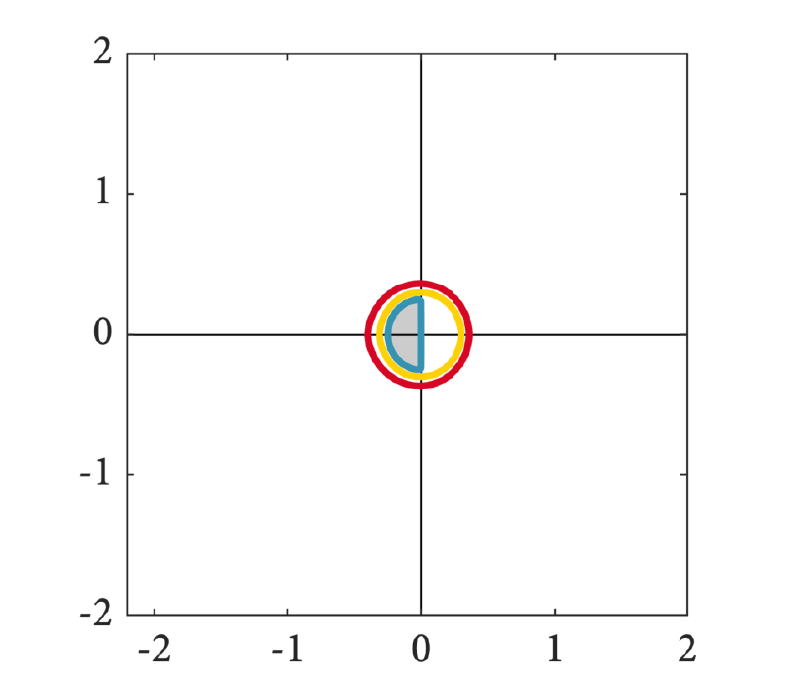}
	\end{minipage}
	\begin{minipage}{0.31\textwidth}
		\includegraphics[width=1\textwidth]{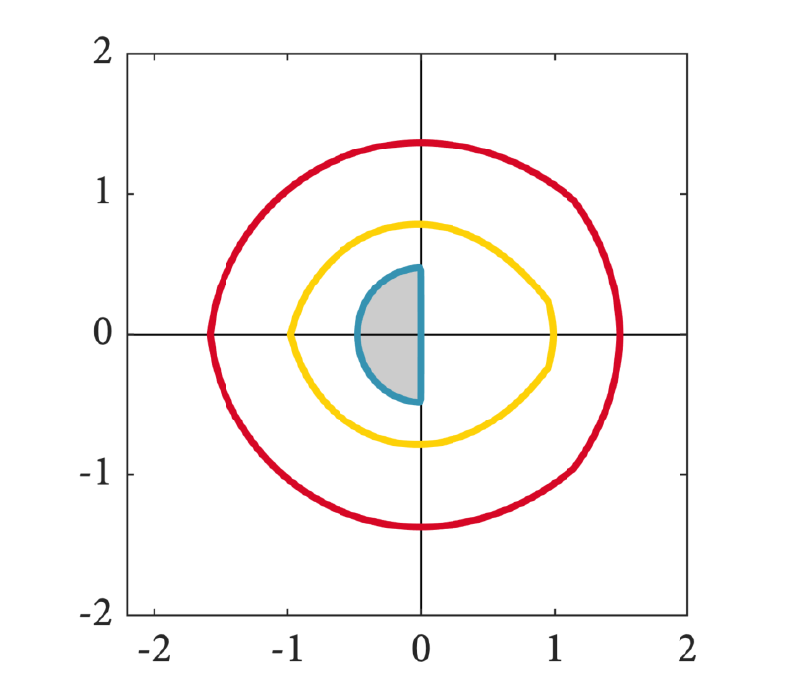}
	\end{minipage}
	\begin{minipage}{0.31\textwidth}
		\includegraphics[width=1\textwidth]{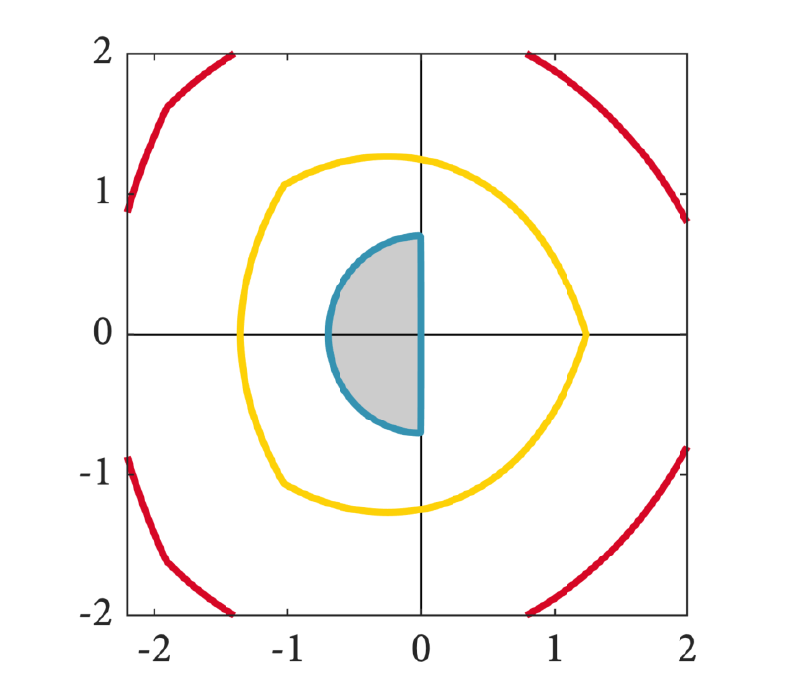}
	\end{minipage}

	\begin{minipage}{0.03\textwidth}
		\begin{footnotesize}
			\rotatebox{90}{$\arg(z_1) = \frac{\pi}{2}$}
		\end{footnotesize}
	\end{minipage}
	\begin{minipage}{0.31\textwidth}
		\includegraphics[width=1\textwidth]{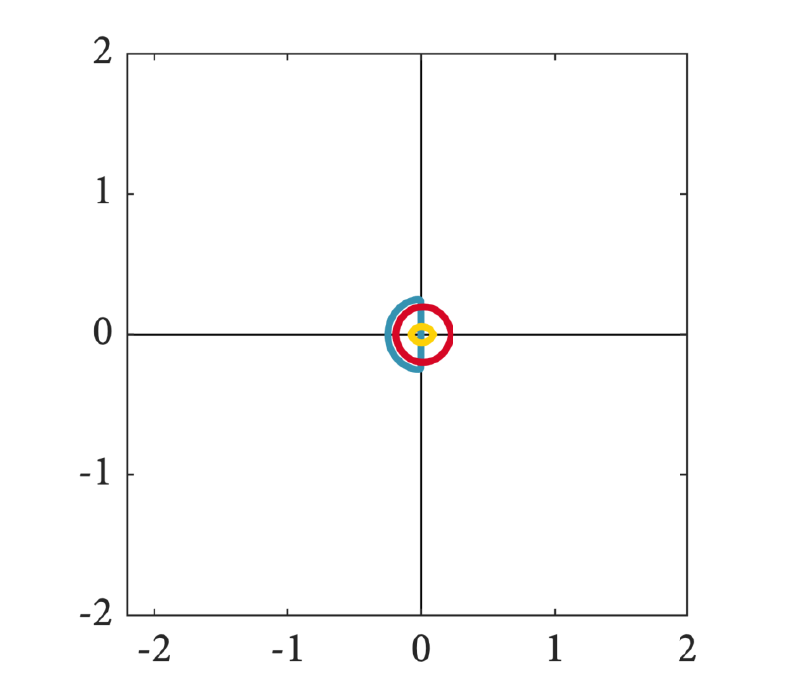}
	\end{minipage}
	\begin{minipage}{0.31\textwidth}
		\includegraphics[width=1\textwidth]{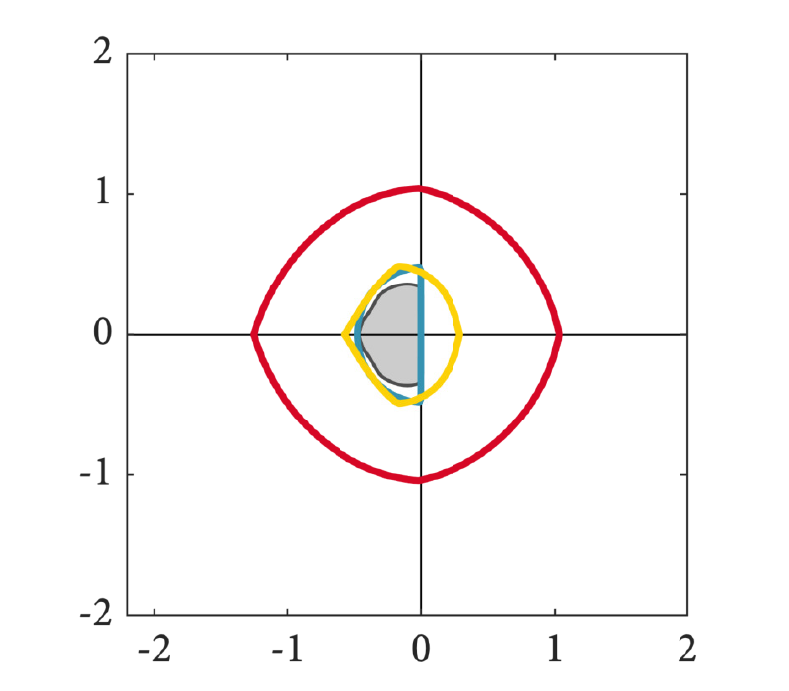}
	\end{minipage}
	\begin{minipage}{0.31\textwidth}
		\includegraphics[width=1\textwidth]{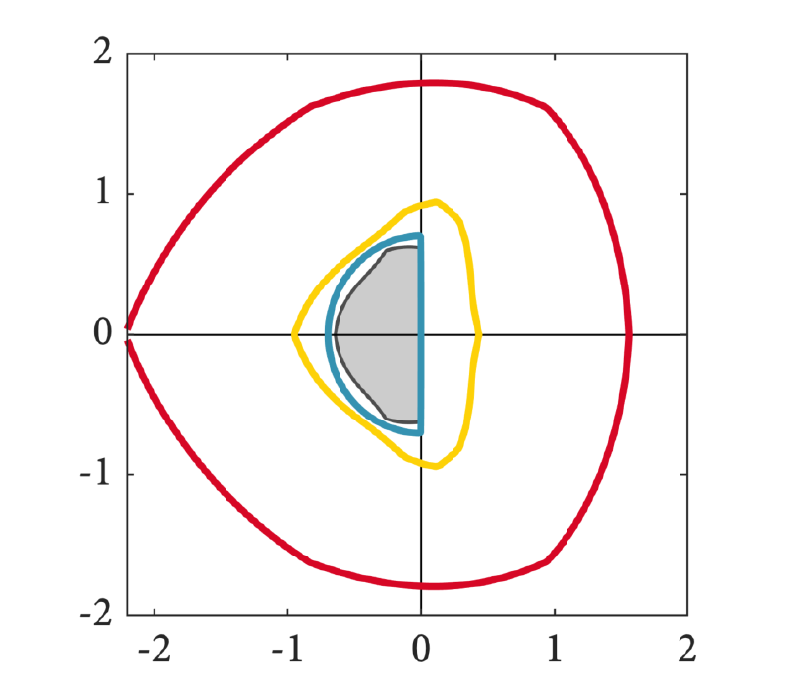}
	\end{minipage}]
	
		\begin{tabular}{c}
{\tiny \textcolor{plot_blue}{\hdashrule[0.2ex]{2em}{2pt}{}} $\left|z_1\right| = 0$ \hspace{1em}}
{\tiny \textcolor{plot_yellow}{\hdashrule[0.2ex]{2em}{2pt}{}} $\left|z_1\right| = 3$ \hspace{1em}}
{\tiny \textcolor{plot_red}{\hdashrule[0.2ex]{2em}{2pt}{}} $\left|z_1\right| = 6$}
\end{tabular}

	\caption{Stability regions for composite FIMEX-Radau*$(q=4,\kappa)$. Each colored contour represents a different stability region $\hat{S}(z_1)$ defined in (\ref{eq:stability-region-symmetric}). The gray region is the stability region $\tilde{S}(\arg(z_1))$ described in (\ref{eq:stability-region-symmetric-rinf}).}
	\label{fig:pradaus-figure}
	
\end{figure}

To present important subsets of the full stability region $S$ we overlay multiple contours of the $\hat{S}(z_1)$ where we take $\arg(z_1) \in \{0, \frac{3\pi}{2}, \frac{\pi}{2}\}$ and $|z_1| \in \{0,3,6\}$. These choices of $\arg(z_1)$ respectively approximate an implicit linear component with diffusion, a mix of diffusion and oscillation, and pure oscillation. In Figure \ref{fig:pradaus-figure} we show these stability regions for the composite method FIMEX-Radau*(4,$\kappa$) from (\ref{eq:composite-pbm}) with $\kappa = 0, 1, 2$. When $\kappa=0$ the FIMEX-Radau*(4,0) method is equivalent to the FIMEX-Radau* method (\ref{eq:poly-imex-radau}) with $q=4$.

The stability regions for FIMEX-Radau* are largest for the diffusive case with $\arg(z_1) = 0$ and smallest for the oscillatory case with $\arg(z_1)=\pi/2$. For each of the three values of $\arg(z_1)$, the FIMEX-Radau*(4,$\kappa$) stability regions always enclose the imaginary axis. Despite the fact that FIMEX-Radau* is equipped with a fully implicit propagator, the stability regions grow slowly in $r$ and are very small when $\arg(z_1)=\frac{\pi}{2}$. However, even a single application of the iterator leads to significantly larger linear stability regions for all three values of $\arg(z_1)$. Lastly, in supplemental materials Figure \ref{supfig:pradau-figure} we also show the equivalent stability regions for the composite FIMEX-Radau method; overall we see that using a lower-order polynomial approximation for the explicit component leads to improved stability.
 
\section{Numerical Experiments}
\label{sec:numerical-experiments}

We investigate the efficiency of FIMEX-Radau and FIMEX-Radau* methods by conducting three types of numerical experiments. First, we study the performance of the integrators on PDEs with periodic domains where inverting the implicit component is trivial. Then, we compare both method families on a finite-element discretization of a non-periodic problem, where solving the fully implicit system requires special care \cite{farrell2020irksome,jiao2020optimal,rana2020new,irk1,irk2}. Lastly, we numerically investigate order reduction on the singularly perturbed Van der Pol equation. 

	For certain spatial discretizations, such as spectral methods, and/or solvers that are difficult to parallelize, the FIMEX-Radau methods allow for certain types of parallelization that are not feasible using existing IMEX-RK methods. For example, when considering periodic domains with Fourier discretizations, the implicit component is trivial to invert, and the majority of the cost is due to the nonlinear function evaluations. In Subsection \ref{subsec:numerical-experiments-pdes-periodic} we demonstrate how parallelization of nonlinear function evaluations can be applied to obtain highly accurate solutions in significantly less time. More generally, there are situations such as when solver performance degrades in parallel, or when finite element construction and nonlinear function evaluations can be very expensive, where the additional parallelization provided by these FIMEX-Radau can be exploited.
	
	Conversely, with finite element discretizations and scalable solvers, spatial parallelism is often very effective. In Subsection \ref{sec:results:dg} we investigate the efficiency of FIMEX-Radau methods with traditional spatial parallelism on a problem where function evaluations are relatively cheap and time-parallelism is not needed.

\subsection{PDEs with periodic boundary conditions} 
\label{subsec:numerical-experiments-pdes-periodic}

In this experiment we evaluate the accuracy and efficiency of composite FIMEX-Radau* methods (\ref{eq:composite-pbm}) on the dispersive Korteweg-De Vries equation with periodic boundary conditions. For comparison we also include results for IMEX-RK integrators of orders one to four; specifically the (1,1,1) and (2,3,2) methods from [Sec. 2.1, 2.5]\cite{ascher1997implicit}, and the ARK3(2)4L[2]SA and ARK4(3)6L[2]SA from \cite{kennedy2003additive}. The  numerical experiment is identical to one found in \cite{buvoli2021epbm,buvoli2019esdc} and is only briefly described below:
	\vspace{0.5em}
\begin{enumerate}[leftmargin=*]
		\item The {\bf Korteweg-de Vries} (KDV) equation from \cite{buvoli2019esdc,buvoli2021epbm}
		\begin{align}
		& \frac{\partial u}{\partial t} = -\left[ \delta \frac{\partial^3 u}{\partial x^3} + \frac{1}{2}\frac{\partial}{\partial x}(u^2) \right] &&
		\begin{aligned}
			& u(x,t=0) = \cos(\pi x), \\
			& x \in [0,2], \nonumber
		\end{aligned}
		\end{align}
		where $\delta = 0.022$. This equation is integrated to time $t=3.6/\pi$ using a 512 point Fourier spectral discretization. 

\end{enumerate}

\subsubsection{Implementation details}

We solve the KDV equation in Fourier space where the initial value problem has the form $\mathbf{y}' = \mathbf{Ly}+N(t, \mathbf{y})$ where $\mathbf{L}$ is an $N \times N$ diagonal matrix that includes the discretized linear differential operator.  For dealiasing we apply the standard 3/2 rule that zeros out the top third of the spectrum. The reference solution $\mathbf{y}_{\text{ref}}$ is computed using a 32nd order exponential spectral deferred correction method \cite{buvoli2019esdc}, and the relative error in our convergence and efficiency studies is defined as $\| \mathbf{y} - \mathbf{y}_{\text{ref}} \|_\infty / \| ]\mathbf{y}_{\text{ref}} \|$.

Since the matrix $\mathbf{L}$ is diagonal, the implicit solve for a fully-implicit propagator or iterator \eqref{eq:imex-radau-coefficient-form} amounts to inverting $N$ decoupled $q\times q$ systems of the form $(\mathbf{I} - r \mathbf{L}_{kk} \mathbf{B}^{\{2\}})\mathbf{x} = \mathbf{b}$ for $k=1,\ldots, N$. To save computational time the matrix inverses for propagators and iterators are precomputed and stored at the first timestep so that an implicit solve can be done using $N$, $q\times q$ matrix multiplications. Since this operation is cheap, the dominant computational cost is due to the nonlinear function evaluations.

FIMEX-Radau integrators benefit from time parallelism since the nonlinear function evaluations and output computations can be computed simultaneously. To quantify the benefits of parallelism, we created both serial Fortran implementations of the integrators and parallel Fortran implementations using OpenMP. A FIMEX-Radau method with $q$ nodes requires $q$ independent function evaluations that can be computed simultaneously using $q$ shared or distributed memory processes / threads. In contrast an IMEX-RK method requires sequential evaluation of the nonlinear term at each stage, therefore this type of parallelization is not possible. However, all methods can benefit from spatial parallelization where the FFT evaluations are computed using multiple threads. For the purposes of this experiment we will only investigate time parallelism and compute all FFT evaluations using one thread.

Our Fortran code can be found in \cite{BuvoliAdditivePBMGithub}, and all the timing results presented in this subsection were produced on a 14 core, 2.0 Ghz Intel Xeon E5-2683 v3 with hyper-threading enabled.

\subsubsection{Results and discussion}

\begin{figure}
	\centering
	\includegraphics[height=0.42\linewidth,trim={0cm 0cm 1cm 0},clip]{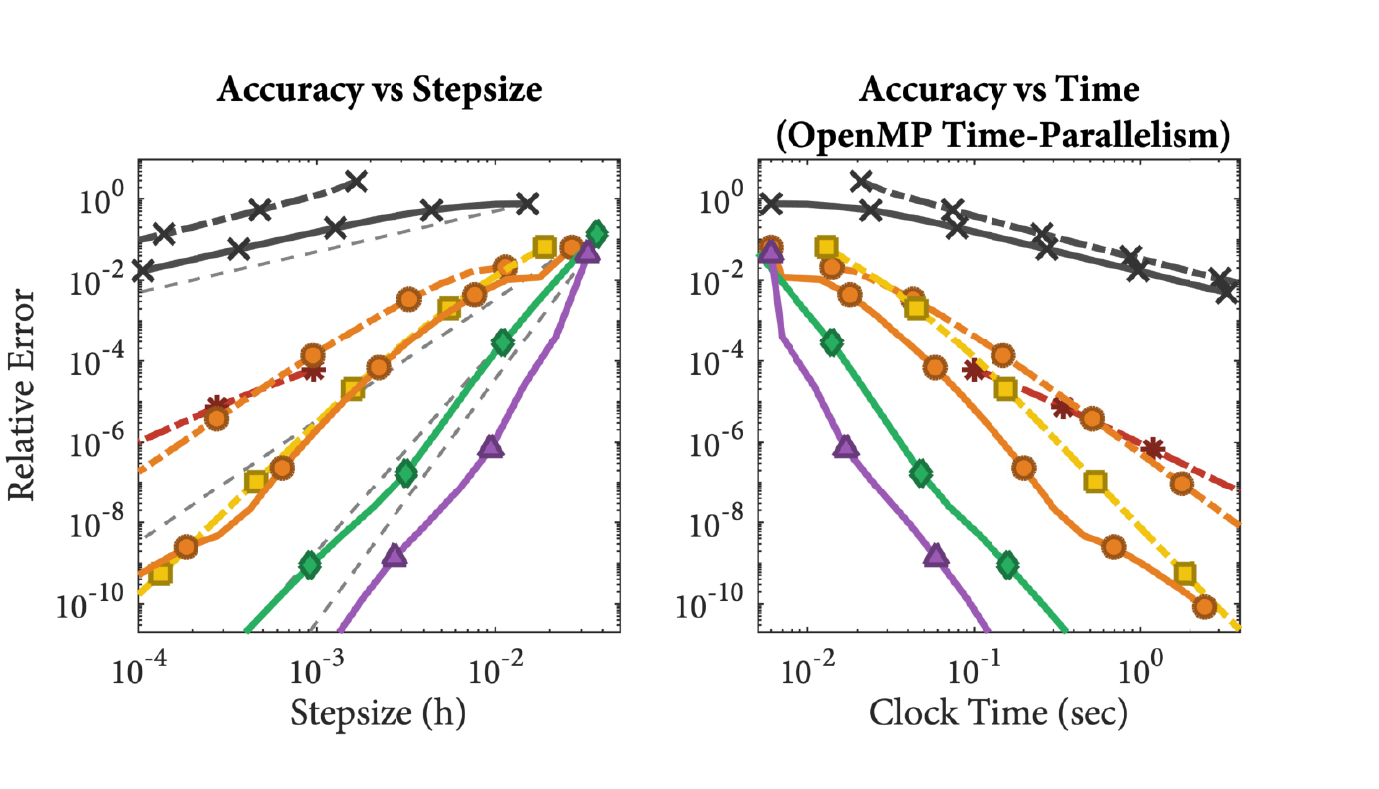} 
	\includegraphics[height=0.42\linewidth,trim={7cm 0cm 1.25cm 0},clip]{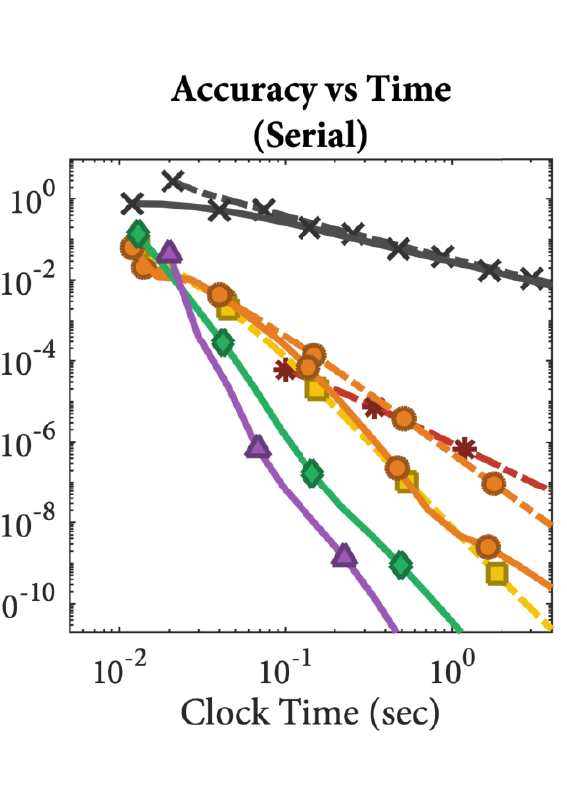}

	\begin{scriptsize}
		\begin{tabular}{rl}
			Order of Accuracy:  & \includegraphics[height=0.025\linewidth,trim={1.72cm 7.05cm 1.9cm .5cm},clip]{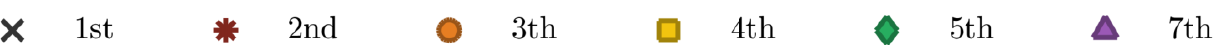} \\[0.25em]
			Integrator Family: 	& \includegraphics[height=0.025\linewidth,trim={5cm 7.05cm 1.9cm .5cm},clip]{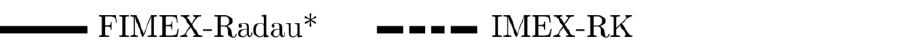}
		\end{tabular}
	\end{scriptsize}
	
	\caption{Accuracy and precision diagrams comparing composite FIMEX-Radau$^*$($q,2$) with ${q=2,3,4,5}$ to IMEX-RK on the KDV equation. All curves start at the first stable timestep. The parallel and serial time plots respectively show run times for composite FIMEX-Radau$^*$($q$,2) methods with and without OpenMP parallelization. The dashed lines of increasing slope in the accuracy plot respectively correspond to first, third, fifth, and seventh order convergence which are the expected orders of convergence (\ref{eq:radau*-convergence-order}) for FIMEX-Radau$^*$.}
	\label{fig:convergence-accuracy-kdv-radau*}
	
\end{figure}

In Figure \ref{fig:convergence-accuracy-kdv-radau*} we present results for the KDV equation solved using composite FIMEX-Radau*($q,2$) methods with $q=2,3,4,5$. We select $\kappa=2$ since the linear stability regions from Figure \ref{fig:pradaus-figure} revealed that FIMEX-Radau*($q$,0), have comparatively poor stability for non-diffusive equations. Overall the FIMEX-Radau*($q$,2) integrators performed excellent, with the serial, seventh-order FIMEX-Radau*(5,2) outperforming all the other methods for any accuracy below $10^{-2}$. When using OpenMP parallelism with five threads, the FIMEX-Radau*(5,2) is the most efficient method across all accuracies, and is capable of obtaining the solution with a relative error of $10^{-2}$ approximately two times faster than any RK method. For more accurate solutions with errors below $10^{-4}$, the difference was even more significant and FIMEX-Radau*(5,2) was able to obtain the solution approximately thirty times faster than IMEX-RK4. For lower-order methods, the third-order FIMEX-Radau*(3,2) method exhibited increased fourth-order convergence throughout most of the stepsize range. This enabled the method to be more efficient than the IMEX-RK4 method, despite only requiring the storage of three solution vectors instead of the five needed by the RK method. Lastly, if very inaccurate solutions are sufficient, then the FIMEX-Radau*(2,1) method was both more efficient and more stable than the IMEX-RK1 (or equivalently FIMEX-Radau(2,0)) method.

\subsection{Numerically investigating order reduction}

In this experiment we numerically investigate order reduction for the composite FIMEX-Radau and FIMEX-Radau$^*$ method (\ref{eq:composite-pbm}) by solving the Van der Pol equation  \cite[p. 403]{hairer1999stiff}
\begin{align}
		\begin{aligned}
			y_1' &= y_2 \\
			y_2'  &= \frac{(1 - y_1^2)y_2 - y_1}{\epsilon} \\
		\end{aligned} &&
		\begin{aligned}
			y_1(0) &=2 \\
			y_2(0) &= -\tfrac{2}{3} + \tfrac{10}{81}\epsilon - \tfrac{292}{2187} \epsilon^2 - \tfrac{1814}{19683} \epsilon^3
		\end{aligned}
		\label{eq:vanderpol}
	\end{align}
integrated to time $t=0.5$. We consider $\epsilon$ values ranging from ${\epsilon = 1}$ to ${\epsilon = 10^{-8}}$. When the stiffness parameter $\epsilon$ is small, this equation is known to cause order-reduction for IMEX-RK methods \cite{ascher1997implicit,boscarino2007error,kennedy2003additive,layton2005implications,izzo2017highly,kennedy2019higher}. 

We integrate the Van der Pol equation using composite methods (\ref{eq:composite-pbm}) where the propagator uses \eqref{eq:imex-radau-L2-poly} or \eqref{eq:imex-radau*-L2-poly}, ${q \in \{3,4,5\}}$ and $\kappa \in \{0, 1, 2\}$. The equation right-hand-side is split in two different ways. A {\em semi-implicit} splitting treats the first component explicitly and the second component implicitly; this same splitting was used in \cite{izzo2017highly, ascher1997implicit,boscarino2007error,layton2005implications,kennedy2003additive,kennedy2019higher}. We also consider the semi-linear splitting (\ref{eq:partition-approximate-jacobian}) where $J_n$ is the exact Jacobian of the right-hand-side at the $n$th timestep. 

To verify FIMEX-Radau($q$,$\kappa$) and FIMEX-Radau$^*$($q$,$\kappa$) methods, we estimate their convergence rates as a function of $\epsilon$ using the approach proposed in \cite{kennedy2019higher}. Specifically, we solve (\ref{eq:vanderpol}) using 30 logarithmically spaced stepsizes ranging from $h=0.25$ to $h=10^{-4}$ and then compute a linear least-squares fit of log(error) versus log(stepsize). Error for an approximate solution $\mathbf{y}$ is measured using $\|\mathbf{y} - \mathbf{y}_{\text{ref}}\|_\infty$ where $\mathbf{y}_{\text{ref}}$ was computed using the \textsc{Matlab} ode15s integrator. For the semi-implicit splitting we set the tolerance of our Newton iteration to $10^{-12}$ and use the \textsc{Matlab} backslash function to solve the associated linear systems. For the linearly-implicit system we also use \textsc{Matlab} backslash function to solve the linear systems at each timestep.

 In Figure \ref{fig:vanderpol-convergence-rate-radaus} we show convergence rate plots and convergence diagrams for the FIMEX-Radau$^*$($q$,$\kappa$) methods with each splitting. Convergence diagrams for different $\epsilon$ values are shown in Supplemental materials \cref{supfig:vanderpol-convergence-rate-imex-radaus-extra}. We also show a convergence rate plot and convergence diagrams for the FIMEX-Radau($q$,$\kappa$) in \cref{supfig:vanderpol-convergence-rate-imex-radau-extra}. There are several important points regarding the results:

\begin{figure}
	\centering

	\begin{minipage}{0.48\textwidth}
		\centering
		\begin{footnotesize}
			{\bf Semi-Implicit Splitting}	
		\end{footnotesize}

		\includegraphics[width=1\linewidth]{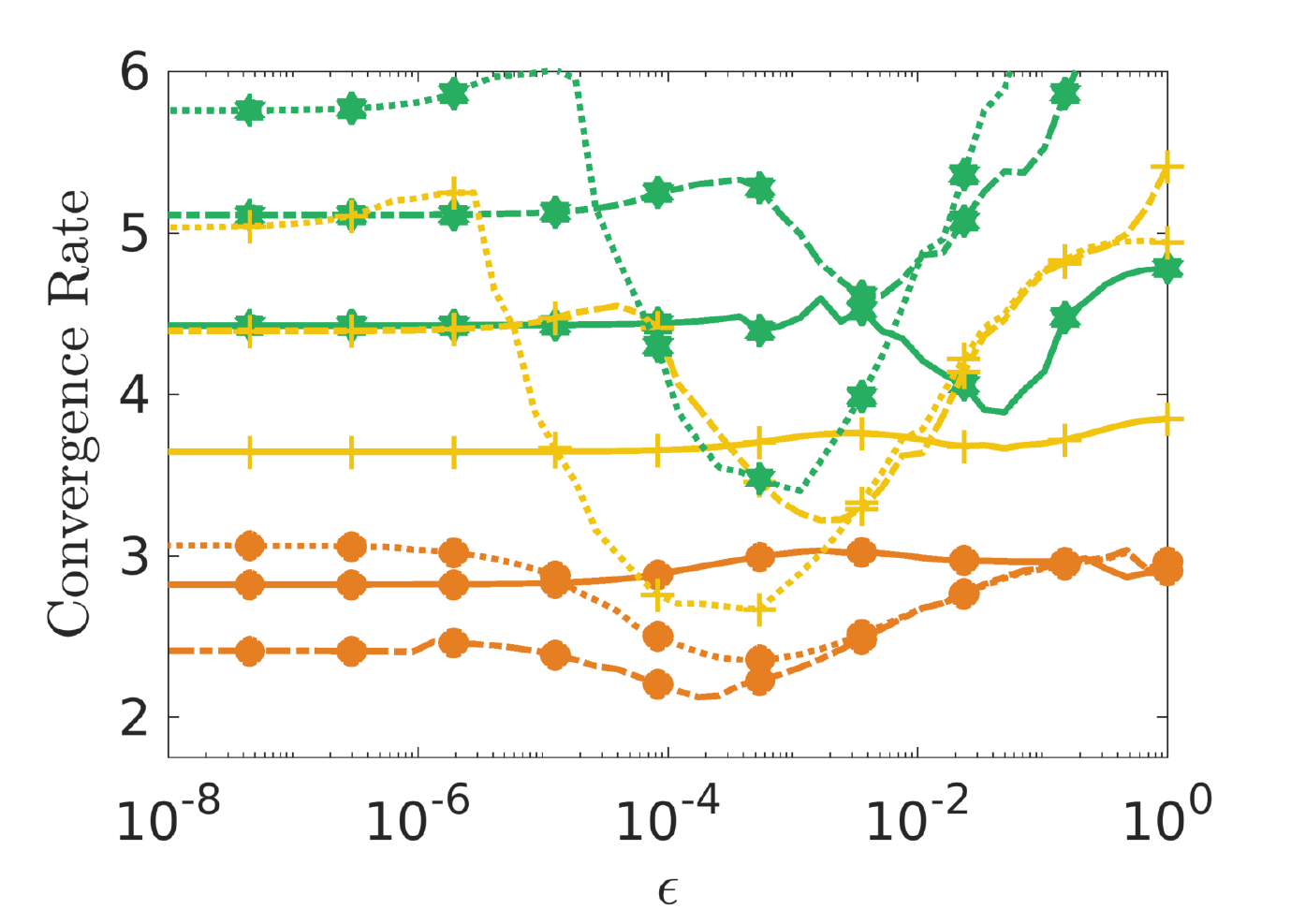}
	\end{minipage}
	\begin{minipage}{0.48\textwidth}
		\centering
		\begin{footnotesize}
			{\bf Linearly-Implicit Splitting}	
		\end{footnotesize}

		\includegraphics[width=1\linewidth]{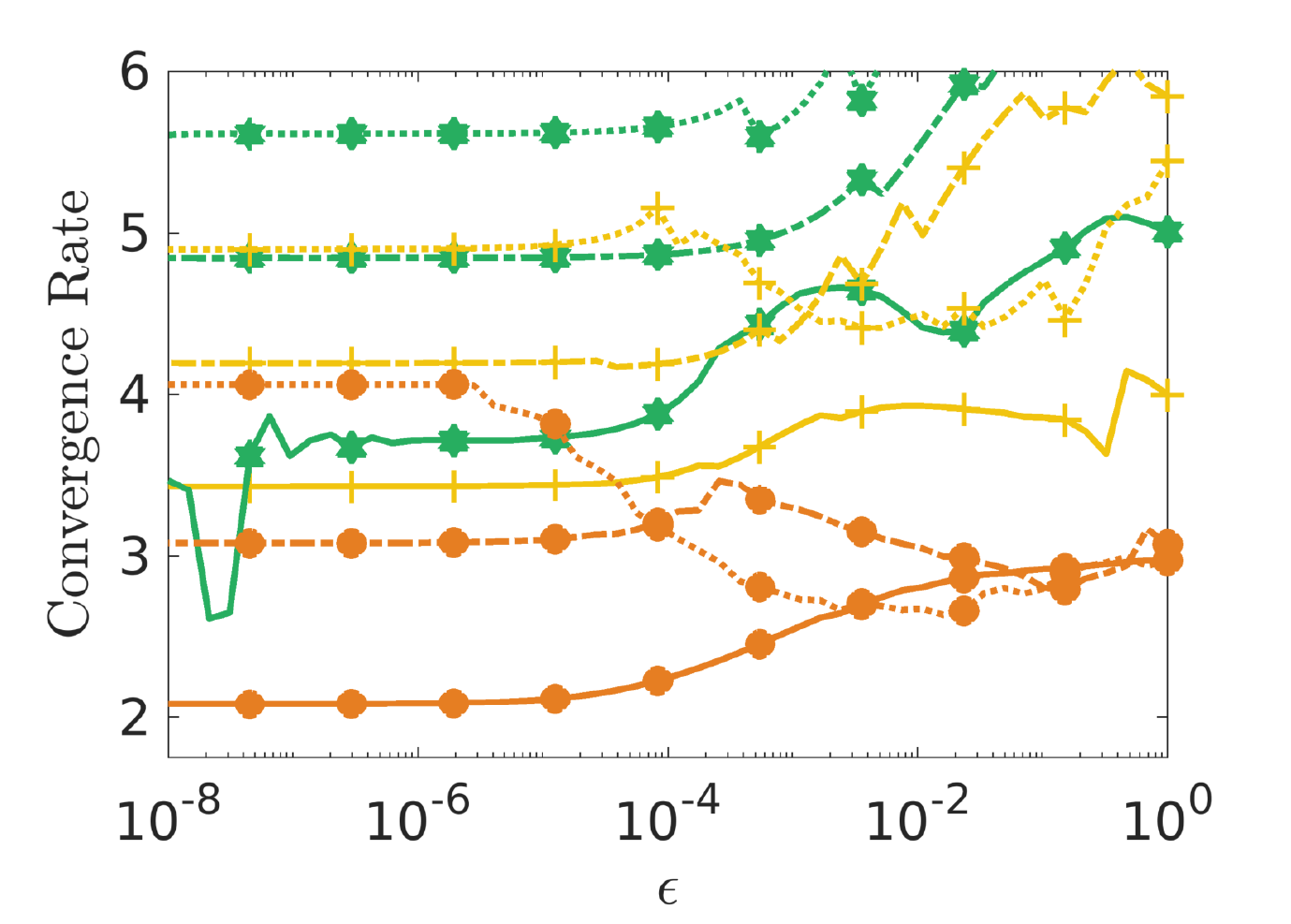}
	\end{minipage}
		
	\includegraphics[width=0.55\textwidth,trim={3.3cm 2.5cm 3.3cm 0},clip]{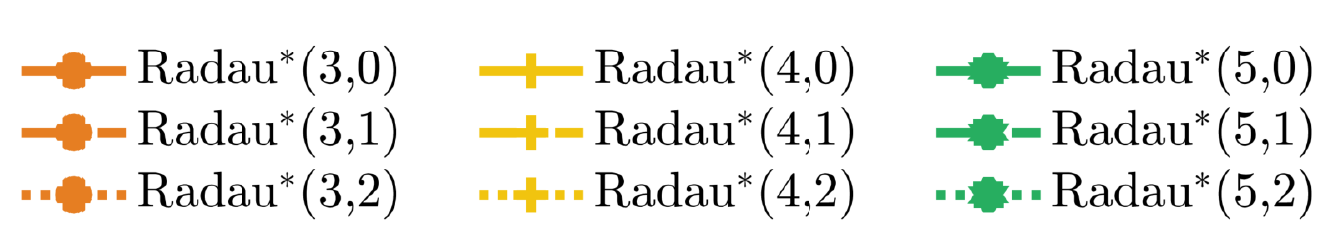}

	\caption{Approximate overall convergence rates for the Van der Pol equation (\ref{eq:vanderpol}) as a function of $\epsilon$ for the FIMEX-Radau$^*$ method with a semi-implicit splitting or a linearly-implicit splitting.}
	\label{fig:vanderpol-convergence-rate-radaus}
\end{figure}

\begin{itemize}[leftmargin=*]
	\item The FIMEX-Radau$^*$ and FIMEX-Radau methods are stable across the full range of stepsizes for both splittings. The linearly-implicit splitting does not require a nonlinear solve at each step, leading to improved computational efficiency. In contrast, all the IMEX-RK methods we tested from \cite{kennedy2019higher,kennedy2003additive}  required substantially smaller timesteps to remain stable with the linearly-implicit splitting (See Supplemental Figure \ref{supfig:vanderpol-convergence-rate-imex-rk}).

	\item For small $\epsilon$, each application of the iterator raises the order of convergence by approximately one. When using a semi-implicit splitting, we see significantly decreased convergence rate for Radau($q$,2) and Radau$^*$($q$,2) methods when $\epsilon$ is between $10^{-2}$ and $10^{-4}$. However the extra iteration improves accuracy significantly. For example, from the convergence diagrams for $\epsilon=10^{-3}$ in \cref{supfig:vanderpol-convergence-rate-imex-radau-extra,supfig:vanderpol-convergence-rate-imex-radaus-extra}, we see that methods with $\kappa=2$ are always more accurate than those with $\kappa=1$ or $\kappa=0$, despite their reduced convergence rates. 
	\item The estimated convergence rates are consistently smaller than the expected orders-of-accuracy \eqref{eq:radau-convergence-order} and \eqref{eq:radau*-convergence-order}. This is due to minor order reduction at stepsizes $h \ge 0.1$ (less than 5 total timesteps).  By including this data in our least squares fits we see lower overall convergence rates. If we only consider $h<0.1$, then we would observe expected convergence (for example, we can see that the convergence curves for $h<0.1$ match the dashed convergence-order lines in the diagrams in \cref{supfig:vanderpol-convergence-rate-imex-radaus-extra,supfig:vanderpol-convergence-rate-imex-radau-extra}). 
	\item The convergence rates for FIMEX-Radau$^*$(3,2) and FIMEX-Radau(3,2) with a linearly-implicit splitting are unusually high. This is due to rapid, convergence at coarse timesteps that leads to higher than expected overall convergence. Nevertheless, for sufficiently small stepsizes, the convergence limits to third-order (e.g. see the convergence diagram for $\epsilon = 10^{-3}$ in \cref{supfig:vanderpol-convergence-rate-imex-radaus-extra,supfig:vanderpol-convergence-rate-imex-radau-extra}).
	\item The FIMEX-Radau$^{*}$(5,0) method with a linearly-implicit splitting did not converge at large timesteps (though it was stable). Increasing $\kappa$ to one resolves the issue.
	\end{itemize}

\subsection{DG advection-diffusion-reaction}\label{sec:results:dg}

Here we consider the time-dependent advection-diffusion-reaction equation
\begin{equation} \label{eq:adv-diff}
	u_t + \nabla \cdot ( \boldsymbol{\beta} u  - \epsilon \nabla u ) + \gamma u^2= f,
\end{equation}
where $\boldsymbol{\beta}(x,y) : = (1,1)^T$
is the velocity field and $\epsilon$ and $\gamma$ constant diffusion and reaction coefficients.

We discretize \eqref{eq:adv-diff} in space using discontinuous
Galerkin finite elements over the spatial domain $\Omega = [0,1] \times [0,1]$.
Dirichlet boundary conditions are weakly enforced on $\partial\Omega$,
advection terms are upwinded \cite{Cockburn2001}, and diffusion terms
are treated with the symmetric interior penalty method
\cite{Arnold1982,Arnold2002}. Let $V_h$ be the DG finite element
space consisting of piecewise polynomials of degree $p$ defined 
locally on elements of the spatial mesh $\mathcal{T}$.
The resulting finite element problem
is to find $u_h \in V_h$ such that, for all $v_h \in V_h$,
{\small
\begin{equation*}
	\begin{split}
	\int_\Omega \partial_t (u_h) v_h \, dx
	- \int_\Omega u_h \boldsymbol{\beta} \cdot \nabla_h v_h \, dx
	+ \int_\Gamma \widehat{u_h} \boldsymbol{\beta} \cdot \llbracket v_h \rrbracket \, ds
	+ \int_\Omega \nabla_h u_h \cdot \nabla_h v_h \, d x 
	+ \int_\Omega \gamma u_h v_h^2 \, dx \\
	- \int_\Gamma \{ \nabla_h u_h \} \cdot \llbracket v_h \rrbracket \, ds
	- \int_\Gamma \{ \nabla_h v_h \} \cdot \llbracket u_h \rrbracket \, ds
	+ \int_\Gamma \sigma \llbracket u_h \rrbracket \cdot \llbracket v_h \rrbracket \, ds
	= \int_\Omega f v_h \, dx,
	\end{split}
\end{equation*}
}where $\nabla_h$ is the broken gradient, $\Gamma$ denotes the skeleton of the mesh,
$\{ \cdot \}$ and $\llbracket \cdot \rrbracket$ denote the average and jump
operators, respectively, and $\widehat{u_h}$ is used to denote the
upwind numerical flux. The parameter $\sigma$ is the \textit{interior
penalty parameter}, important for a stable discretization
\cite{Arnold2002}, which we set to $\sigma = (p+1)^2/h$.

The discretization is implemented in the MFEM finite element library
\cite{Anderson2020}. Classical algebraic multigrid (AMG) in the \emph{hypre}
library \cite{hypre} (BoomerAMG) is used to solve the implicit equations,
with Falgout coarsening, classical interpolation, a strength tolerance
of 0.2, and hybrid parallel Gauss-Seidel relaxation. For the
Radau schemes, we must solve an implicit system analogous to that
which arises in fully implicit Runge Kutta methods. We wrap the
AMG solver with the block preconditioning for fully
implicit Runge Kutta methods developed in \cite{irk1,irk2} to solve
the Radau stage equations. We use 4th-order elements in space, with
mesh spacing $h_x\approx 0.0078$ leading to expected spatial accuracy
$\sim\mathcal{O}(10^{-9})$. We consider two problems and splittings:
\begin{enumerate}
	\item No reaction, treat advection explicitly and diffusion and source
	term implicitly (\Cref{sec:results:dg:ad}). 
	\item $\gamma=10$; treat nonlinear reaction explicitly and
	advection-diffusion implicitly (\Cref{sec:results:dg:adr}).
\end{enumerate}
All simulations are run on 4 dual socket Intel Xeon Gold 6152 22-core processors (i.e., 44 cores/node).

\subsubsection{Advection-diffusion}\label{sec:results:dg:ad}

This section considers no reaction ($\gamma = 0$) and an additive splitting of \eqref{eq:adv-diff} treating the advection explicitly and diffusion implicitly. We choose the forcing function $f(x,y,t)$ such that the analytical solution is given by $u_*(x,y,t) = \sin(2\pi x(1-y)(1+2t))\sin(2\pi y(1-x)(1+2t))$. Here we solve all linear systems to relative residual tolerance of $10^{-12}$. All simulations are run to final (simulation) time tf = 2 120 total MPI processes.

Results in this section indicate superior stability and accuracy of FIMEX-Radau* methods over traditional IMEX-RK. However, results are somewhat artificial in this setting in that (i) there exist fast AMG solvers for advection-diffusion, even in the advection-dominated regime \cite{air2,sivas2021air}, and (ii) due to stability issues with the additive splitting of this discretization shown below, even classical AMG in hypre \cite{hypre} works well as-is applied to the advection-diffusion ratios considered here. That is, although FIMEX-Radau* methods show advantages over classical IMEX-RK, in practice one would probably treat this particular problem and discretization implicitly. Other related advection-diffusion equations however, such as a DG-BDM discretization of incompressible Navier Stokes, can be very difficult to solve implicitly for non-trivial advection, so it is possible FIMEX-Radau* methods have a use in such advection-diffusion problems, where implicit treatment is difficult or not viable. 

\begin{figure}[!htb]
  \centering
  \begin{subfigure}[b]{0.25\textwidth}
    \includegraphics[width=\textwidth]{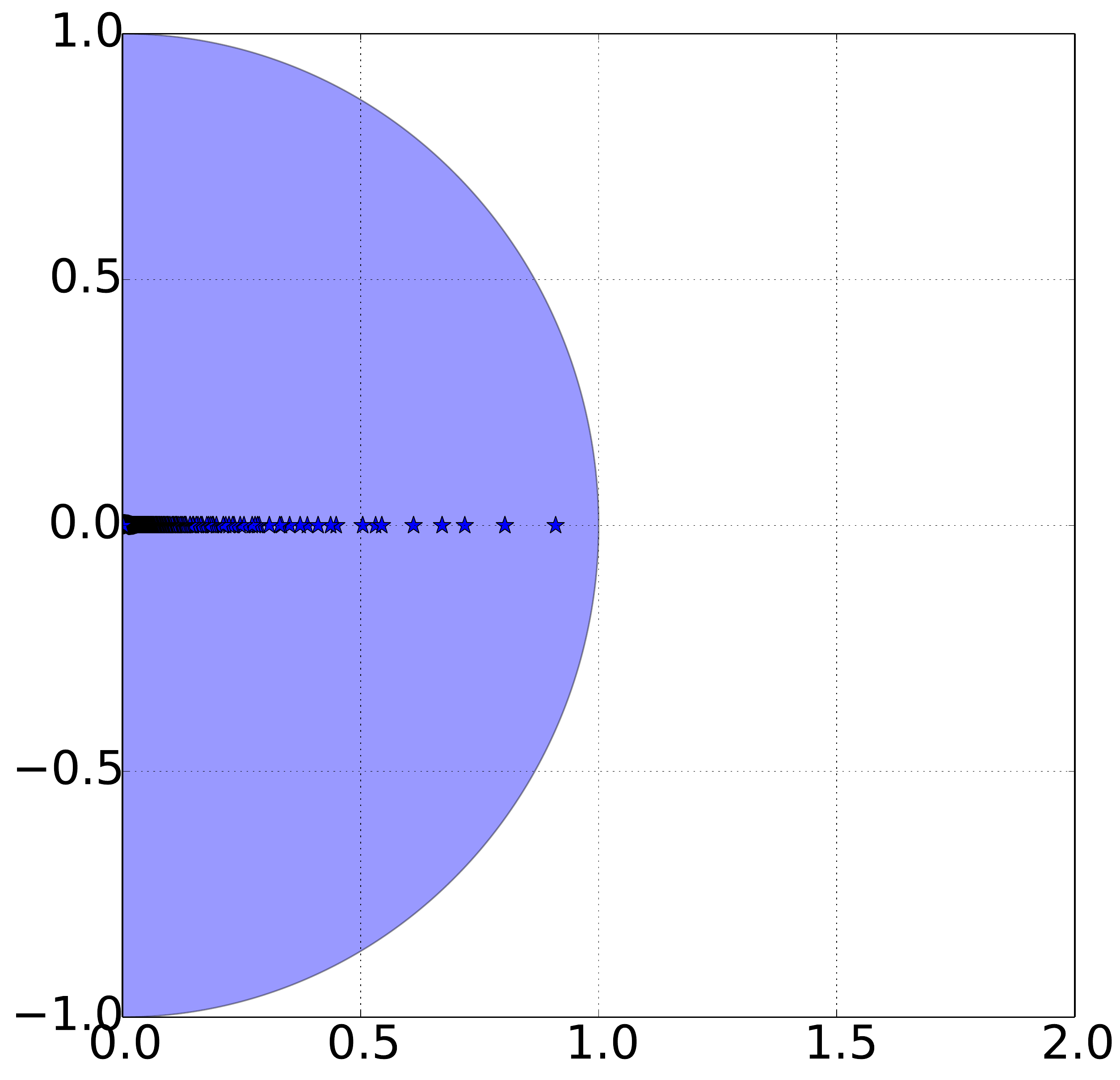}
		\caption{IMEX(1,1,1)}
  \end{subfigure}
   \begin{subfigure}[b]{0.25\textwidth}
    \includegraphics[width=\textwidth]{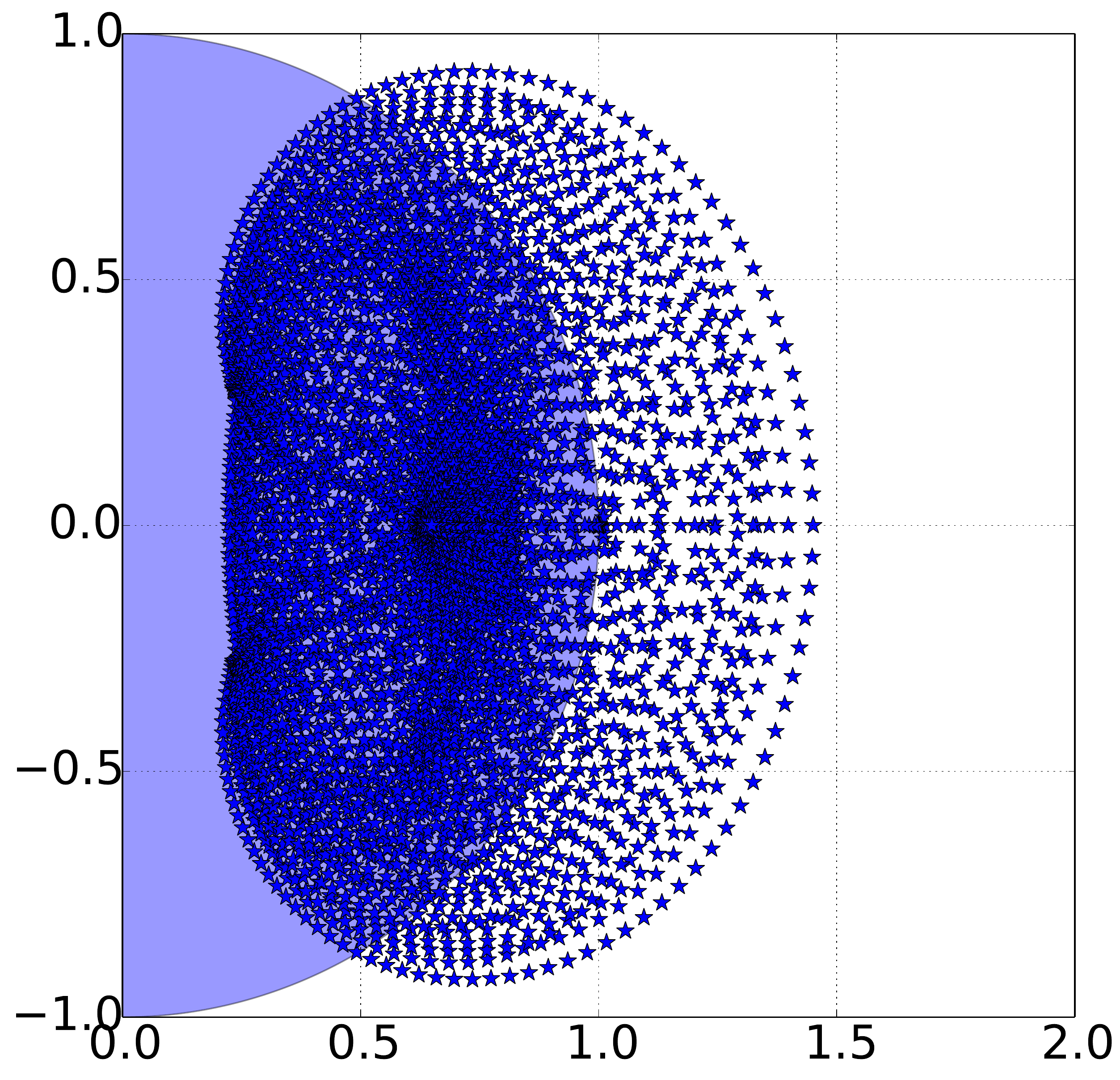}
		\caption{IMEX(1,2,1)}
  \end{subfigure}
			\caption{Eigenvalues of propagation operator for two first-order IMEX-RK
			schemes from \cite{ascher1997implicit}, with time-step $ h =0.004876$ given
			by the advective stability limit. The unit circle is shaded; to ensure stability, eigenvalues
			must be inside of shaded region.}
      \label{fig:eig}
\end{figure}

First note that stability of IMEX integration for
this discretization is challenging. Consider two variations in first-order IMEX-RK
presented in \cite{ascher1997implicit}, the IMEX(1,1,1) and IMEX(1,2,1)
schemes.  Note that a necessary condition for stability is having the eigenvalues of the propagation operator bounded $<1$ in magnitude. Thus, we consider mesh spacing $h_x\approx 0.03$ and directly construct the propagation operators for the IMEX(1,1,1)-
and IMEX(1,2,1)-schemes. Letting $\epsilon = 1$ and $ h =0.004876$
be the (numerically determined) forward Euler advective stability
limit, eigenvalues for the two propagation operators are shown in
\Cref{fig:eig}. Note that for this moderately small time step,
RK(1,2,1) is very unstable; numerical tests confirm this in practice,
with the solution rapidly growing by orders of magnitude.
In fact, many of the IMEX-RK schemes proposed in \cite{ascher1997implicit}
suffer from similar instabilities with this problem
and discretization. The \emph{only} schemes from \cite{ascher1997implicit}
that have demonstrated reasonable stability here are schemes in which
\emph{both} the implicit and explicit method are stiffly accurate.
Moving forward, we only consider these schemes, namely IMEX(1,1,1),
RK(2,2,2), and IMEX(4,4,3) in the notation of \cite{ascher1997implicit}.
In addition, we consider ESDIRK ARK3(2)4L[2]SA from \cite{kennedy2003additive}, which we abbreviate ARK(4,3), and demonstrates stability for sufficiently small $ h $. 

\begin{figure}[!htb]
  \centering
	\hspace*{1em} \includegraphics[width=0.8\textwidth]{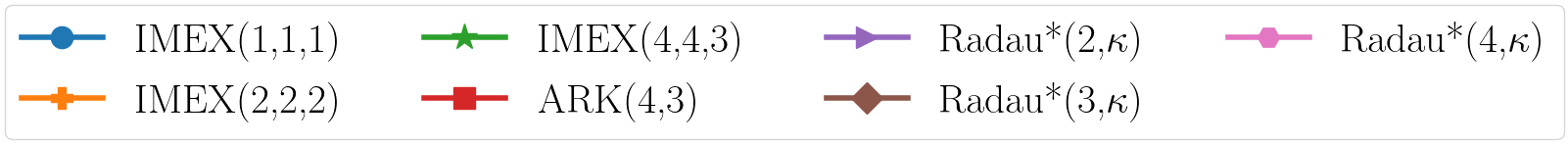}\\
	\hspace*{1em} \includegraphics[width=0.48\textwidth]{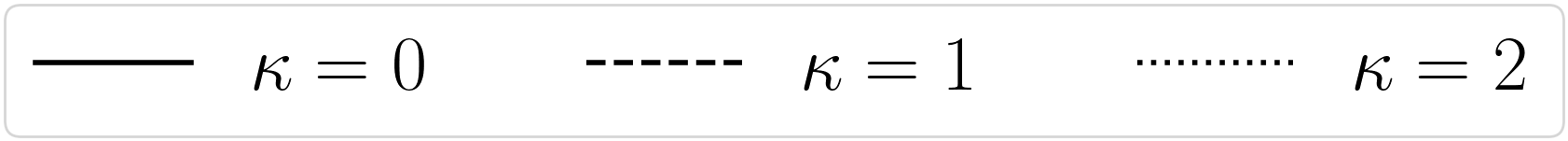}\\
  \begin{subfigure}[b]{0.4\textwidth}
    \hspace*{.001\textwidth}\includegraphics[width=.969\textwidth]{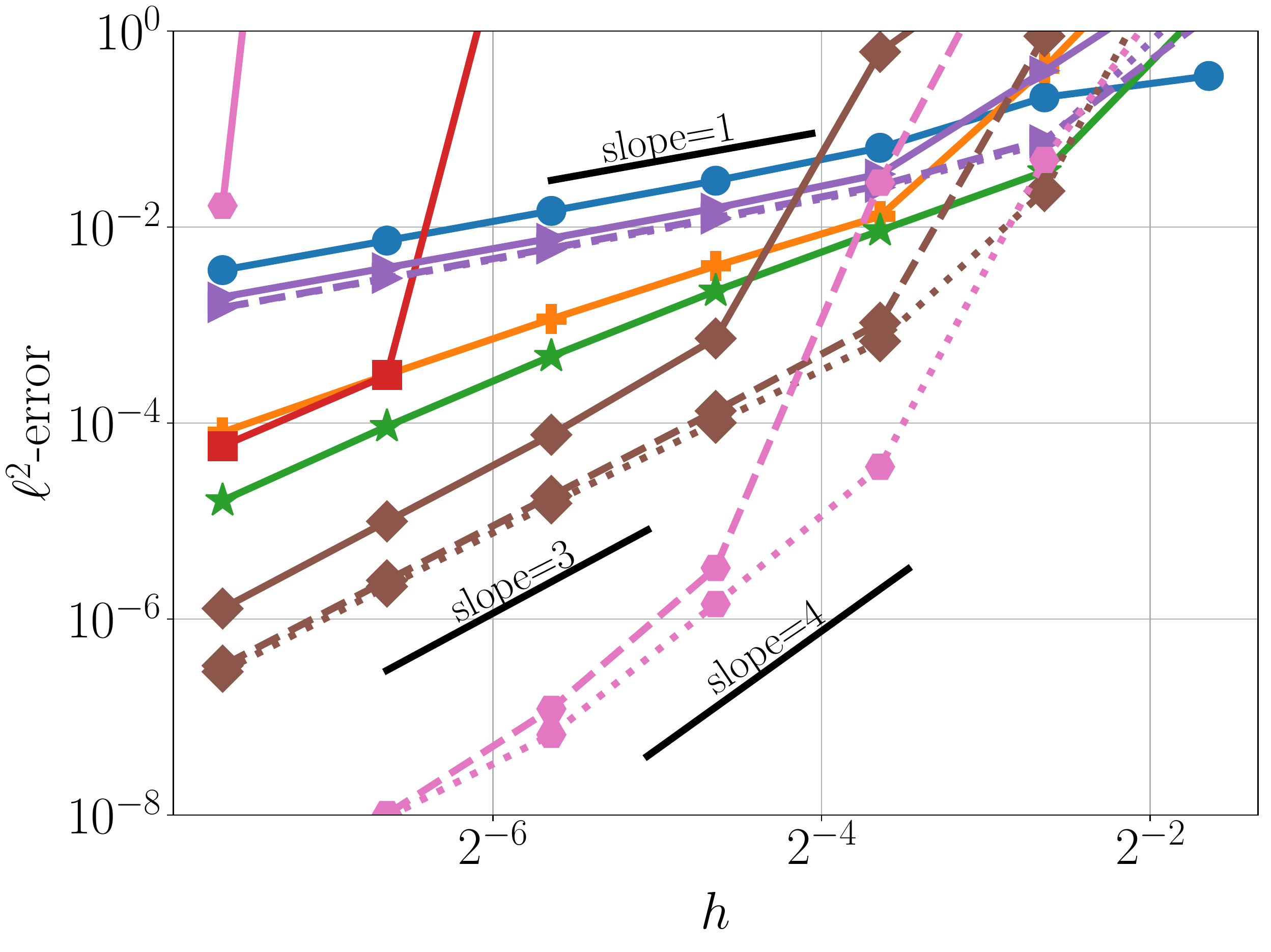} \\
    \hspace*{.01\textwidth}\includegraphics[width=.98\textwidth]{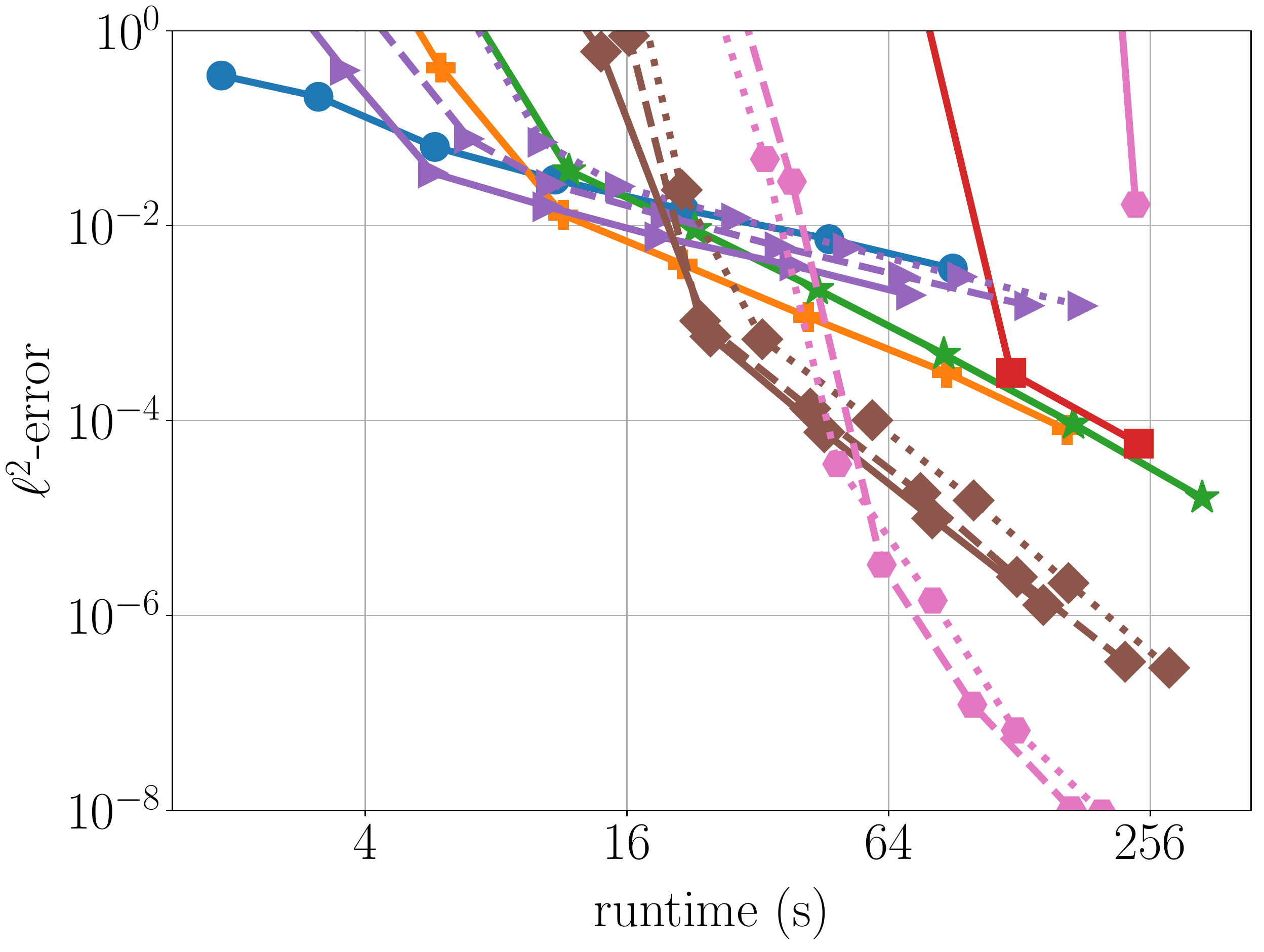}
    \caption{$\epsilon = 0.1$}\label{fig:eps01}
  \end{subfigure}
   \begin{subfigure}[b]{0.4\textwidth}
    \hspace*{.001\textwidth}\includegraphics[width=.97\textwidth]{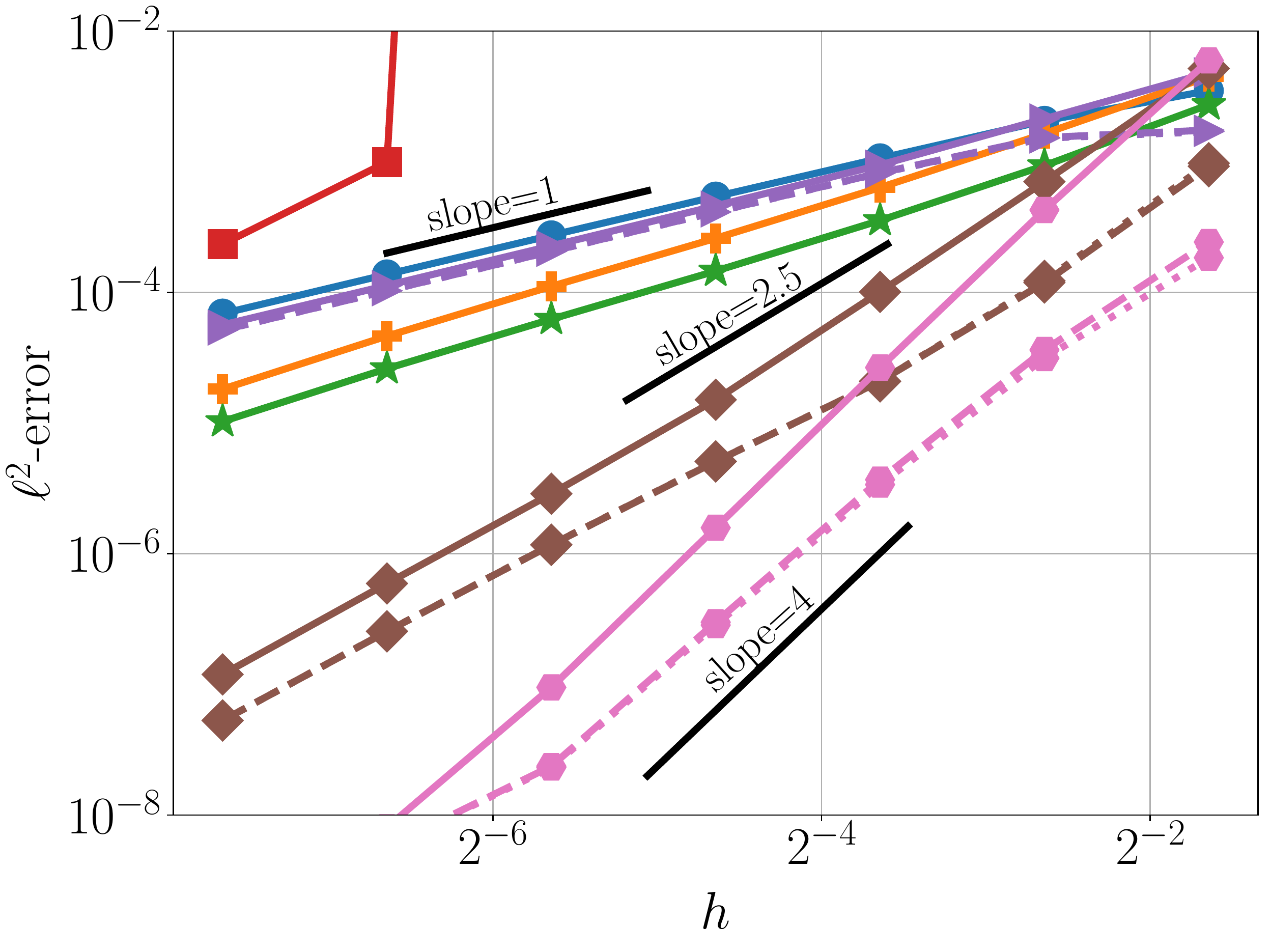}\\
    \hspace*{.01\textwidth}\includegraphics[width=.99\textwidth]{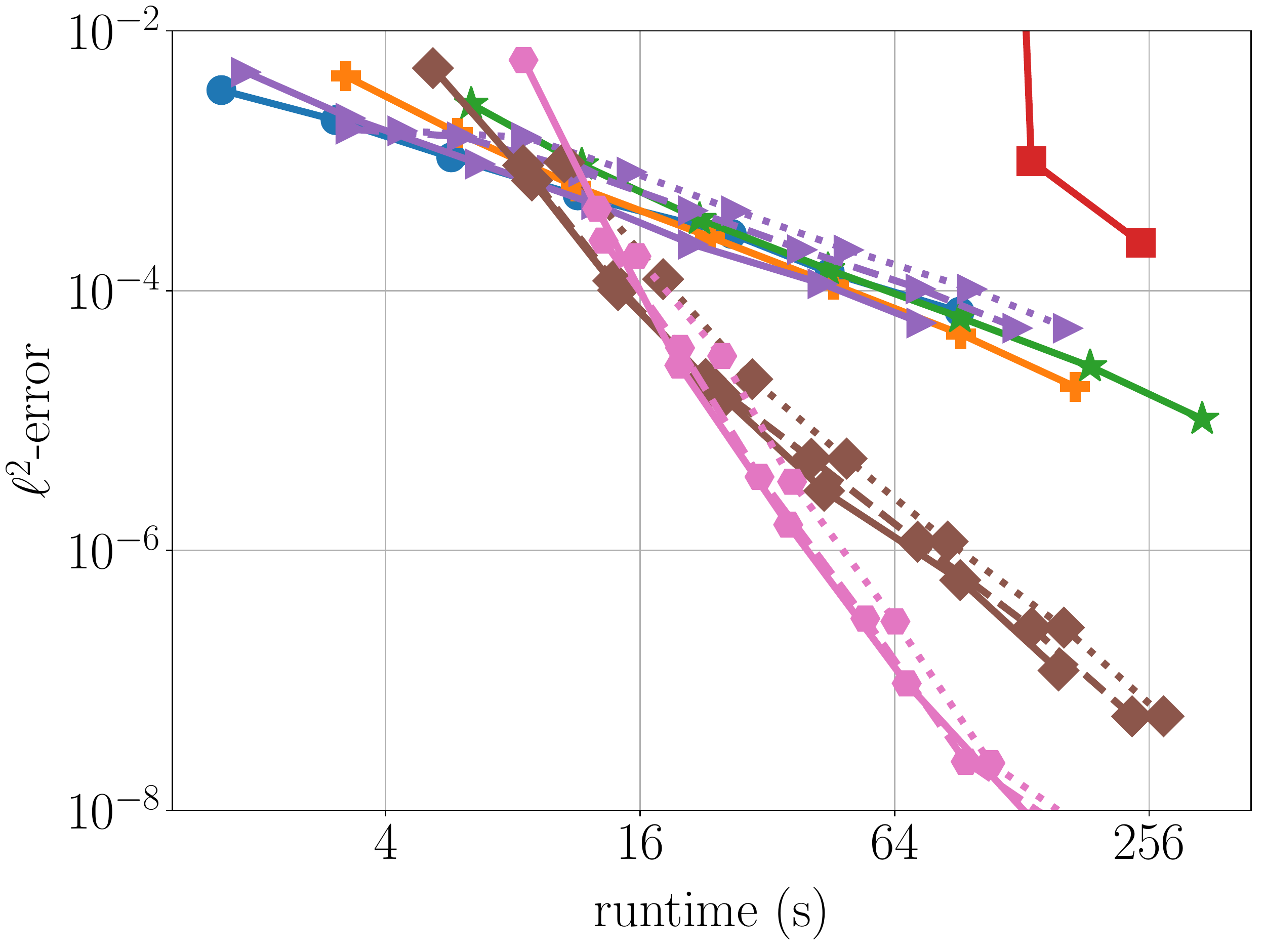}
    \caption{$\epsilon = 10$}\label{fig:eps10}
  \end{subfigure}
      \caption{$\ell^2$-accuracy as a function of $ h $ and
      wallclock time to solution.}
      \label{fig:dg}
\end{figure}

\Cref{fig:dg} compares IMEX(1,1,1), IMEX(2,2,2), IMEX(4,4,3)
\cite{ascher1997implicit}, and ARK(4,3) \cite{kennedy2003additive}, with
FIMEX-Radau*($q,\kappa$), for $q\in\{2,3,4\}$ and $\kappa\in
\{0,1,2\}$. Plots present solution accuracy with respect to time-step
size, $ h $,  and total wallclock time to solution with respect to
accuracy, for two different diffusion coefficients.
FIMEX-Radau methods are the most efficient in terms of error as a
function of runtime to solution for accuracy $\lessapprox 10^{-3}$.
Moreover, for $\epsilon=0.1$ and $\epsilon =10$, the FIMEX-Radau*
methods obtain $\sim10^{-5}$ accuracy respectively $\approx16\times$ and
$\approx16\times$ faster than the best IMEX-RK method, IMEX(4,4,3).
In terms of order of convergence, we see that all methods suffer
from some order reduction, particularly for $\epsilon = 10$, but
FIMEX-Radau* methods retain much smaller error constants than the
IMEX-RK methods.

In terms of wallclock time, each iterator application tends to
decrease error roughly proportional with its computational cost.
However, they also provide improved stability. E.g., for 
$\epsilon=0.1$ and FIMEX-Radau*(4,$\kappa$), the
method largely diverges for $\kappa =0$ (except at the smallest
$ h $ considered, $ h =0.005$), but applying one or two
applications of the iterator yields a stable, high-order 
method. Now consider $\epsilon=0.01$ and $ h  = 0.01$. Results
for the \emph{only} IMEX-RK methods that did not diverge are shown
in \Cref{tab:radau}, along with a selection of FIMEX-Radau and
FIMEX-Radau* methods that converge. Note, the one exception is
FIMEX-Radau*(2,0), which again provides improved accuracy over
IMEX-Euler, at a relatively trivial cost, and FIMEX-Radau(2,1),
which improves the accuracy of IMEX-Euler via the application of
an iterator. Note that the FIMEX-Radau* methods again provide
several orders of magnitude smaller error than IMEX(4,4,3), for
a comparable wallclock time.

\begin{table}[!ht]
  \centering
    \begin{tabular}{c c | c | c}
		\hline
		& IMEX(1,1,1) & Radau*(2,0) & Radau(2,1) \\\hline
		$\ell^2$-error & 0.018 & 0.009 & 0.009 \\
		runtime (seconds) & 35 & 37 & 71 \\\hline
	\end{tabular}
	\begin{tabular}{c c | c | c | c}
		\hline
		& IMEX(4,4,3) & Radau(3,1) & Radau*(3,2) & Radau*(4,2) \\\hline
		$\ell^2$-error & $7.4\cdot10^{-5}$ & $2.1\cdot10^{-5}$ & $2.6\cdot10^{-6}$ & $2.9\cdot10^{-8}$ \\
		runtime (seconds) & 165 & 122 & 151 & 184 \\\hline
	\end{tabular}
	\caption{$\ell^2$-error and runtime for all  stable integrators
	 with $\epsilon=0.01$ and $ h  = 0.01$. }
	\label{tab:radau}
\end{table}

\subsubsection{Advection-diffusion-reaction}\label{sec:results:dg:adr}

Here we consider a more practical splitting of explicit treatment for
a nonlinear reaction with $\gamma=10$, and implicit treatment for the
advection-diffusion equation, with $\epsilon = 0.1$. By treating the
reaction explicitly, we do not have to rebuild the implicit AMG solver
every time step or nonlinear iteration, which for sparse algebraic
solvers is often a nontrivial expense. We also compare with classical
A-stable and L-stable implicit DIRK methods, which use a simplified
Newton iteration linearized about the beginning of each time step
(we found this to be most efficient on average, compared with a
full Newton iteration which requires rebuilding the solver every
iteration, or a Picard iteration that lags the nonlinearity). DIRK
methods are prefixed with the stability (A or L) and suffixed with
the order; e.g., LDIRK3 means a 3rd-order L-stable DIRK method. The
methods are used as implemented in the MFEM library \cite{Anderson2020}.
The forcing function is chosen for an exact solution $u_*(x,y,t) = 1.0 /
(e^{10(x + y - t)} + 1)$ for $x,y,t\in[0,1]$.

For a given time integration scheme, we solve the implicit equations
to tolerance $10^{-3}$ times the order of expected spatial and
temporal accuracy based on $h_x$ and $h$ to the appropriate power.
This is to ensure systems are solved sufficiently accurately to
achieve discretization error, but avoid oversolving (e.g.,
achieving ten digits of residual accuracy for two digits of
physical accuracy). All simulations are run with 64
total MPI processes. The implicit simplified Newton iteration
typically requires 4-10 iterations to converge, where each
iteration inverts the (simplified) Jacobian to three digits
relative residual tolerance. 

\begin{figure}[!htb]
  \centering
  \begin{subfigure}[b]{0.4\textwidth}
    \hspace*{.001\textwidth}\includegraphics[width=.969\textwidth]{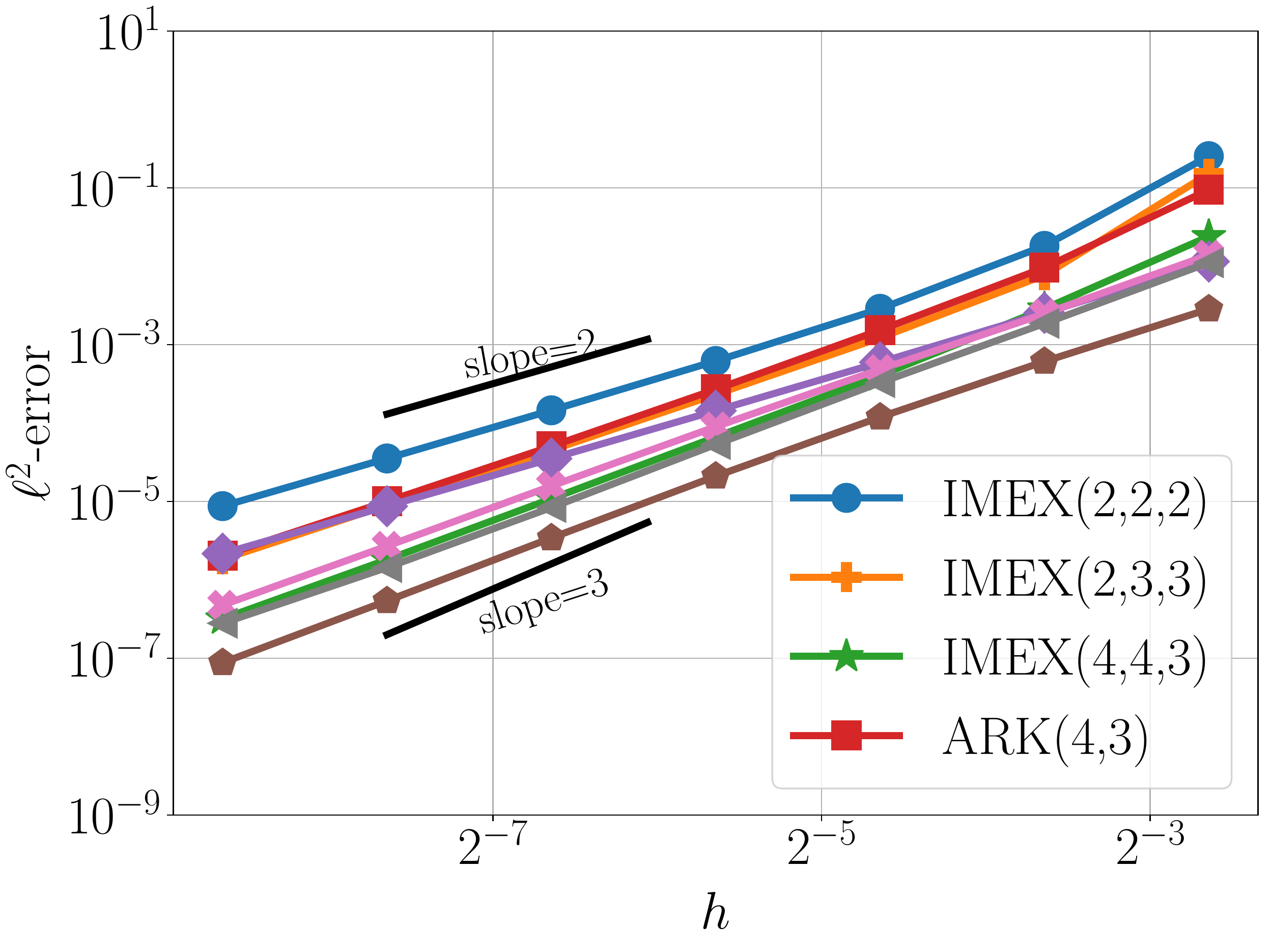}
    \hspace*{.001\textwidth}\includegraphics[width=.97\textwidth]{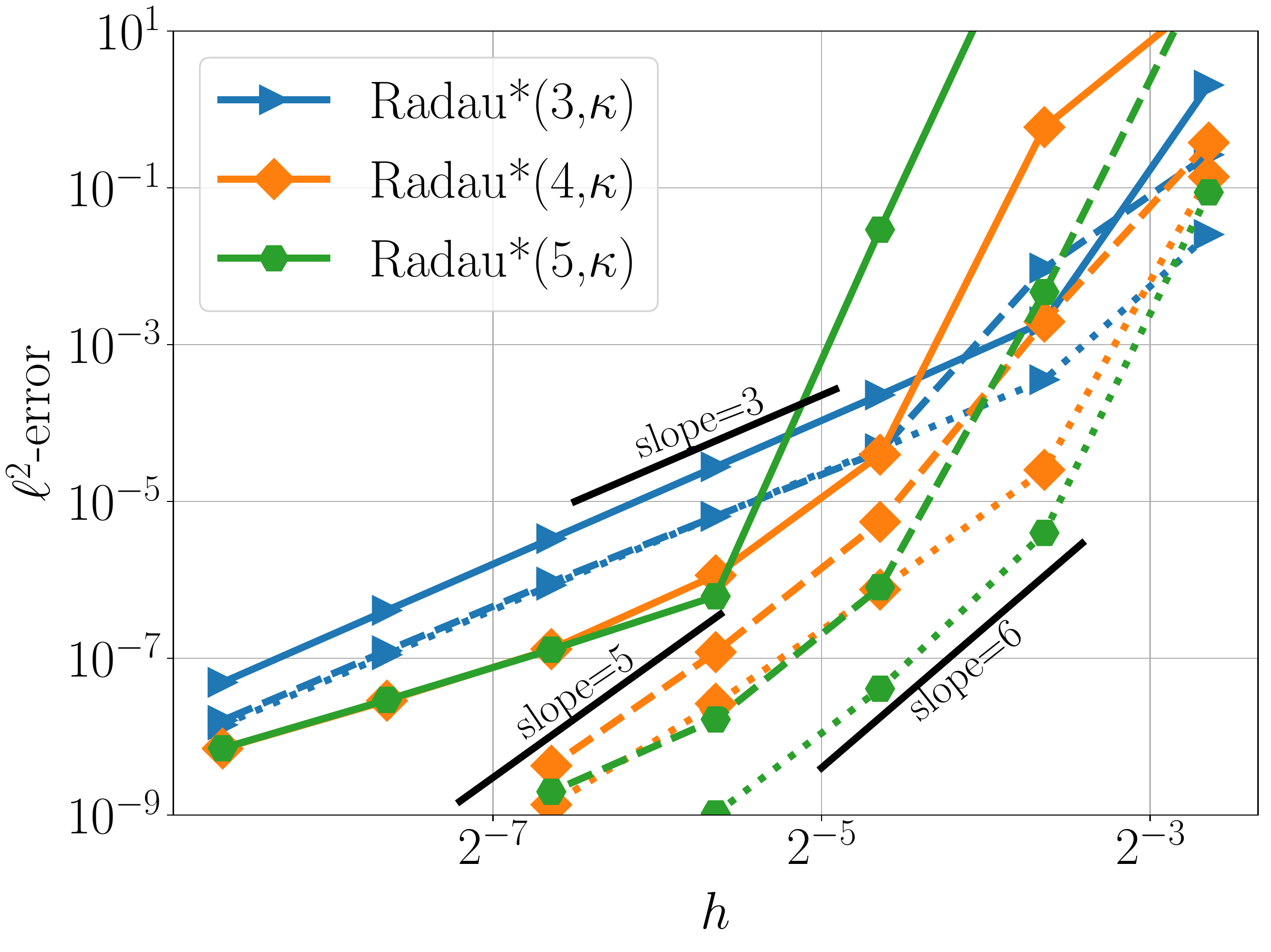}
    \caption{}\label{fig:dg-radau1}
  \end{subfigure}
   \begin{subfigure}[b]{0.4\textwidth}
    \hspace*{.01\textwidth}\includegraphics[width=.98\textwidth]{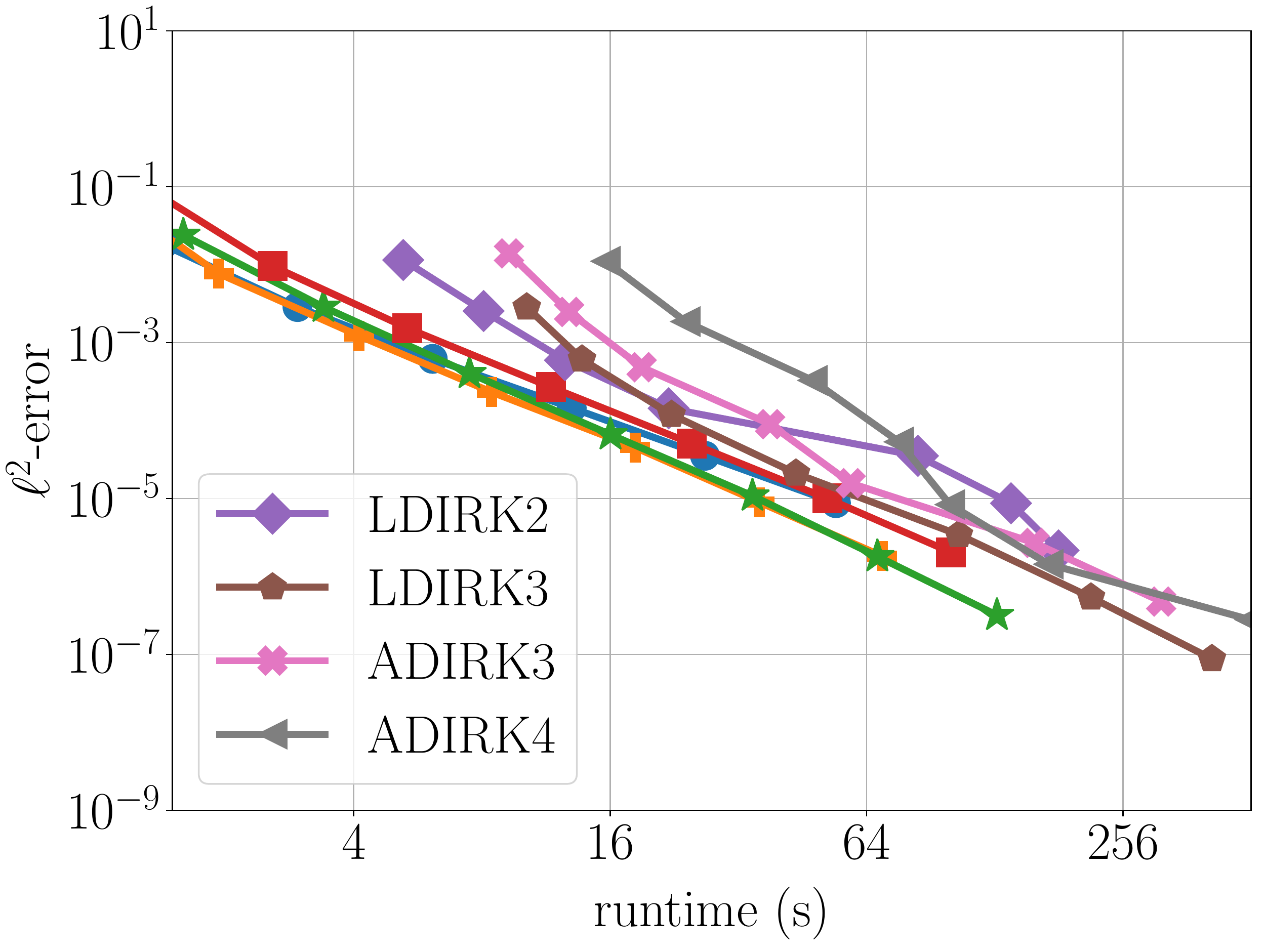}
    \hspace*{.01\textwidth}\includegraphics[width=.99\textwidth]{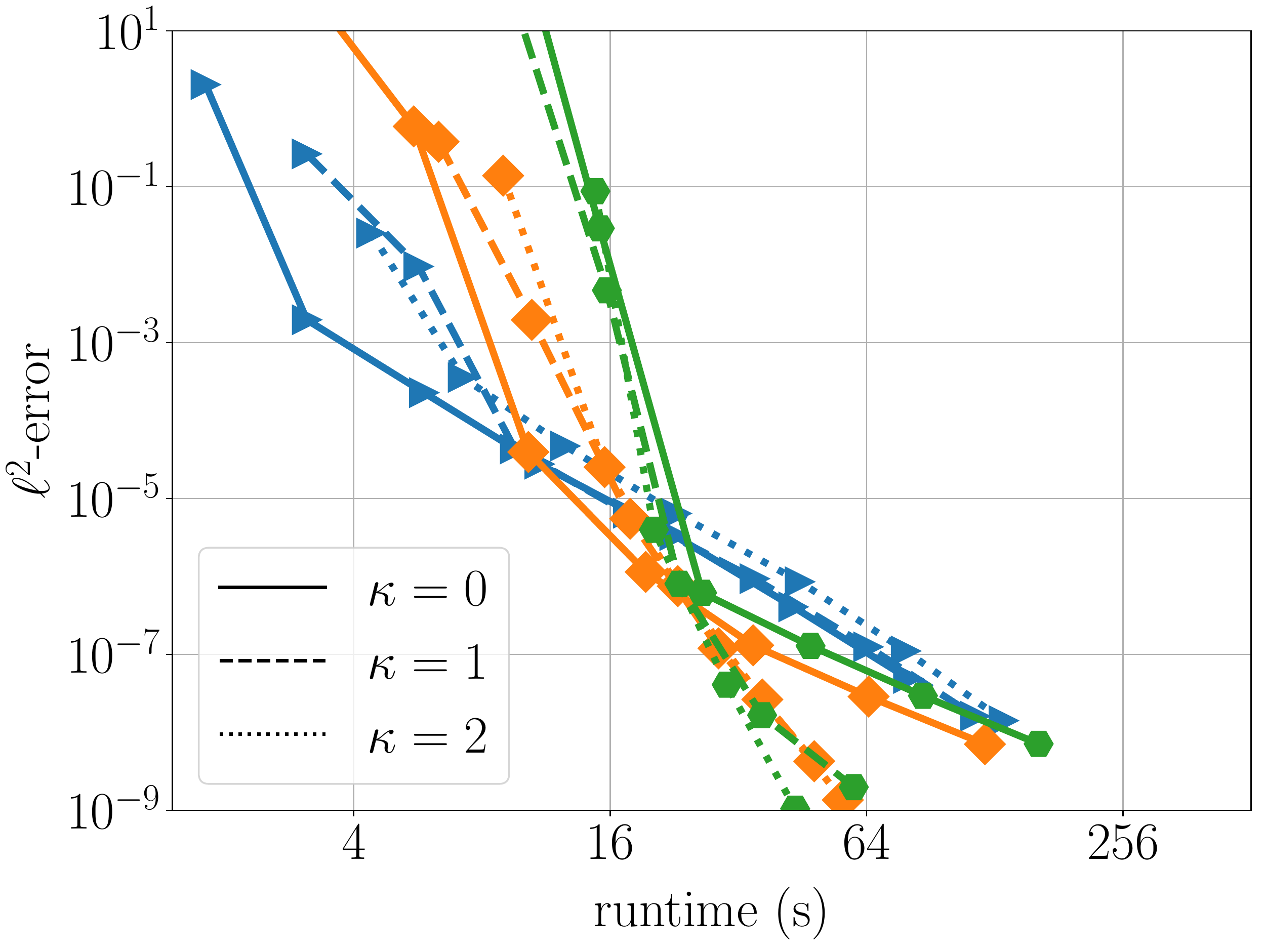}
    \caption{}\label{fig:dg-radau2}
  \end{subfigure}
      \caption{$\ell^2$-accuracy as a function of $h$ and
      wallclock time to solution. Legends apply to both images
      in a row.}
      \label{fig:dg-adr}
\end{figure}

Error as a function of $h$ and total wallclock time are shown
in \Cref{fig:dg-adr}. Here we see that most schemes achieve
their expected accuracy (recall the order of Radau$^*(q,\kappa)
= q+\kappa$, up to the order of Radau, given by $2q-3$).
Radau$^*(q,0)$ is less predictable (e.g., see $q=4,5$ for $\kappa=0$),
but the error constants remain small once a reasonably small step
size is chosen. It is possible the order reduction is due to time
dependent boundary conditions, a known problem with RK methods,
as we can also see that IMEX-RK and DIRK methods don't quite
reach expected 3rd-order converence. In general, we see that
IMEX schemes are advantageous for this problem over DIRK methods,
typically achieving better accuracy for a fixed wallclock
time. The difference is small for classical IMEX-RK methods,
but more noticable with FIMEX-Radau$^*$ methods.
For moderate to high accuracies, roughly lower than
$10^{-4}$, the FIMEX-Radau$^*$ methods are the best solution
with respect to wallclock time. It is also worth pointing
out that the asymptotic error constants for FIMEX-Radau$^*$ are
quite small; e.g., for 3rd-order methods with a fixed $h$,
FIMEX-Radau$^*$ error is $10-100\times$ smaller than
comparable-order IMEX-RK schemes.

\section{Summary and Conclusion}
\label{sec:conclusion}

In the first half of this paper we introduced a general framework for constructing multivalued or one-step additive polynomial time integrators with any degree of implicitness. By utilizing interpolating polynomials, the framework enables the derivation of high-order methods without requiring the solution of nonlinear order conditions (even if the underlying equation has an arbitrary number of additive partitions). Moreover, by selecting high-order \ODEps{}, one trivially achieves high stage order that prevents order reduction on stiff problems. Lastly, the framework also includes iterators that can be used to compute initial conditions and to construct composite PBMs with improved stability and accuracy properties.

The second-half of this paper focused exclusively on FIMEX methods that utilize Adams \ODEps{}. Specifically we introduced a new class of FIMEX-Radau integrators, which use a Radau IIA fully implicit component coupled with a parallel block method for the explicit component. In Section \ref{sec:numerical-experiments} we demonstrated the significant potential of these integrators in practice. In each of our experiments, FIMEX-Radau methods consistently yielded the most accurate solution compared with other methods from \cite{ascher1997implicit,kennedy2003additive}. Moreover, thanks to recent developments in solvers for fully implicit Runge-Kutta methods \cite{irk1,irk2}, they were also the fastest methods (in terms of total runtime) to obtain a given accuracy. Furthermore, by using iterators, the composite FIMEX-Radau methods offer a simple way to develop very high-order FIMEX integrators with low storage cost. We believe these methods are most applicable for problems where moderate to high accuracy is desired, and one can benefit from the strong accuracy and stability provided by combining multistep methods and fully implicit methods (e.g., for time-dependent boundary conditions, or semi-stiff nonlinearities that one does not want to treat implicitly). 

Thus far, we have only introduced one class of additive polynomial integrators. In part two of this work we will introduce new families of diagonally implicit IMEX-Polynomial integrators that are based on the BDF \ODEps{}. In addition, we point out that Radau and other fully implicit Runge-Kutta methods offer a number of unique advantages and properties that cannot be obtained with other classes of integrators. Future work will study the potential of FIMEX-Radau and related integrators on topics such as differential algebraic equations and conservation of invariants which fully implicit Runge-Kutta methods are uniquely suited for. We also plan to develop strategies for adaptively controlling the stepsize and the number of iterations $\kappa$ for polynomial fully-implicit-explicit collocation methods.

\section*{Acknowledgements} Los Alamos National Laboratory report number LA-UR-21-28709, and the National Science Foundation, Computational Mathematics Program DMS-2012875.
	
\bibliographystyle{siamplain}
\bibliography{references_nourl}

\appendix

\end{document}


\maketitle

\section{FIMEX-Radau and FIMEX-Radau*}

This section contains additional information about the polynomial FIMEX-Radau and FIMEX-Radau* methods including: coefficients for methods with $q=3$ and $q=4$, links to a code repository that contains scripts for deriving the coefficients for higher-order methods, and an additional stability figure for FIMEX-Radau that is not contained in the paper. 

\subsection{Coefficients for $q=2,3,4$}
\label{sup:imex-radau-coefficients}

FIMEX-Radau and FIMEX-Radau$^*$ propagator and iterator methods can be expressed as
\begin{align}
	& \text{propagator:}
	& \mathbf{y}^{[n+1]} &= \mathbf{A} \mathbf{y}^{[n]} + r \mathbf{B}^{\{1\}} \mathbf{f}^{\{1\}[n+1]} + r \mathbf{B}^{\{2\}} \mathbf{f}^{\{2\}[n]} \\
	& \text{iterator:}
	& \mathbf{y}^{[n+1]} &= \tilde{\mathbf{A}} \mathbf{y}^{[n]} + r \mathbf{B}^{\{1\}} \left( \mathbf{f}^{\{1\}[n+1]} + \mathbf{f}^{\{2\}[n]}\right).
\end{align}
where the $q\times q$ matrices $\mathbf{A}$, $\tilde{\mathbf{A}}$ are
\begin{align}
	A_{ij} = 
	\begin{cases}
 		1 & $j = q$ \\
 		0 & \text{otherwise}
 	\end{cases}	&&
 	\tilde{A}_{ij} = 
	\begin{cases}
 		1 & $j = 1$ \\
 		0 & \text{otherwise}
 	\end{cases}	
\end{align}
and the $q\times q$ matrices $\mathbf{B}^{\{1\}}$, and $\mathbf{B}^{\{2\}}$ for $q=3$, $q=4$ are
\begin{itemize}[leftmargin=*]
	\item $q=2$ (for FIMEX-Radau$^*$ replace $\mathbf{B}^{\{2\}}$ with $\mathbf{B}^{*\{2\}}$)
		\begin{small}
			\begin{align*}
				\mathbf{B}^{\{1\}} &=
				\left[
					\begin{array}{ccc}
					 0 & 0 \\
					 0 & 2 \\
					\end{array}
				\right]	 &
				\mathbf{B}^{\{2\}} &=
				\left[
					\begin{array}{ccc}
					 0 & 0 \\
					 0 & 2
				\end{array}
				\right] &
				\mathbf{B}^{*\{2\}} &=
				\left[
					\begin{array}{ccc}
					 0 & 0 \\
					 -1 & 3
					\end{array}
				\right]
			\end{align*}	
		\end{small}
	\item $q=3$ (for FIMEX-Radau$^*$ replace $\mathbf{B}^{\{2\}}$ with $\mathbf{B}^{*\{2\}}$)
		\begin{small}
			\begin{align*}
				\mathbf{B}^{\{1\}} &=
				\left[
					\begin{array}{ccc}
					 0 & 0 & 0 \\
					 0 & \frac{5}{6} & -\frac{1}{6} \\
					 0 & \frac{3}{2} & \frac{1}{2} \\
					\end{array}
				\right]	 &
				\mathbf{B}^{\{2\}} &=
				\left[
					\begin{array}{ccc}
					 0 & 0 & 0 \\
					 0 & -\frac{1}{6} & \frac{5}{6} \\
					 0 & -\frac{3}{2} & \frac{7}{2} \\
					\end{array}
				\right] &
				\mathbf{B}^{*\{2\}} &=
				\left[
					\begin{array}{ccc}
					 0 & 0 & 0 \\
					 \frac{8}{27} & -\frac{11}{18} & \frac{53}{54} \\
					 4 & -\frac{15}{2} & \frac{11}{2} \\
					\end{array}
				\right]
			\end{align*}	
		\end{small}
	\item $q=4$	(for FIMEX-Radau$^*$ replace $\mathbf{B}^{\{2\}}$ with $\mathbf{B}^{*\{2\}}$)
	\begin{small}
		\begin{align*}
			\mathbf{B}^{\{1\}} &=\left[
			\begin{array}{cccc}
			 0 & 0 & 0 & 0 \\
			 0 & \frac{1}{180} \left(88-7 \sqrt{6}\right) & \frac{1}{900} \left(296-169 \sqrt{6}\right) & \frac{2}{225} \left(-2+3 \sqrt{6}\right) \\
			 0 & \frac{1}{900} \left(296+169 \sqrt{6}\right) & \frac{1}{180} \left(88+7 \sqrt{6}\right) & \frac{1}{225} (-2) \left(2+3 \sqrt{6}\right) \\
			 0 & \frac{8}{9}-\frac{1}{3 \sqrt{6}} & \frac{1}{18} \left(16+\sqrt{6}\right) & \frac{2}{9} \\
			\end{array} 
			\right] \\
			\mathbf{B}^{\{2\}} &= \left[
			\begin{array}{cccc}
			 0 & 0 & 0 & 0 \\
			 0 & \frac{1}{180} \left(-272+113 \sqrt{6}\right) & \frac{1}{900} \left(296-169 \sqrt{6}\right) & \frac{2}{225} \left(223-72 \sqrt{6}\right) \\
			 0 & \frac{1}{900} \left(296+169 \sqrt{6}\right) & \frac{1}{180} \left(-272-113 \sqrt{6}\right) & \frac{2}{225} \left(223+72 \sqrt{6}\right) \\
			 0 & \frac{1}{18} \left(-56+41 \sqrt{6}\right) & \frac{1}{18} \left(-56-41 \sqrt{6}\right) & \frac{74}{9} \\
			\end{array}
			\right] \\
			\mathbf{B}^{*\{2\}} &=\left[
				\begin{array}{cccc}
				 0 & 0 & 0 & 0 \\
				 \frac{1}{500} \left(-1091+424 \sqrt{6}\right) & \frac{1}{360} \left(-1198+517 \sqrt{6}\right) & \frac{32402-14053 \sqrt{6}}{9000} & \frac{12193-4152 \sqrt{6}}{4500} \\
				 \frac{1}{500} \left(-1091-424 \sqrt{6}\right) & \frac{32402+14053 \sqrt{6}}{9000} & \frac{1}{360} \left(-1198-517 \sqrt{6}\right) & \frac{12193+4152 \sqrt{6}}{4500} \\
				 -16 & \frac{5}{18} \left(8+37 \sqrt{6}\right) & \frac{1}{18} (-5) \left(-8+37 \sqrt{6}\right) & \frac{122}{9} \\
				\end{array}
			\right]
		\end{align*}
	\end{small} 
\end{itemize}
For larger $q$ the coefficients can be obtained using the MATLAB script that found in the repository \url{https://github.com/pipack/paper-additive-pbms}. For convenience we also include the script in \cref{sup:matlab-radau-script}.

\subsection{Matlab Script for Radau Coefficients}
\label{sup:matlab-radau-script}

\begin{verbatim}
	function [Bi, Be] = radauCoeff(q, star)
    if(nargin == 1)
        star = false;
    end

    Bi = zeros(q);
    Be = zeros(q);
    z  = nodes(q);
    Bi(:,2:end) = transpose(quadW(z(2:end), z(1) * ones(q,1), z));
    if(star)
        Be = transpose(quadW(z, z(q) * ones(q,1), z + 2));
    else
        Be(:,2:end) = transpose(quadW(z(2:end), z(q) * ones(q,1), z + 2));
    end
end

function w = quadW(x, a, b)
%QUADW Returns quadrature weights for the nodes x where
%
%   \int^b(i)_a(i) f(x) dx \approx \sum_{j=1}^n w(j,i) f(x(i))
%
% == Parameters ===============================================================
%   x (vector) - nodes
%   a (vector) - left endpoints
%   b (vector) - right endpoints
% =============================================================================

V = transpose(fliplr(vander(x)));
p = (1:length(x))';
b = ((b(:)').^p - (a(:)').^p ) ./ p;
w = V \ b;
end

function z = nodes(q)
%QUADW Returns nodes for radau method
% == Parameters ===============================================================
%   q (vector) - right endpoints
% =============================================================================

if(q > 8) % use chebfun if q > 8
    if(isempty(which('radaupts')))
        error('q > 8 requires Chebfun (https://www.chebfun.org/)');
    end
    z = [-1; -flip(radaupts(q-1))];
    return;
end

nodes = {
    [1]'
    [-0.333333333333333,  1]'
    [-0.689897948556636,  0.289897948556636, 1]'
    [-0.822824080974592, -0.181066271118531,  0.575318923521694,  1]'
    [-0.885791607770965, -0.446313972723752,  0.167180864737834,  ...
    0.720480271312439,  1]'
    [-0.920380285897062, -0.603973164252784, -0.124050379505228,  ...
    0.390928546707272,  0.802929828402347,  1]'
    [-0.941367145680430, -0.703842800663031, -0.326030619437691,  ...
    0.117343037543100,  0.538467724060109,  0.853891342639482,  1]'
}; 
z = [-1; nodes{q-1}];

end
\end{verbatim}

\section{Additional stability regions for FIMEX-Radau}
\label{sup:imex-radau-stability}

\begin{figure}[h!] 
	\centering

	\begin{minipage}{0.03\textwidth}
	\end{minipage}
	\begin{minipage}{0.31\textwidth}
		\centering
		\begin{footnotesize}
			\hspace{2em} $\kappa = 0$
		\end{footnotesize}
	\end{minipage}
	\begin{minipage}{0.31\textwidth}
		\centering
		\begin{footnotesize}
			\hspace{2em} $\kappa = 1$
		\end{footnotesize}
	\end{minipage}
	\begin{minipage}{0.31\textwidth}
		\centering
		\begin{footnotesize}
			\hspace{2em} $\kappa = 2$
		\end{footnotesize}
	\end{minipage} \\[1em]
	
	\begin{minipage}{0.03\textwidth}
		\begin{footnotesize}
			\rotatebox{90}{$\theta = 0$}
		\end{footnotesize}
	\end{minipage}
	\begin{minipage}{0.31\textwidth}
		\includegraphics[width=1\textwidth]{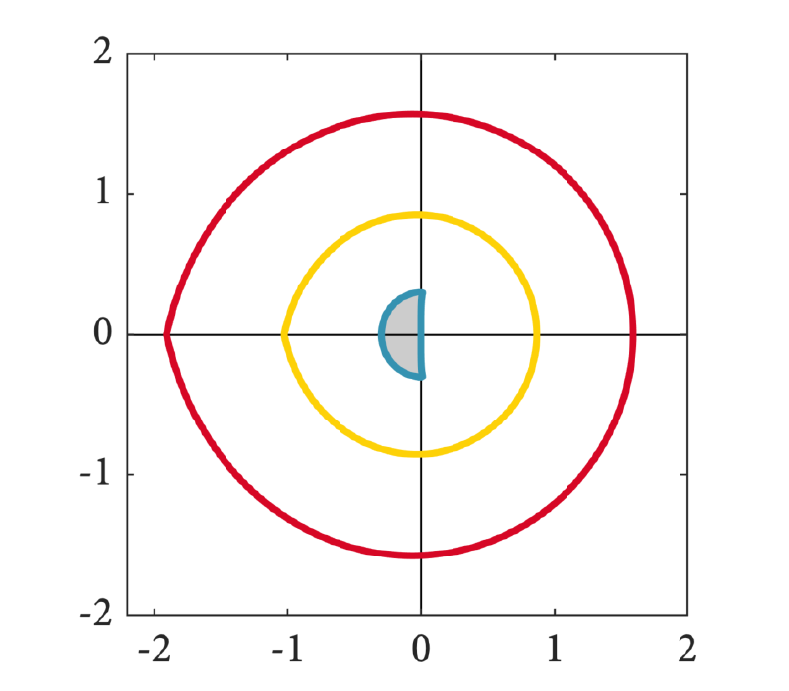}
	\end{minipage}
	\begin{minipage}{0.31\textwidth}
		\includegraphics[width=1\textwidth]{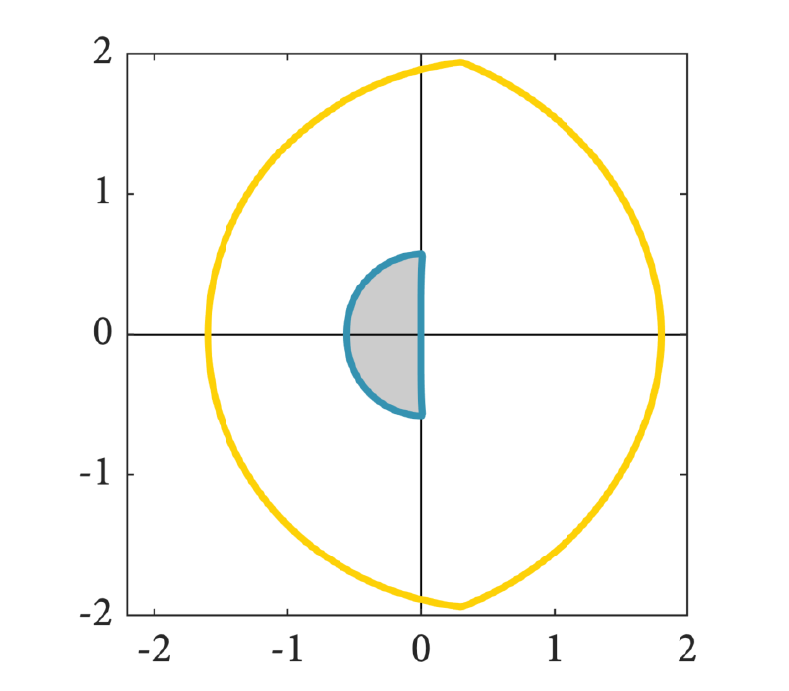}
	\end{minipage}
	\begin{minipage}{0.31\textwidth}
		\includegraphics[width=1\textwidth]{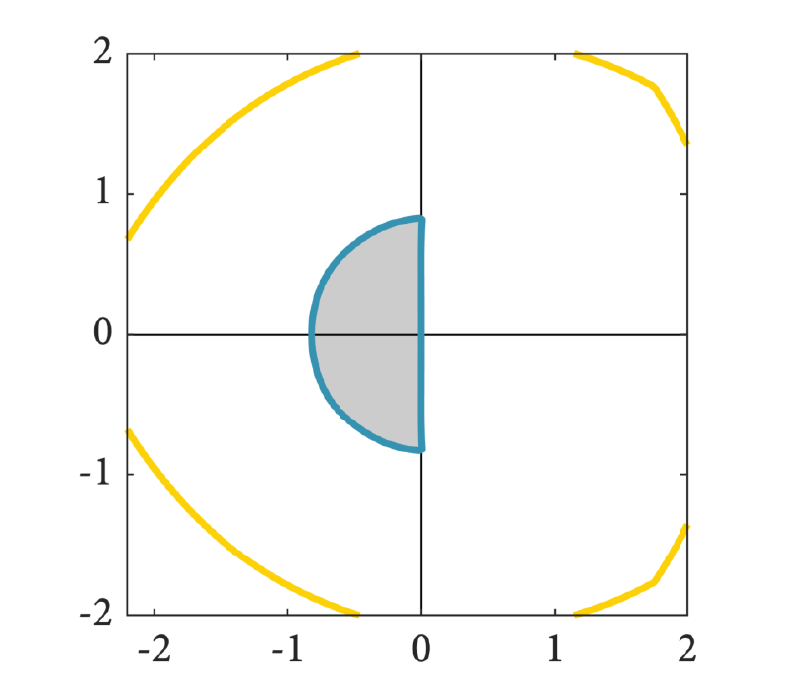}
	\end{minipage}
	
	\begin{minipage}{0.03\textwidth}
		\begin{footnotesize}
			\rotatebox{90}{$\theta = \frac{3\pi}{2}$}
		\end{footnotesize}
	\end{minipage}
	\begin{minipage}{0.31\textwidth}
		\includegraphics[width=1\textwidth]{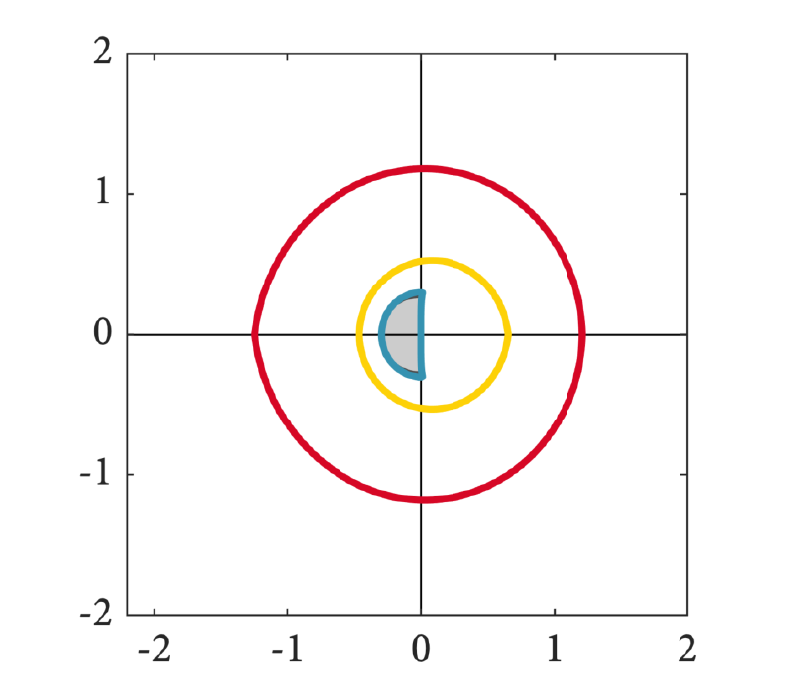}
	\end{minipage}
	\begin{minipage}{0.31\textwidth}
		\includegraphics[width=1\textwidth]{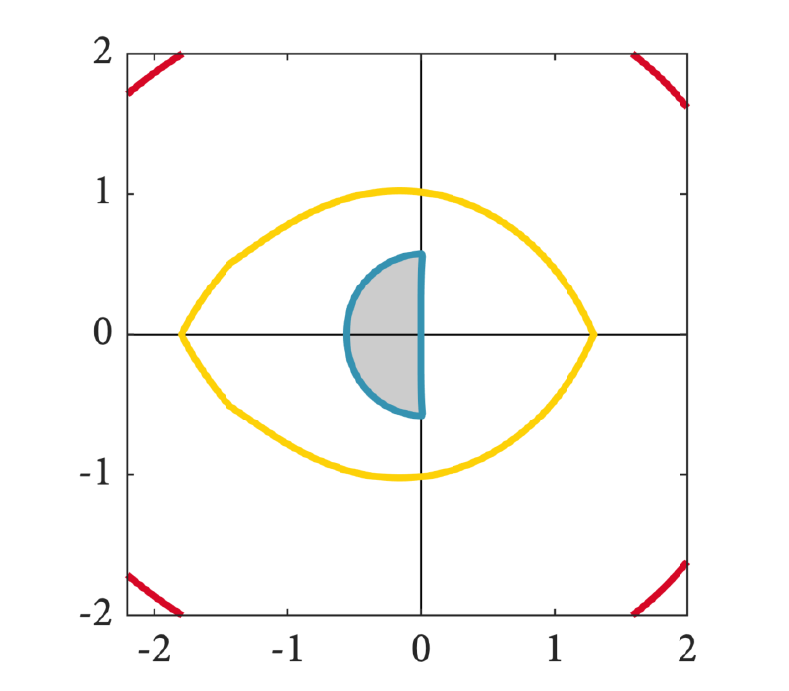}
	\end{minipage}
	\begin{minipage}{0.31\textwidth}
		\includegraphics[width=1\textwidth]{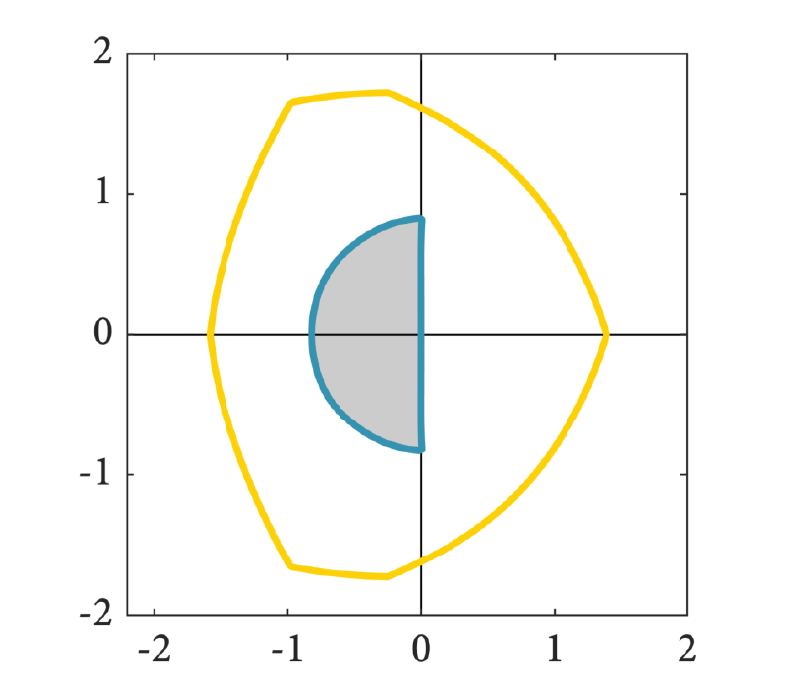}
	\end{minipage}
	
	\begin{minipage}{0.03\textwidth}
		\begin{footnotesize}
			\rotatebox{90}{$\theta = \frac{\pi}{2}$}
		\end{footnotesize}
	\end{minipage}
	\begin{minipage}{0.31\textwidth}
		\includegraphics[width=1\textwidth]{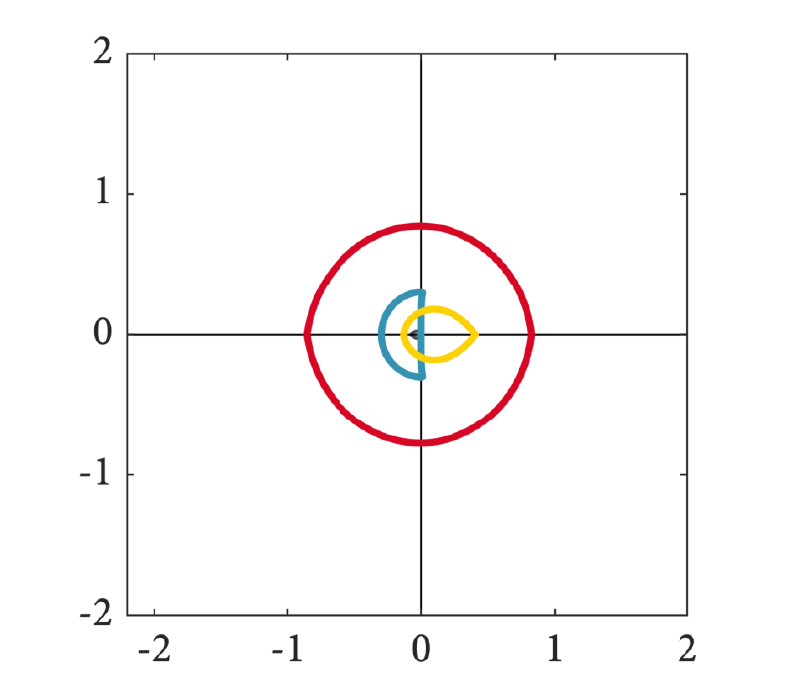}
	\end{minipage}
	\begin{minipage}{0.31\textwidth}
		\includegraphics[width=1\textwidth]{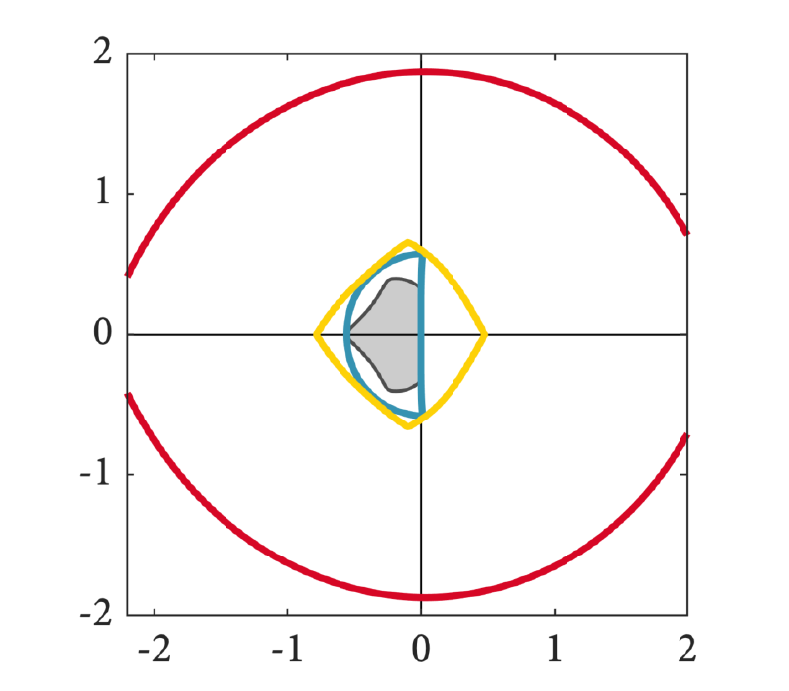}
	\end{minipage}
	\begin{minipage}{0.31\textwidth}
		\includegraphics[width=1\textwidth]{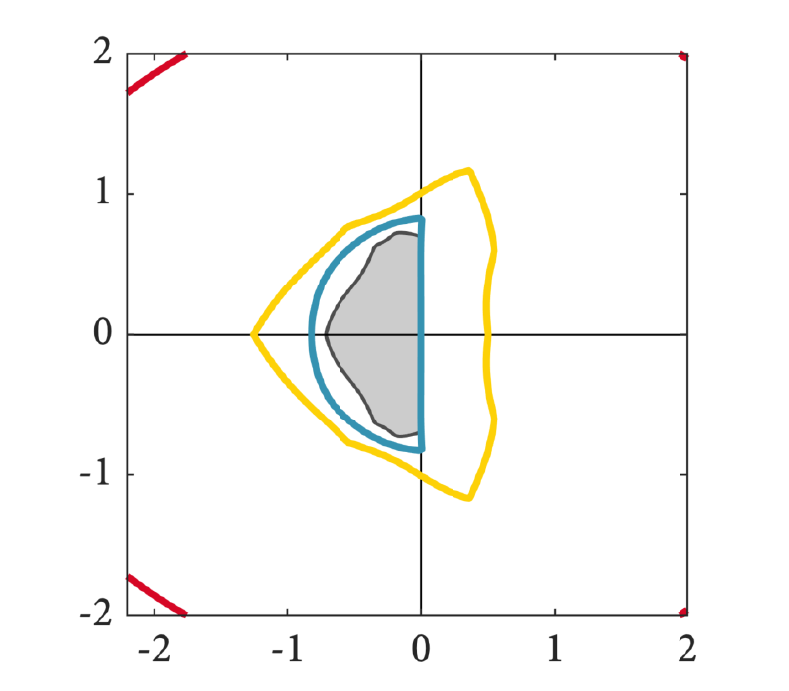}
	\end{minipage}]
	
		\begin{tabular}{c}
{\tiny \textcolor{plot_blue}{\hdashrule[0.2ex]{2em}{2pt}{}} $\left|r\right| = 0$ \hspace{1em}}
{\tiny \textcolor{plot_yellow}{\hdashrule[0.2ex]{2em}{2pt}{}} $\left|r\right| = 3$ \hspace{1em}}
{\tiny \textcolor{plot_red}{\hdashrule[0.2ex]{2em}{2pt}{}} $\left|r\right| = 6$}
\end{tabular}

	\caption{Stability regions for composite FIMEX-Radau with $q=4$. Each colored contour represents a different stability region $R_{r,\theta}$ defined in (\ref{eq:stability-region-symmetric}). The gray region is the stability region $R_\theta$ described in (\ref{eq:stability-region-symmetric-rinf}).}
	\label{supfig:pradau-figure}
\end{figure}

In Figure \ref{supfig:pradau-figure} we show stability regions for a composite FIMEX-Radau method with $q=4$. These can be be directly compared with the stability regions for the composite FIMEX-Radau* with $q=4$ that are shown in Figure \ref{fig:pradaus-figure}. In short we can see that the higher-order polynomial approximation used by FIMEX-Radau* for the explicit component leads to decreased in stability.

\section{Van Der Pol Results for IMEX-RK}
\label{sup:van-der-pol-imex-rk}

\begin{figure}[h!]
	\centering
	
	\begin{minipage}{0.48\textwidth}
		\centering
		\begin{footnotesize}
			{\bf Semi-Implicit Splitting}	
		\end{footnotesize}
		\vspace{1em}

		\begin{footnotesize}
			Approximate Convergence Rate
		\end{footnotesize}
		\includegraphics[width=1\linewidth]{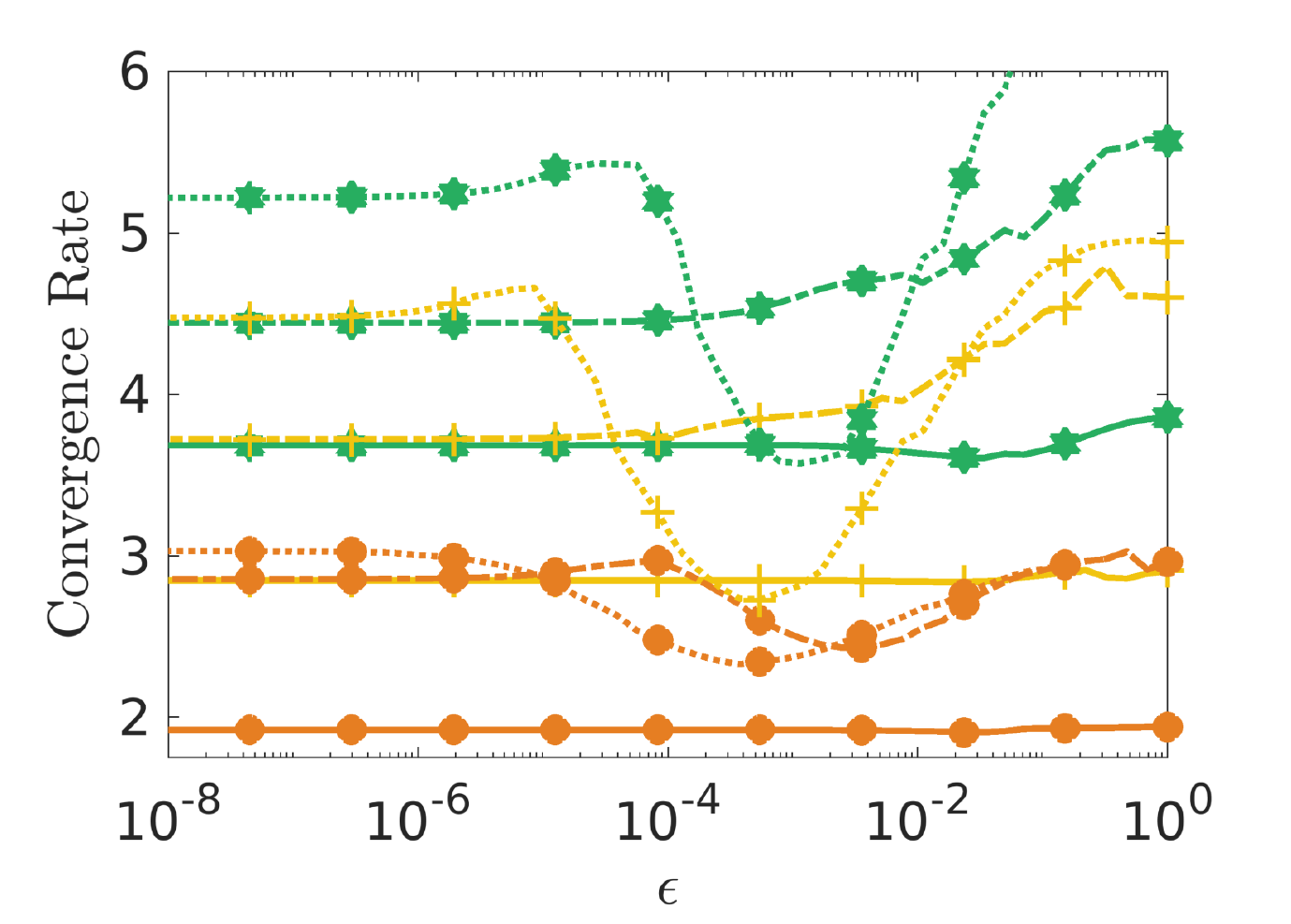}
	\end{minipage}
	\begin{minipage}{0.48\textwidth}
		\centering
		\begin{footnotesize}
			{\bf Linearly-Implicit Splitting}	
		\end{footnotesize}
		\vspace{1em}
		
		\begin{footnotesize}
			Approximate Convergence Rate
		\end{footnotesize}
		\includegraphics[width=1\linewidth]{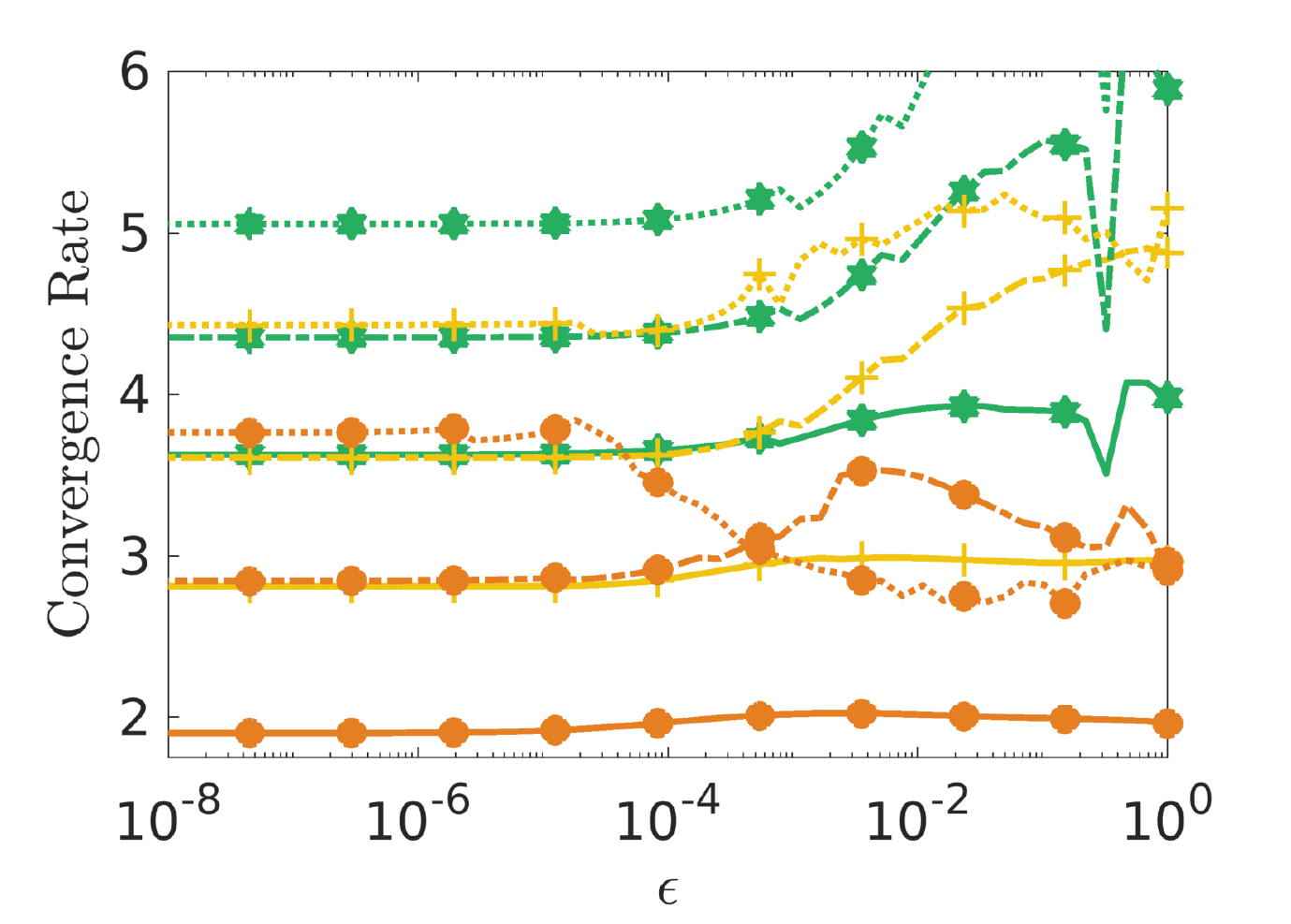}
	\end{minipage}
	
	\begin{minipage}{0.48\textwidth}
		\centering
		
		\begin{footnotesize}
			Convergence Diagram ($\epsilon = 10^{-3}$)
		\end{footnotesize}

		\includegraphics[width=1\linewidth]{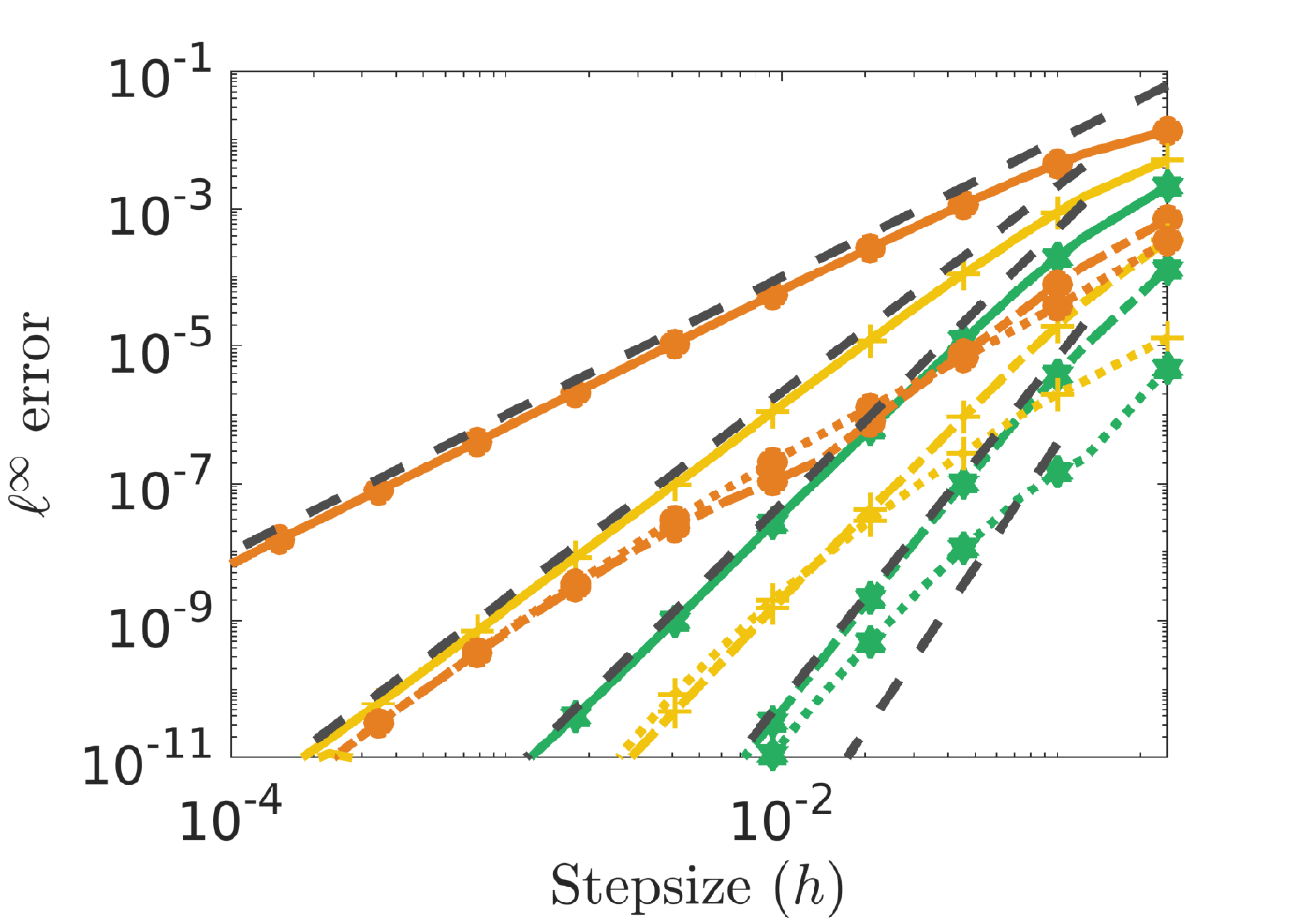}
		
		\begin{footnotesize}
			Convergence Diagram ($\epsilon = 10^{-5}$)
		\end{footnotesize}
		
		\includegraphics[width=1\linewidth]{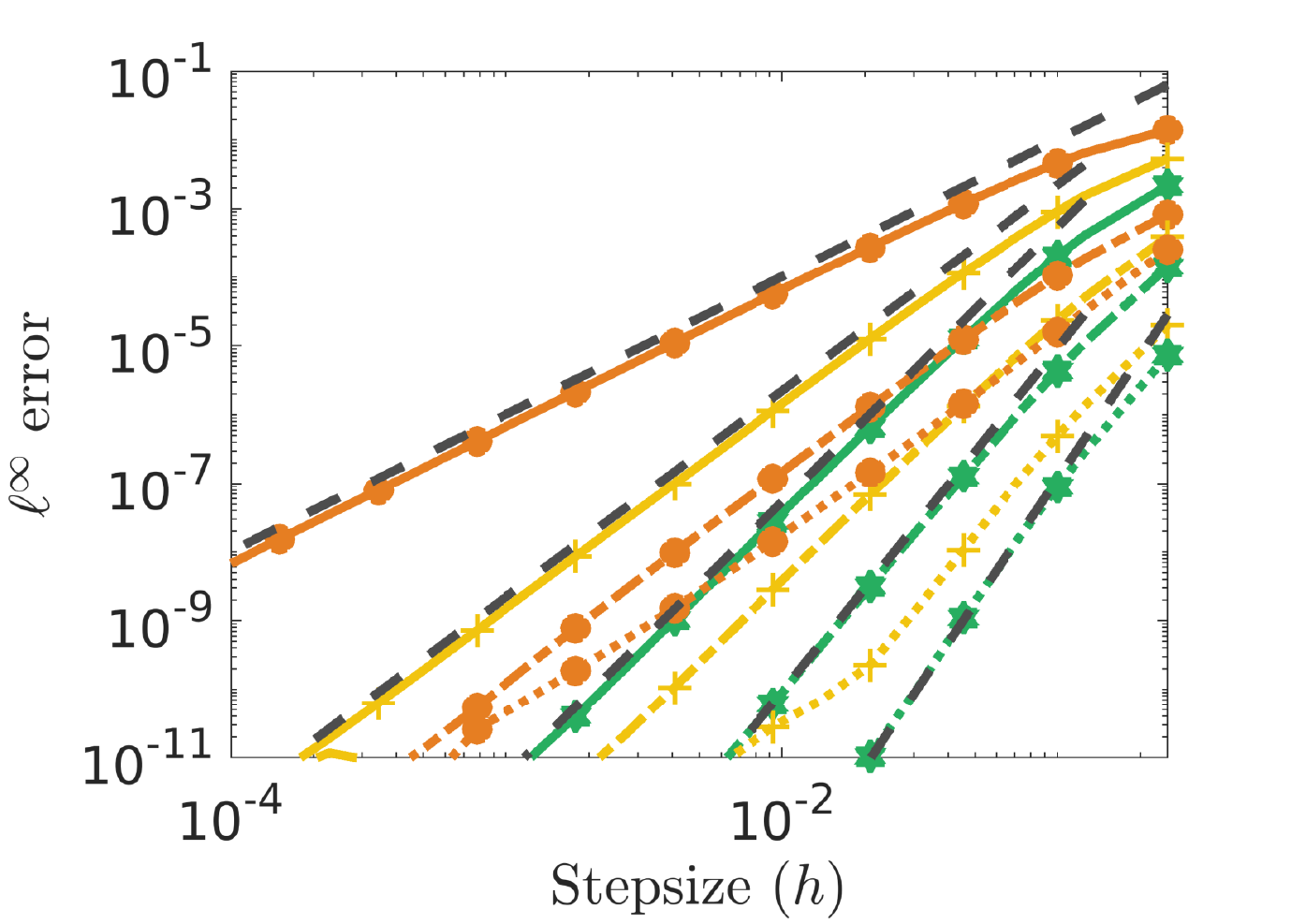}
		
		\begin{footnotesize}
			Convergence Diagram ($\epsilon = 10^{-8}$)
		\end{footnotesize}
		
		\includegraphics[width=1\linewidth]{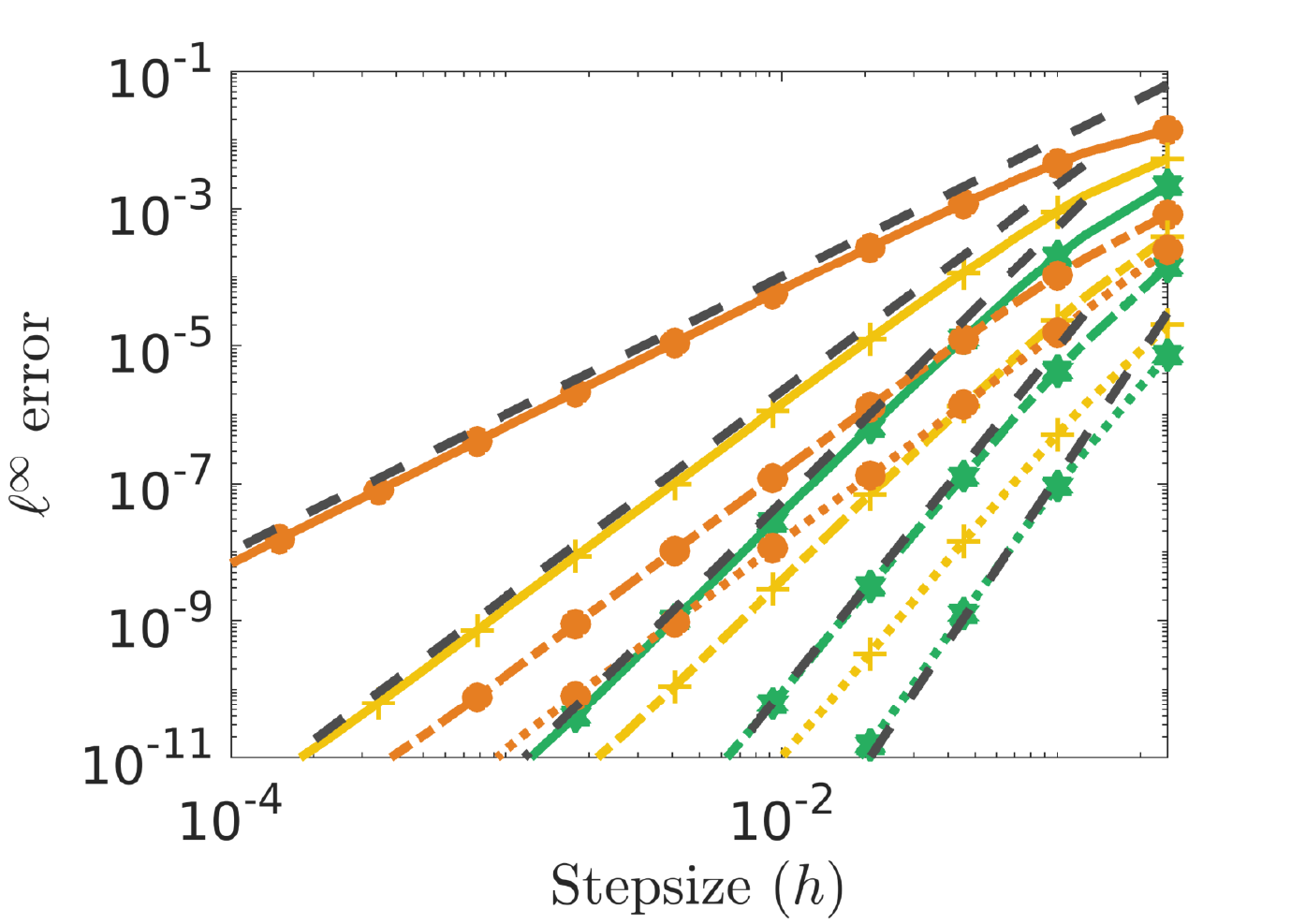}
	\end{minipage}
	\begin{minipage}{0.48\textwidth}
		\centering
		
		\begin{footnotesize}
			Convergence Diagram ($\epsilon = 10^{-3}$)
		\end{footnotesize}

		\includegraphics[width=1\linewidth]{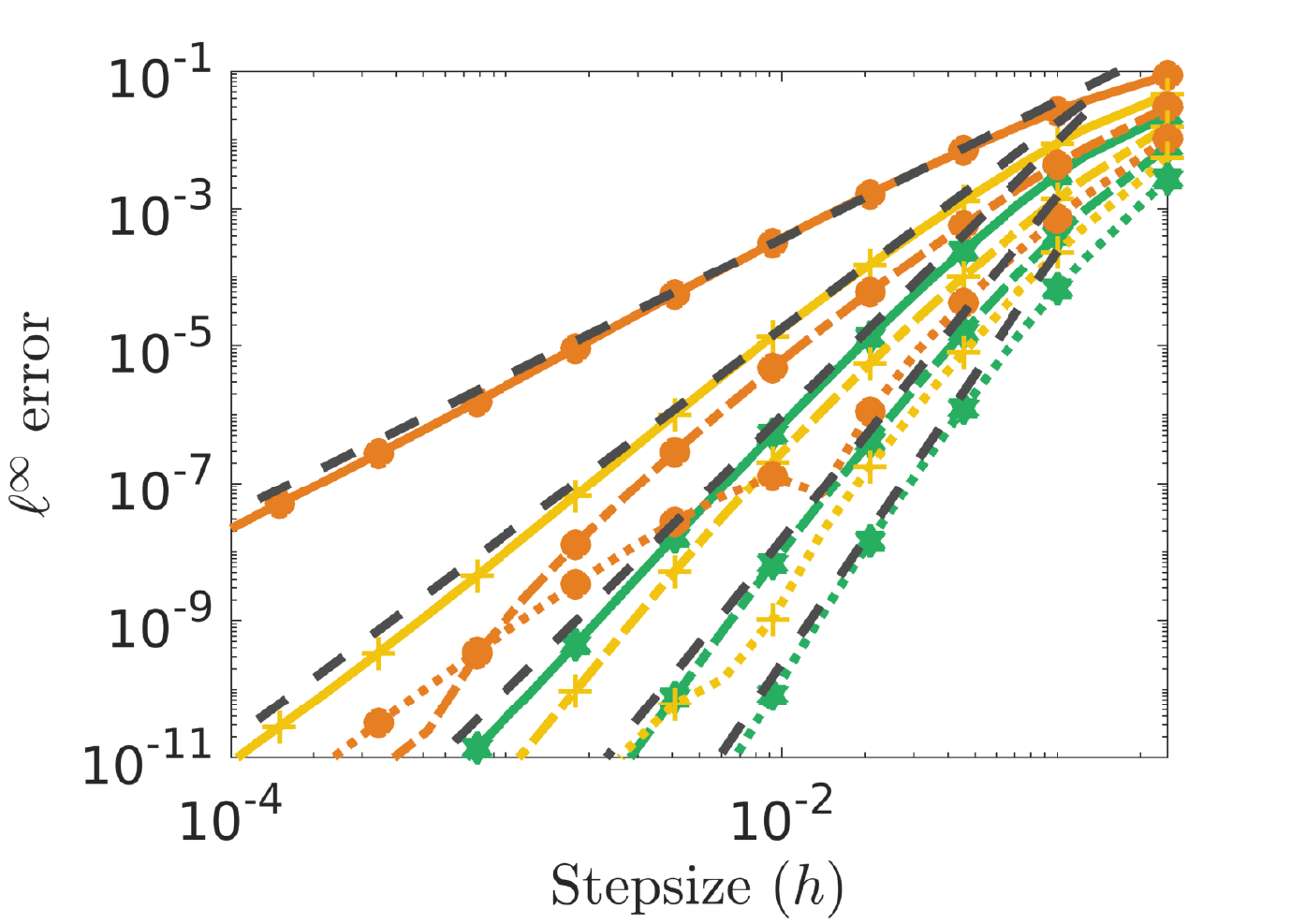}
		
		\begin{footnotesize}
			Convergence Diagram ($\epsilon = 10^{-5}$)
		\end{footnotesize}
		
		\includegraphics[width=1\linewidth]{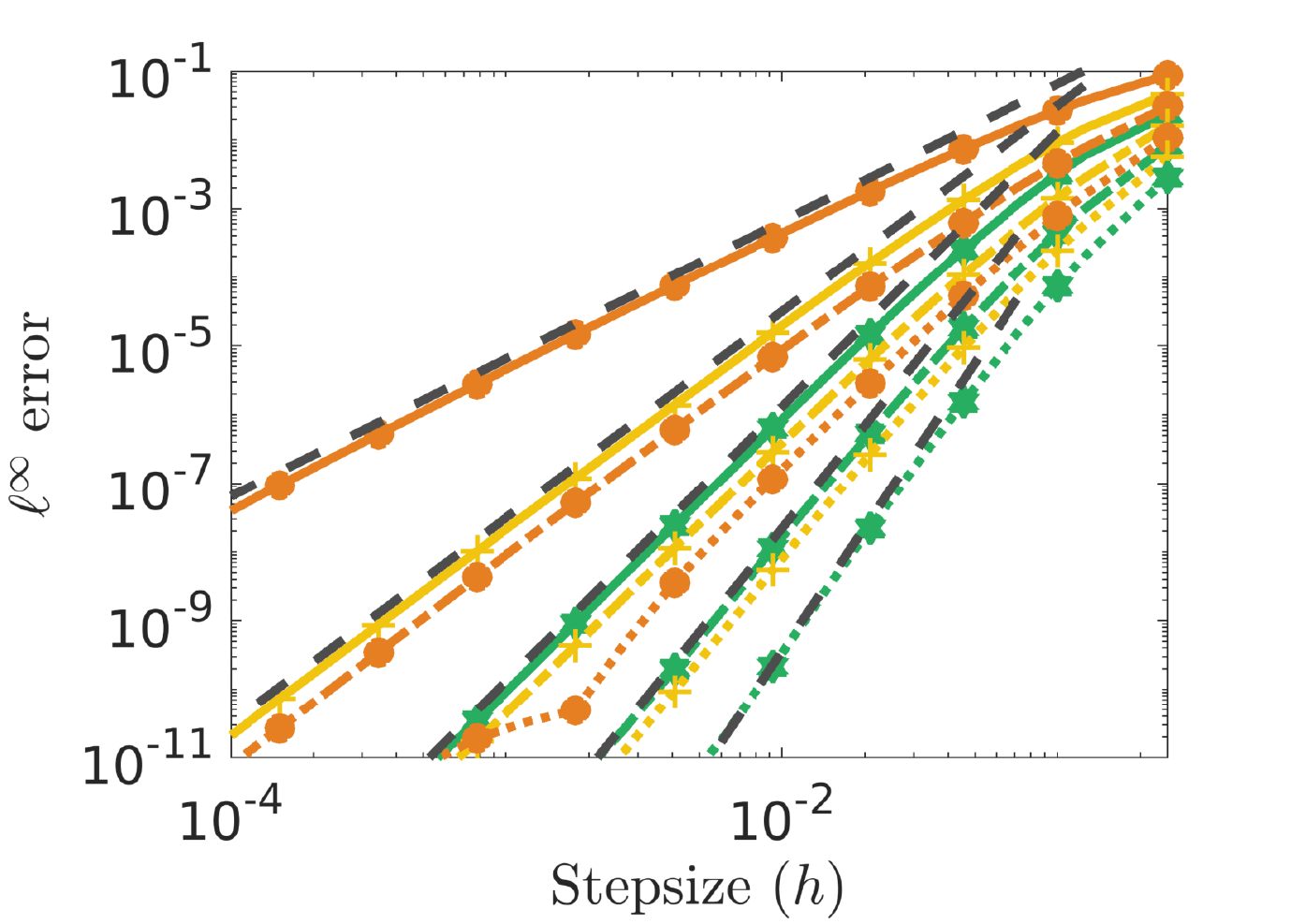}
		
		\begin{footnotesize}
			Convergence Diagram ($\epsilon = 10^{-8}$)
		\end{footnotesize}
		
		\includegraphics[width=1\linewidth]{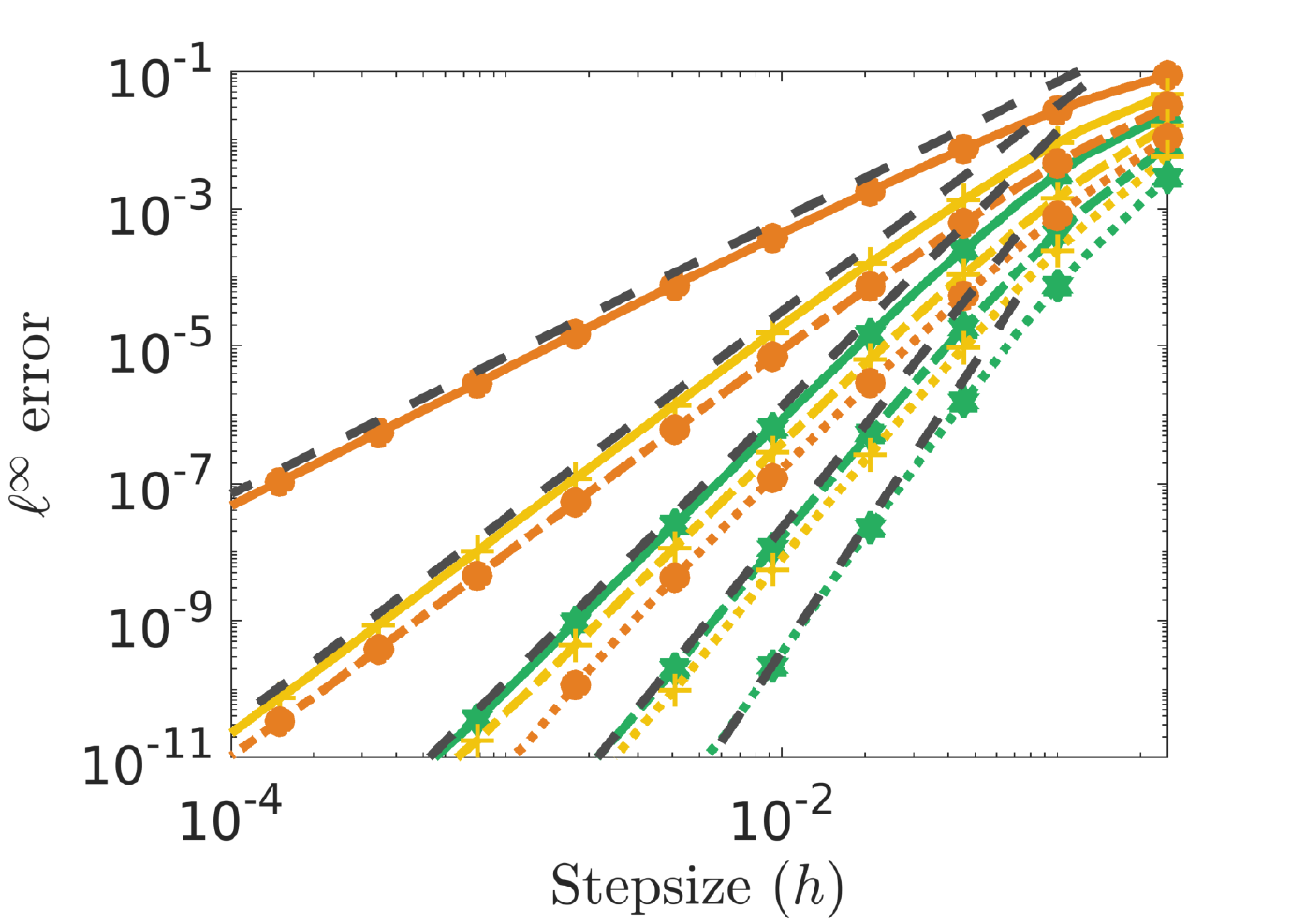}
	\end{minipage}
	\vspace{0.5em}
	
	\includegraphics[width=0.55\textwidth,trim={3.3cm 2.5cm 3.3cm 0},clip]{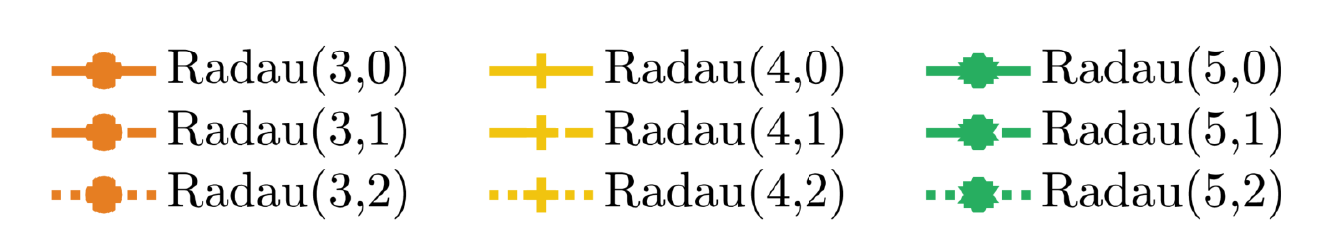}

	\caption{Estimated convergent rate (first row) and convergence diagrams (rows two through four) for FIMEX-Radau methods on the Van der Pol equation. The left and right columns respectively contain results for the semi-implicit and linearly implicit splittings.}
	\label{supfig:vanderpol-convergence-rate-imex-radau-extra}
\end{figure}

\begin{figure}[h!]
	\centering
	
	\begin{minipage}{0.48\textwidth}
		\centering
		\begin{footnotesize}
			{\bf Semi-Implicit Splitting}	
		\end{footnotesize}
		\vspace{1em}
		
		\begin{footnotesize}
			Convergence Diagram ($\epsilon = 10^{-3}$)
		\end{footnotesize}

		\includegraphics[width=1\linewidth]{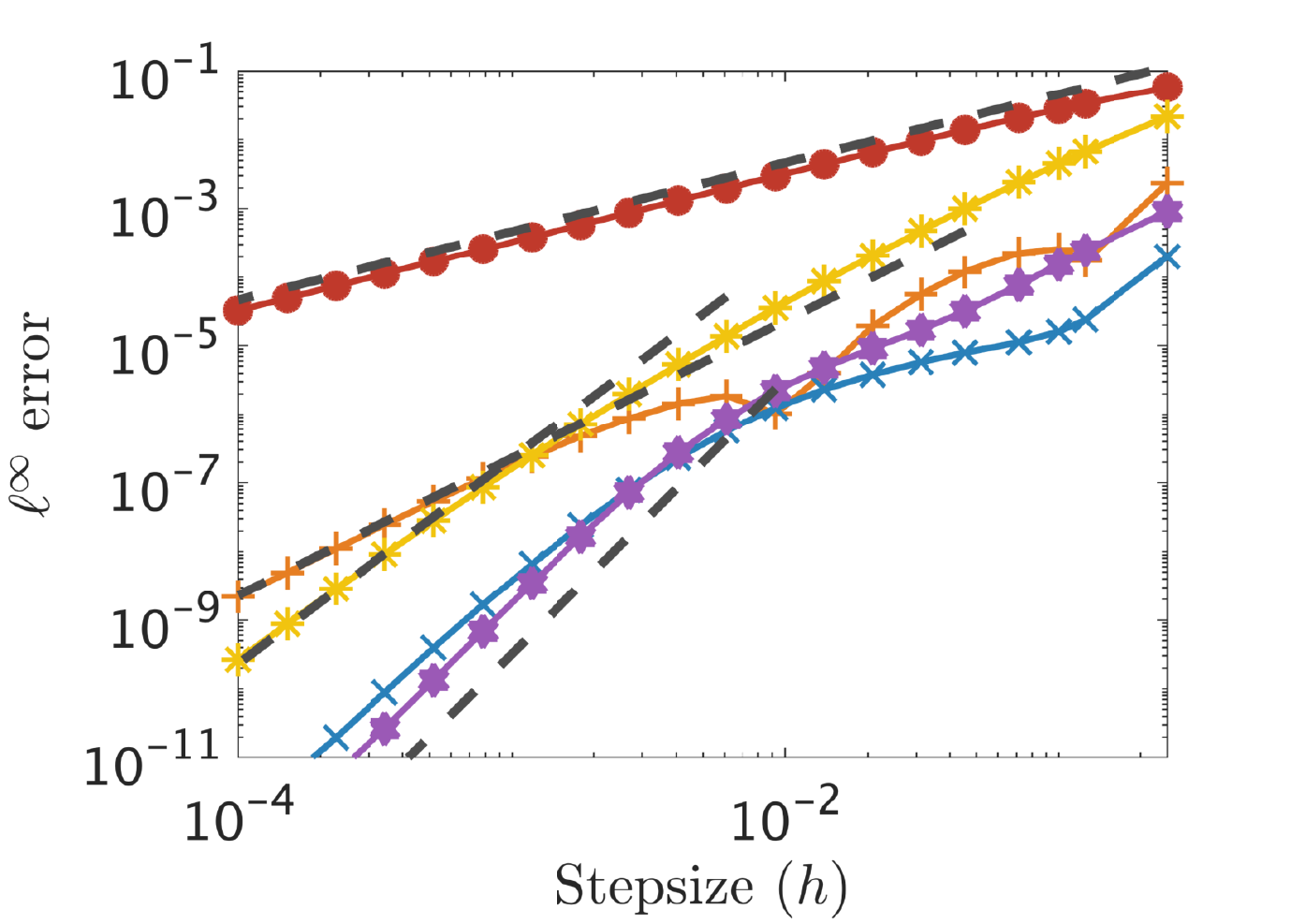}
		
		\begin{footnotesize}
			Convergence Diagram ($\epsilon = 10^{-5}$)
		\end{footnotesize}

		\includegraphics[width=1\linewidth]{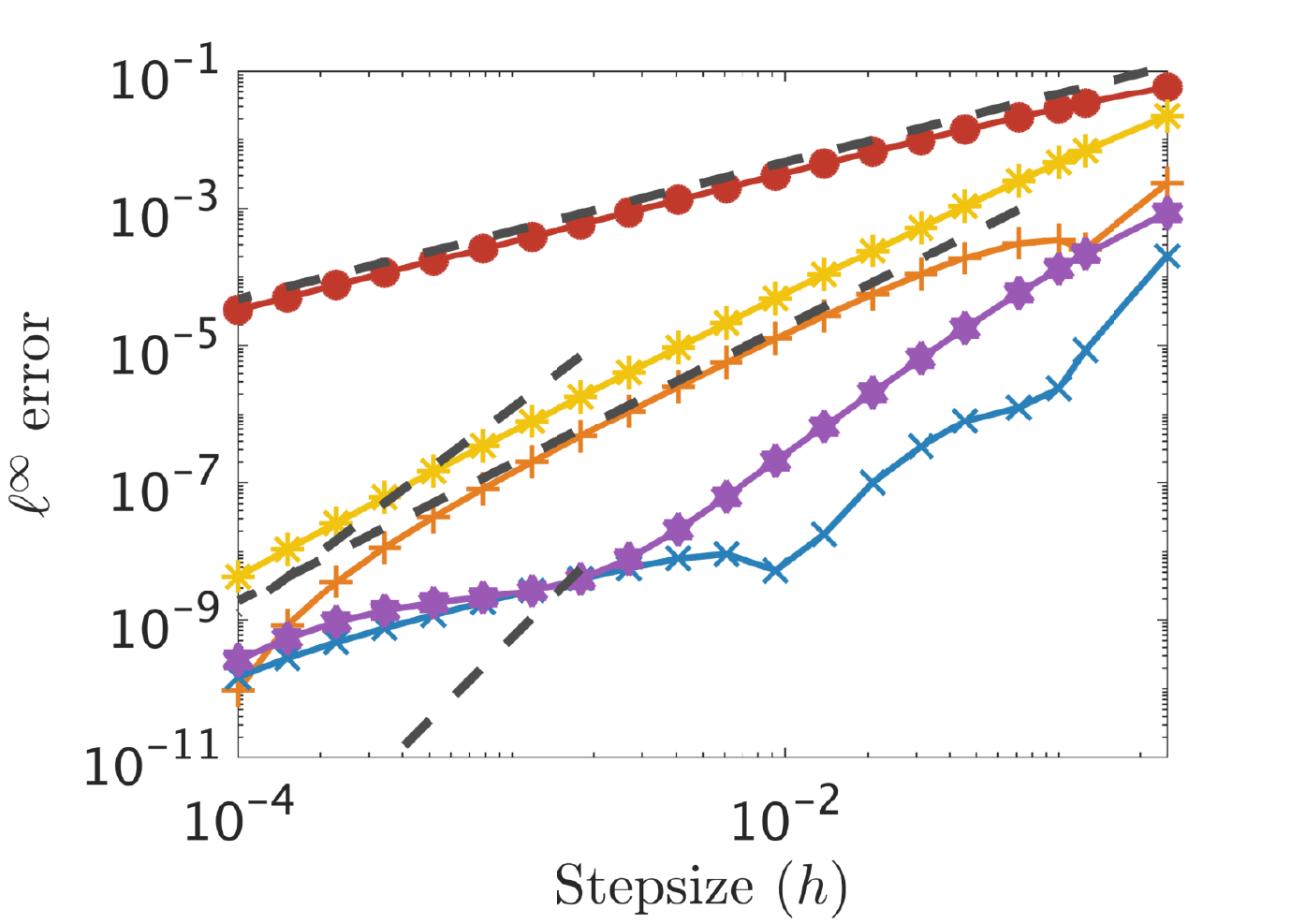}
		
		\begin{footnotesize}
			Convergence Diagram ($\epsilon = 10^{-8}$)
		\end{footnotesize}
		
		\includegraphics[width=1\linewidth]{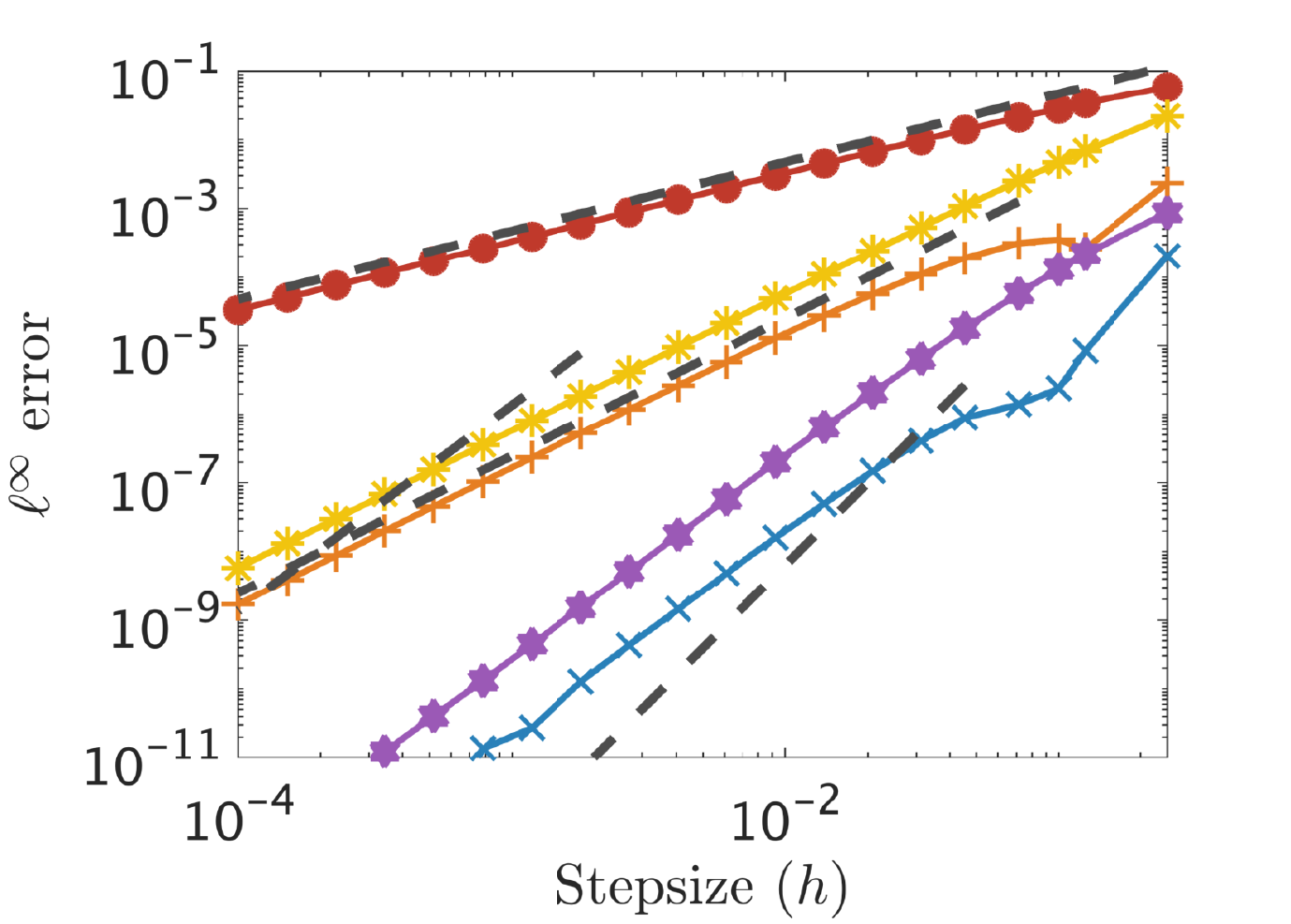}
	\end{minipage}
	\begin{minipage}{0.48\textwidth}
		\centering
		\begin{footnotesize}
			{\bf Linearly-Implicit Splitting}	
		\end{footnotesize}
		\vspace{1em}
		
		\begin{footnotesize}
			Convergence Diagram ($\epsilon = 10^{-3}$)
		\end{footnotesize}

		\includegraphics[width=1\linewidth]{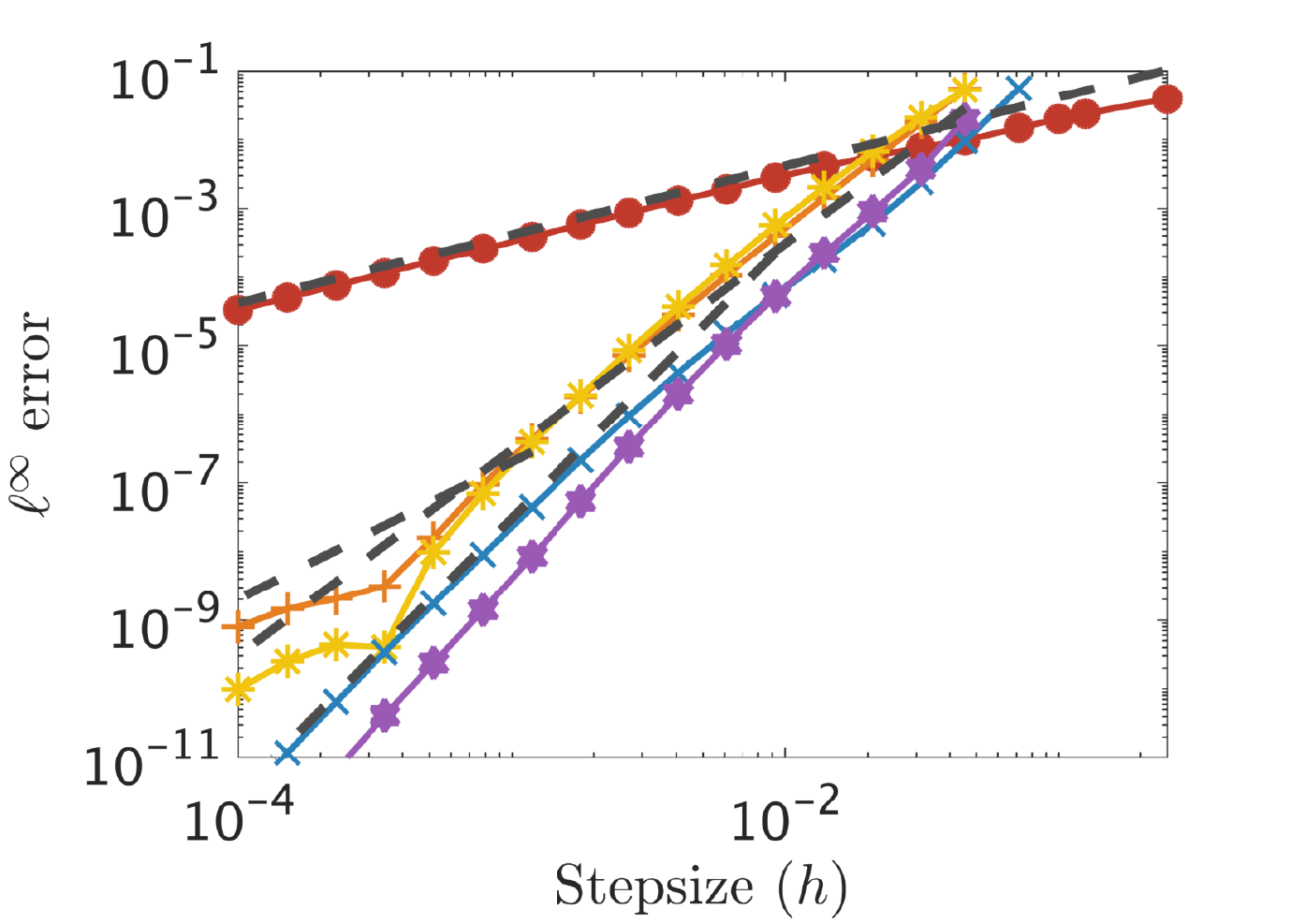}
		
		\begin{footnotesize}
			Convergence Diagram ($\epsilon = 10^{-5}$)
		\end{footnotesize}
		
		\includegraphics[width=1\linewidth]{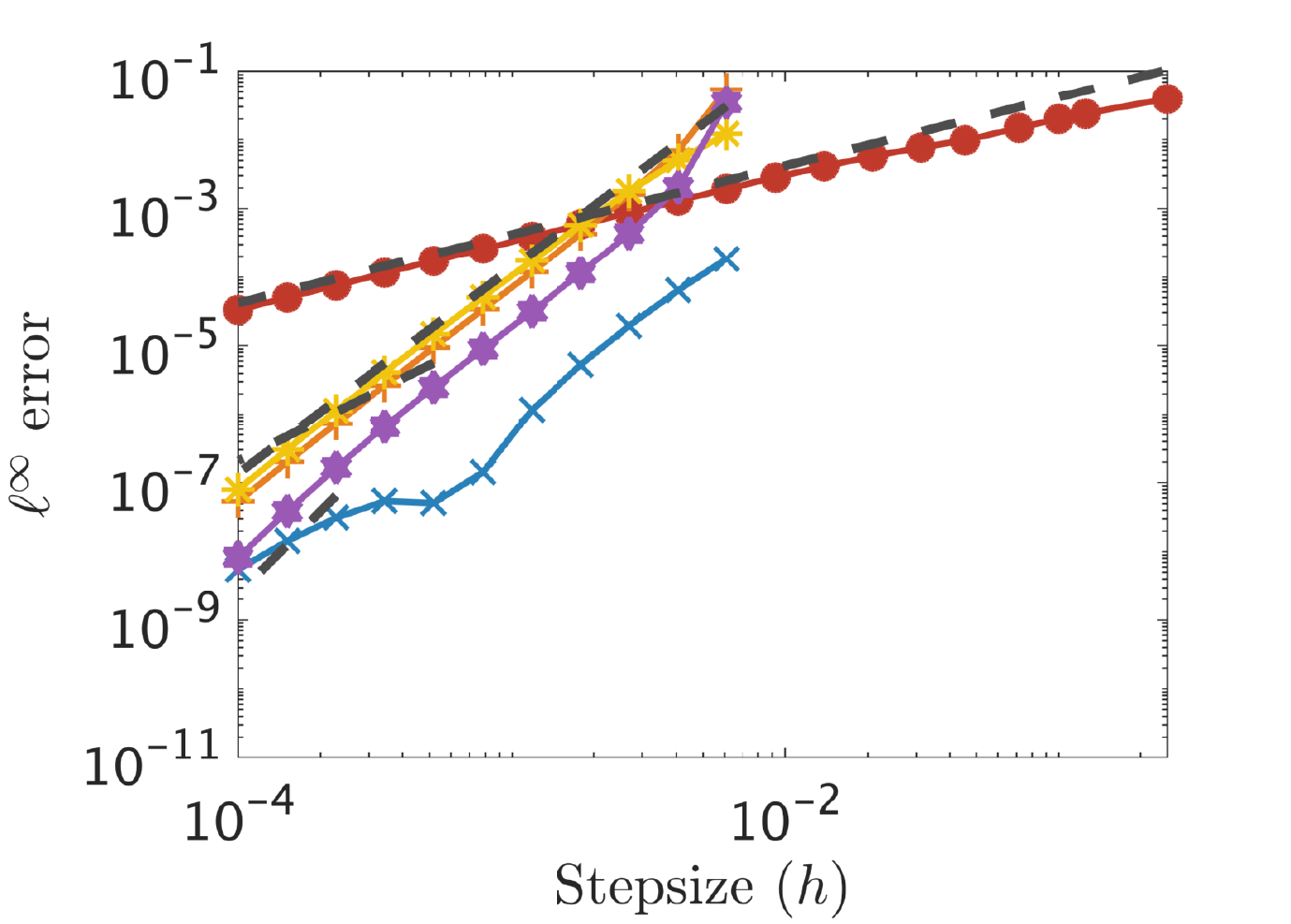}
		
		\begin{footnotesize}
			Convergence Diagram ($\epsilon = 10^{-8}$)
		\end{footnotesize}
		
		\includegraphics[width=1\linewidth]{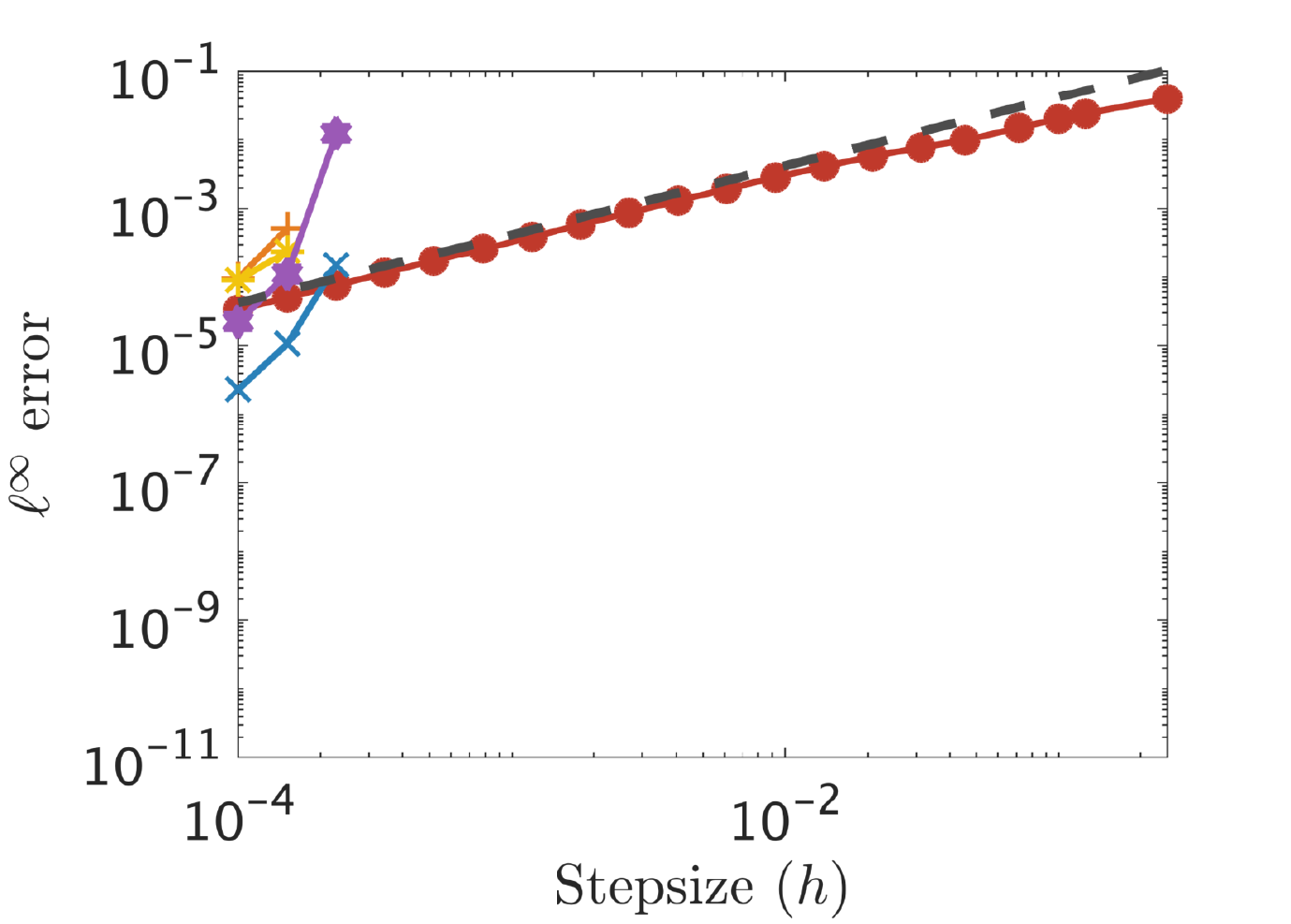}
		
	\end{minipage}
	\vspace{0.5em}
	
	\begin{minipage}{0.75\textwidth}
		\centering
		\includegraphics[width=1\linewidth,trim={.1cm 7.5cm 1.95cm 0},clip]{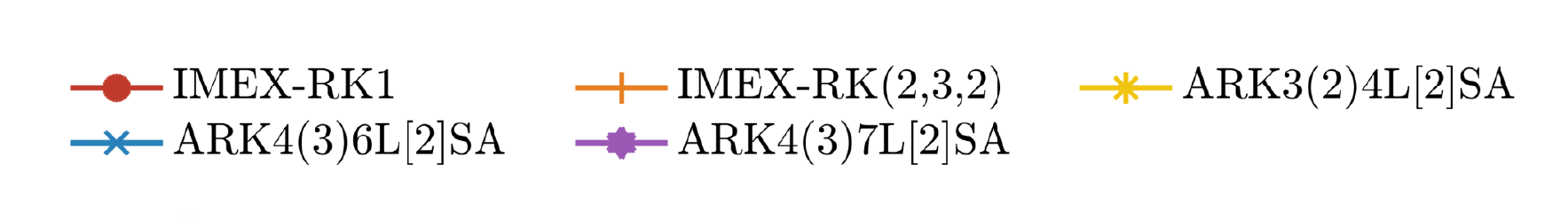}
	\end{minipage}
	\vspace{1em}

	\caption{Four example convergence diagrams for IMEX-RK methods on the Van der Pol equation using the semi-implicit and linearly implicit splittings. IMEX-RK methods require very small timesteps to remain stable when using the linearly-implicit splitting for small $\epsilon$.}
	\label{supfig:vanderpol-convergence-rate-imex-rk}
\end{figure}

This subsection contains additional convergence diagrams for FIMEX-Radau and IMEX-RK method on the Van der Pol equation. These diagrams supplement those shown in Figure \ref{fig:vanderpol-convergence-rate-radau}. In Figure \ref{supfig:vanderpol-convergence-rate-imex-radau-extra} we show four additional figures for FIMEX-Radau, and in Figure \ref{supfig:vanderpol-convergence-rate-imex-rk} we show the corresponding diagrams for IMEX-RK methods. The results clearly demonstrate the superior stability properties of the FIMEX-Radau methods for the linearly implicit splitting when $\epsilon$ is small.

\begin{figure}[h!]
	\centering
	
	\begin{minipage}{0.48\textwidth}
		\centering
		\begin{footnotesize}
			{\bf Semi-Implicit Splitting}	
		\end{footnotesize}
		\vspace{1em}

		\begin{footnotesize}
			Approximate Convergence Rate
		\end{footnotesize}
		\includegraphics[width=1\linewidth]{figures/experiments/vanderpol/SI-RadauS-conv}
	\end{minipage}
	\begin{minipage}{0.48\textwidth}
		\centering
		\begin{footnotesize}
			{\bf Linearly-Implicit Splitting}	
		\end{footnotesize}
		\vspace{1em}
		
		\begin{footnotesize}
			Approximate Convergence Rate
		\end{footnotesize}
		\includegraphics[width=1\linewidth]{figures/experiments/vanderpol/LI-RadauS-conv}
	\end{minipage}
	
	\begin{minipage}{0.48\textwidth}
		\centering
		
		\begin{footnotesize}
			Convergence Diagram ($\epsilon = 10^{-3}$)
		\end{footnotesize}

		\includegraphics[width=1\linewidth]{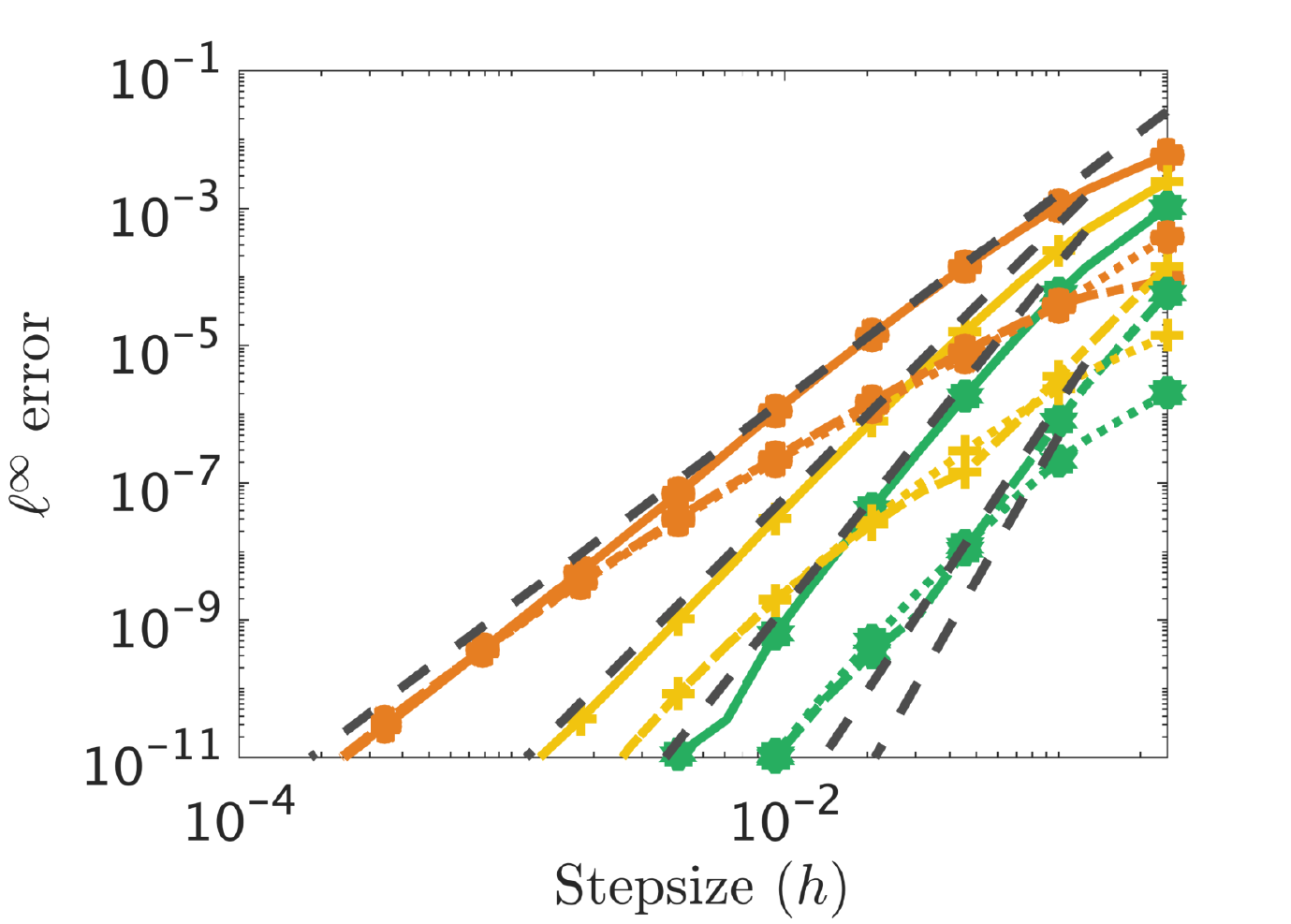}
		
		\begin{footnotesize}
			Convergence Diagram ($\epsilon = 10^{-5}$)
		\end{footnotesize}
		
		\includegraphics[width=1\linewidth]{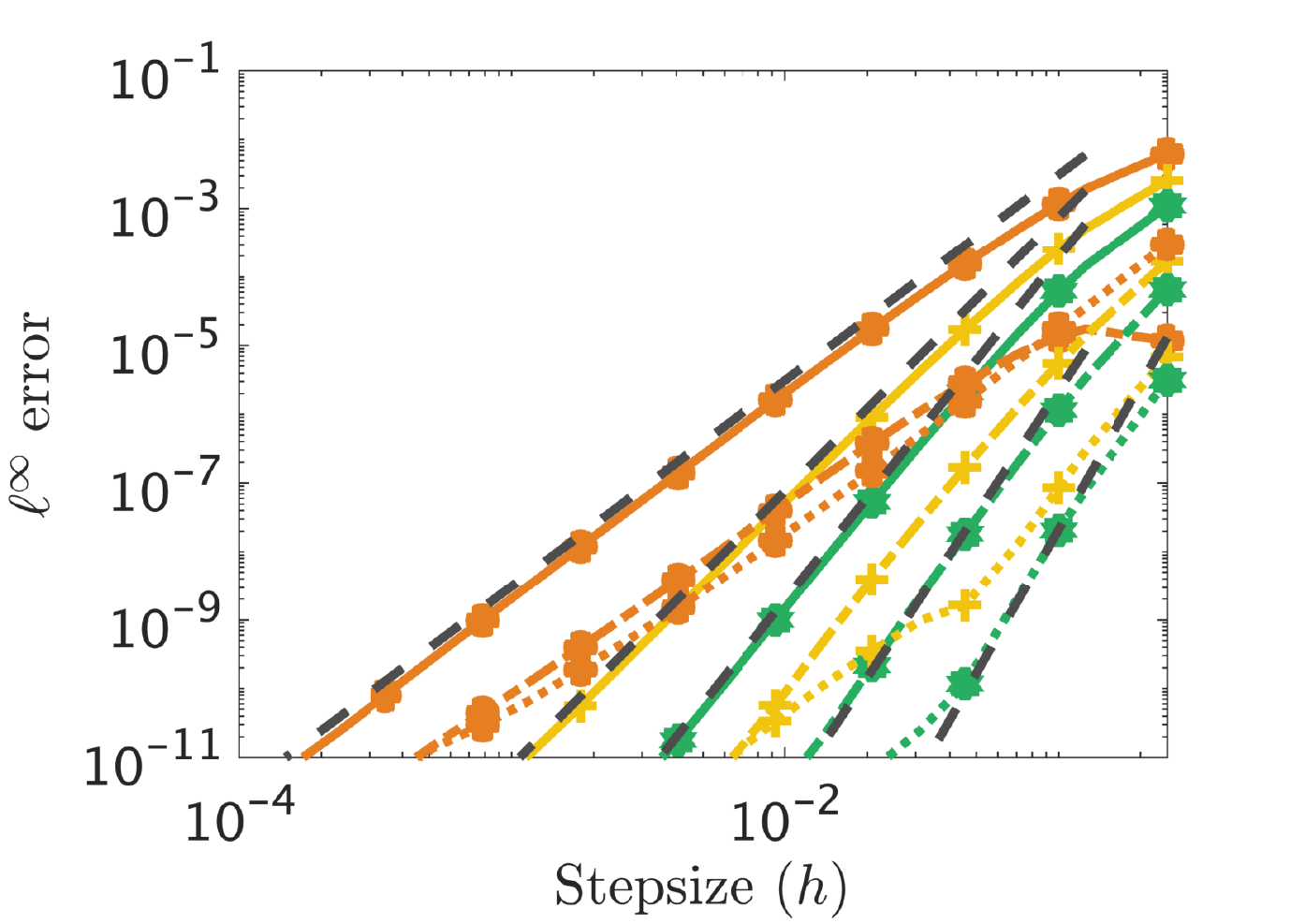}
		
		\begin{footnotesize}
			Convergence Diagram ($\epsilon = 10^{-8}$)
		\end{footnotesize}
		
		\includegraphics[width=1\linewidth]{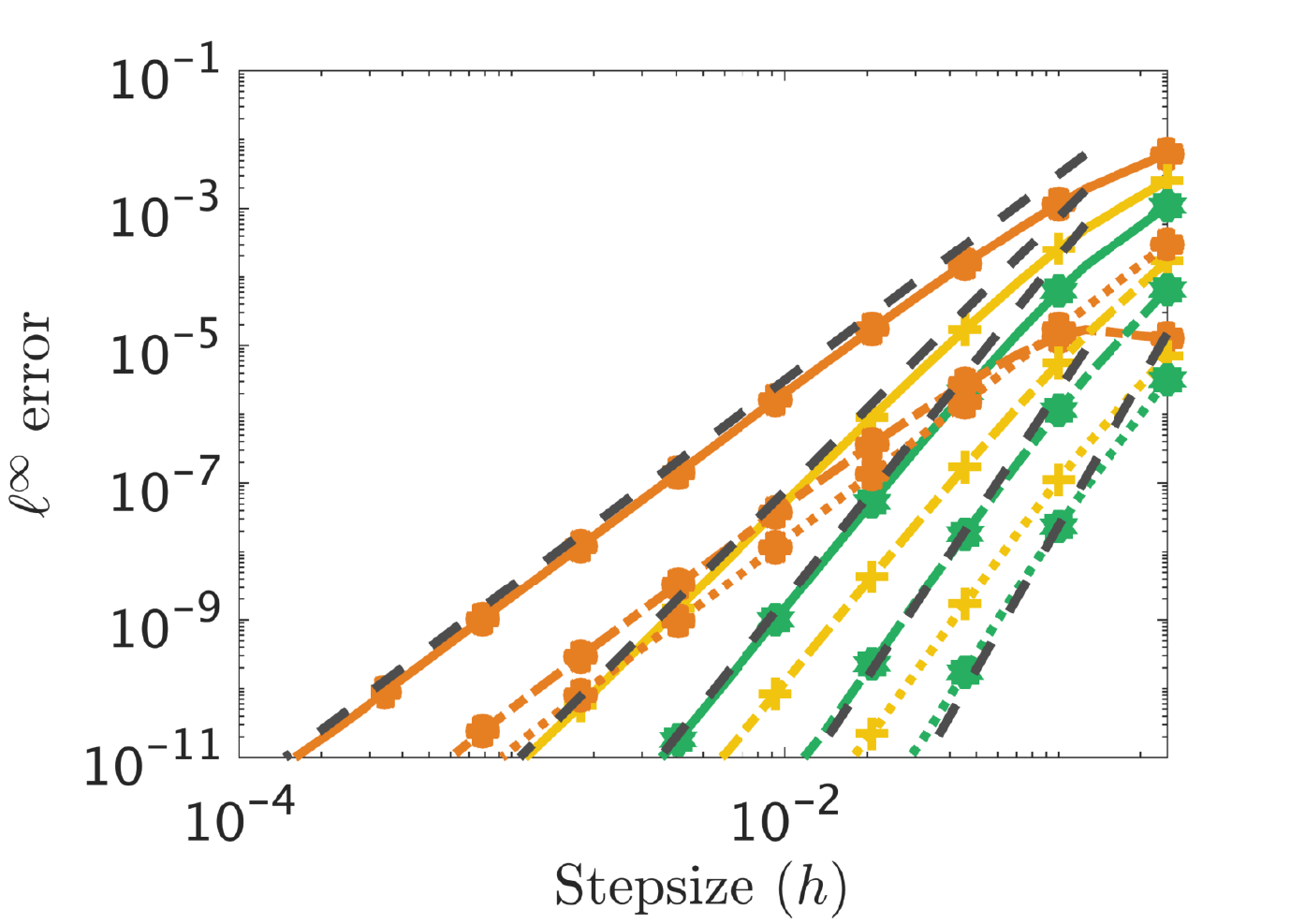}
	\end{minipage}
	\begin{minipage}{0.48\textwidth}
		\centering
		
		\begin{footnotesize}
			Convergence Diagram ($\epsilon = 10^{-3}$)
		\end{footnotesize}

		\includegraphics[width=1\linewidth]{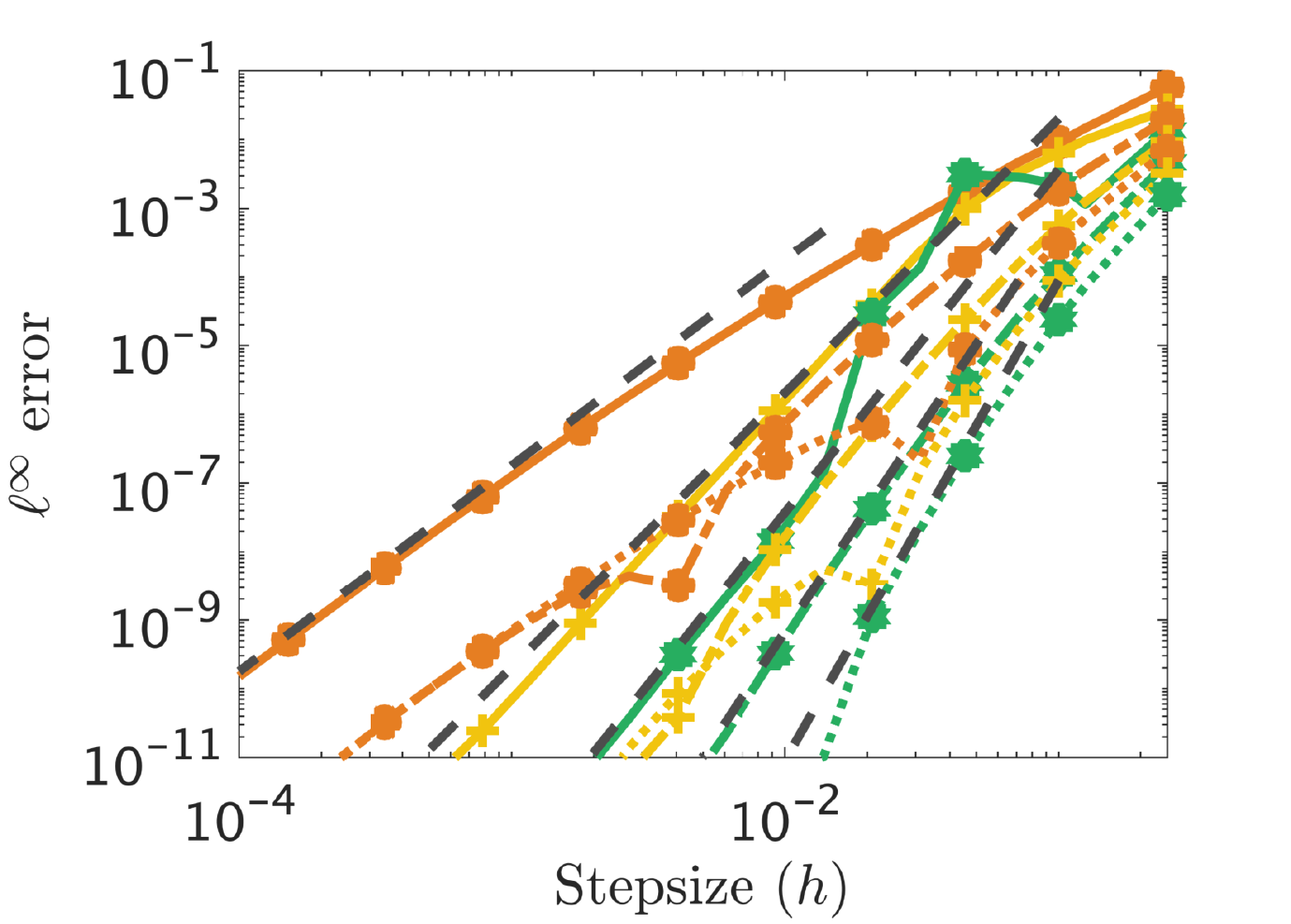}
		
		\begin{footnotesize}
			Convergence Diagram ($\epsilon = 10^{-5}$)
		\end{footnotesize}
		
		\includegraphics[width=1\linewidth]{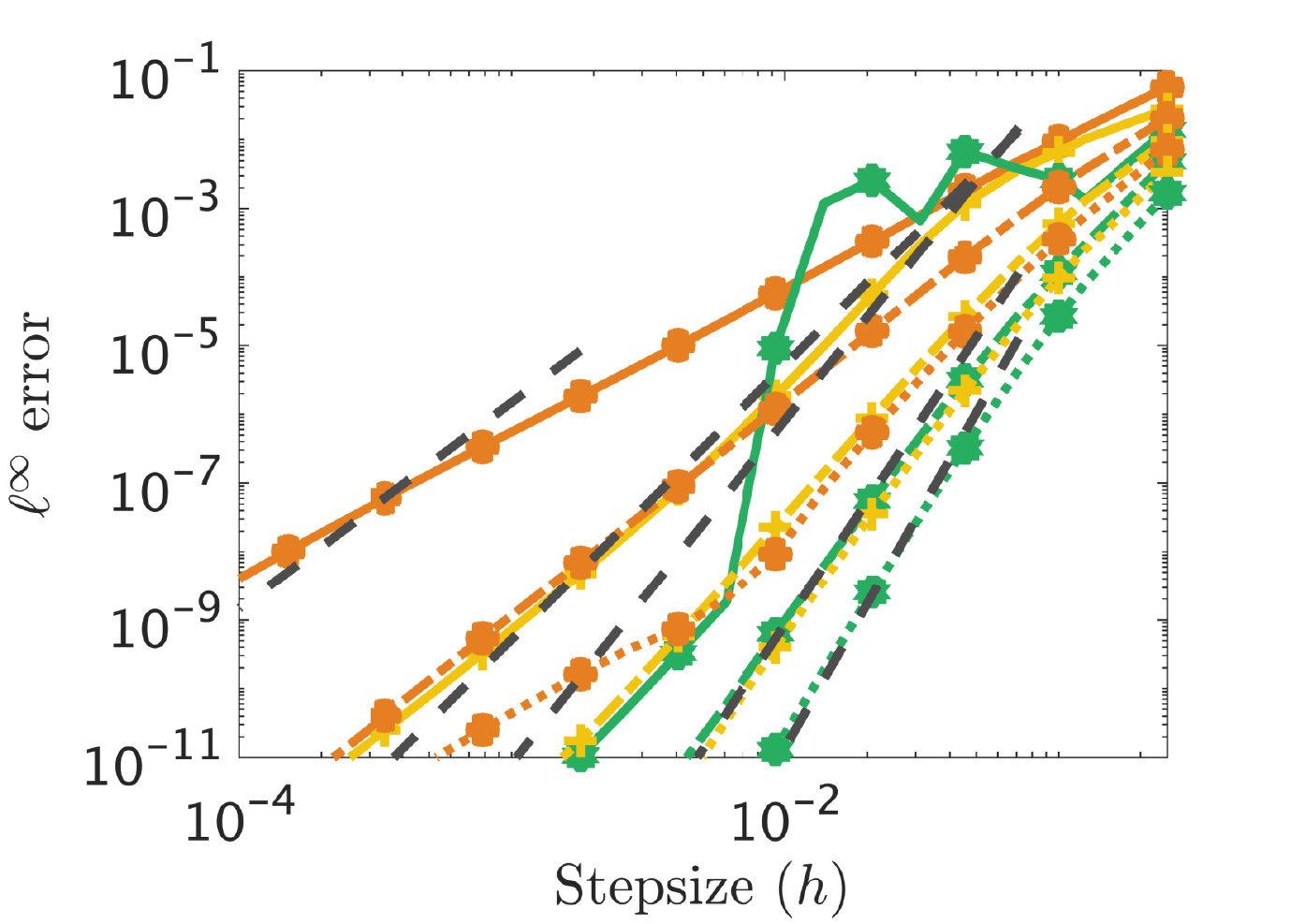}
		
		\begin{footnotesize}
			Convergence Diagram ($\epsilon = 10^{-8}$)
		\end{footnotesize}
		
		\includegraphics[width=1\linewidth]{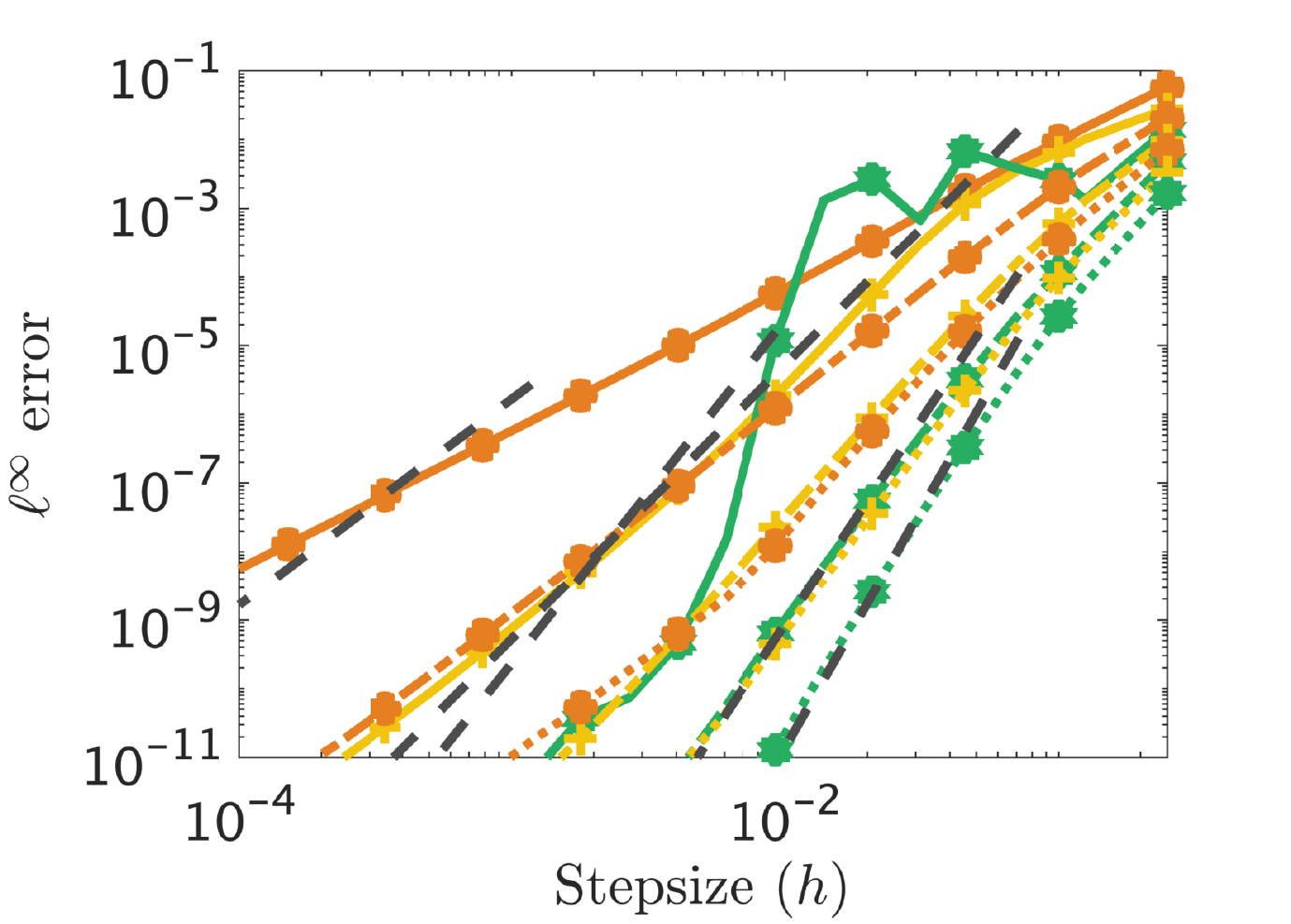}
	\end{minipage}
	\vspace{0.5em}
	
	\includegraphics[width=0.55\textwidth,trim={3.3cm 2.5cm 3.3cm 0},clip]{figures/experiments/vanderpol/SI-RADAUS-Kappa-eps-0_legend}

	\caption{Estimated convergent rate (first row) and convergence diagrams (rows two through four) for FIMEX-Radau methods on the Van der Pol equation. The left and right columns respectively contain results for the semi-implicit and linearly implicit splittings. The approximate convergent rate diagrams are identical to those shown in Figure \ref{fig:vanderpol-convergence-rate-radaus}.}
	\label{supfig:vanderpol-convergence-rate-imex-radaus-extra}
\end{figure}

\section{Polynomial IMEX PBMs as parametrized IMEX GLMS}
\label{sup:additive-polynomial-as-additive-glm}

Here we show the equivalence between an IMEX polynomial block method (i.e. \cref{eq:block_coefficient_partitioned} with $m=2$ and $\mathbf{C}^{\{2\}}$ strictly lower triangular) and a parametrized IMEX general linear method (GLM). Consider the the doubly partitioned equation
\begin{align}
	\mathbf{y}' = f^{\{1\}}(\mathbf{y}) + f^{\{2\}}(\mathbf{y}).
\end{align}
An IMEX general linear method that treats $f^{\{1\}}$ implicitly and $f^{\{2\}}$ explicitly is
\begin{align}
	Y_i &= \sum_{j=1}^q u_{ij} y_j^{[n]} + h \sum_{j=1}^{s} a^{\{1\}}_{ij} f^{\{1\}}(Y_j) + h \sum_{j=1}^{i-1} a^{\{2\}}_{ij} f^{\{2\}}(Y_j) & i&=1,\ldots, s\\
	y_i^{[n+1]} &= \sum_{j=1}^q v_{ij} y_j^{[n]} + h\sum_{j=1}^s \left( b^{\{1\}}_{j} f^{\{1\}}(Y_j) + b_{j} f^{\{2\}}(Y_j) \right) & i&=1,\ldots, q
\end{align}
where $s$ denotes the number of stages and $q$ is the number of inputs. Any IMEX method can be  represented compactly using the Butcher tablue
\begin{center}
	\begin{tabular}{l|l|l|l}	
		$c$ 	& $A^{\{1\}}$ & $A^{\{2\}}$ & $U$ \\ \hline
				& $B^{\{1\}}$ & $B^{\{2\}}$ & $V$
	\end{tabular}
\end{center}
where $A^{\{k\}} = [a^{\{k\}}_{ij}]$, $b^{\{k\}} = [b^{\{k\}}_{ij}]$, $U = [u_{ij}]$, $V = [v_{ij}]$.

All IMEX polynomial block methods can be recast as IMEX GLMs where the matrices are parametrized in terms of the extrapolation parameter $\alpha$. We will use bold face to denote the matrices of the polynomial block method and regular font for the GLM matrices. 

We start from an IMEX polynomial PBM with $\mathbf{q}$ inputs
\begin{align}
	\mathbf{y}^{[n+1]} = \mathbf{A}(\alpha) \mathbf{y}^{[n]} + r \sum_{k=1}^2 \mathbf{B}^{\{k\}}(\alpha) \mathbf{f}^{\{k\}[n]} + r \sum_{k=1}^2  \mathbf{C}^{\{k\}}(\alpha) \mathbf{f}^{\{k\}[n+1]}
	\label{eq:imex-pbm-app}
\end{align}
where we take $\mathbf{C}^{\{2\}}$ strictly lower triangular so that we are explicit in $f^{\{2\}}$. Next define a GLM with inputs
\begin{align}
	y^{[n]} = 
	\begin{bmatrix}
		\mathbf{y}^{[n]} \\
		\mathbf{f}^{\{1\}[n]} \\
		\mathbf{f}^{\{2\}[n]}
	\end{bmatrix}
\end{align}
and parametrized coefficient matrices
\begin{align*}
	U(\alpha) &= \left[ \mathbf{A}(\alpha) ~|~ \mathbf{B}^{\{1\}}(\alpha) ~|~ \mathbf{B}^{\{2\}}(\alpha) \right], &
	A^{\{1\}}(\alpha) &= \mathbf{C}^{\{1\}}(\alpha), &
	A^{\{2\}}(\alpha) &= \mathbf{C}^{\{2\}}(\alpha), \\[1em]
	V(\alpha) &= 
	\begin{bmatrix}
 		\mathbf{A}(\alpha) & \mathbf{B}^{\{1\}}(\alpha) & \mathbf{B}^{\{2\}}(\alpha)\\
 		\mathbf{0} & \mathbf{0} & \mathbf{0}\\
 		\mathbf{0} & \mathbf{0} & \mathbf{0}
 	\end{bmatrix}, &
 	B^{\{1\}}(\alpha) &= 
	\begin{bmatrix}
 		\mathbf{C}^{\{1\}}(\alpha) \\
 		\mathbf{I} \\
 		\mathbf{0}
 	\end{bmatrix}, &	
 	B^{\{2\}}(\alpha) &= 
	\begin{bmatrix}
 		\mathbf{C}^{\{2\}}(\alpha) \\
 		\mathbf{0} \\
 		\mathbf{I}
 	\end{bmatrix}.
\end{align*}
These choices lead to a GLM with $q=3\mathbf{q}$ inputs and $s=\mathbf{q}$ stages that is equivalent to and IMEX PBM \cref{eq:imex-pbm-app}.

\vfill
\begin{center}
	\includegraphics[width=1em]{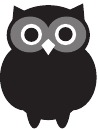}	
\end{center}

